\documentclass[12pt,leqno]{amsart} 
\theoremstyle{plain} 

\newtheorem{global-theorem}{Theorem}
\newtheorem{theorem}{Theorem}[section]
\newtheorem{lemma}[theorem]{Lemma}
\newtheorem{scholium}[theorem]{Scholium}

\newtheorem{corollary}[theorem]{Corollary}
\newtheorem{conjecture}[theorem]{Conjecture}

\newtheorem{proposition}[theorem]{Proposition}
\newtheorem{convention}[theorem]{Convention}
\newtheorem{problem}[theorem]{Problem}
\newtheorem{question}[theorem]{Question}

\newtheorem{exercise}{Exercise}

\newtheorem{prop-def}[theorem]{Proposition-Definition}
\newtheorem{lemma-def}[theorem]{Lemma-Definition}

\newcommand{\eop}{\ \hfill $\Box$}

\newcommand{\delbar}{\overline{\partial}}

\newcommand{\invisible}[1]{}

\numberwithin{equation}{section}

\newcommand{\cc}{{\mathbb C}}
\newcommand{\pp}{{\mathbb P}}
\newcommand{\rr}{{\mathbb R}}

\newcommand{\zz}{{\mathbb Z}}

\newcommand{\hh}{{\mathbb H}}
\newcommand{\ggg}{{\mathbb G}}

\newcommand{\aaa}{{\mathbb A}}
\newcommand{\LL}{{\mathbb L}}

\newcommand{\Cc}{{\mathcal C}}

\newcommand{\Dd}{{\mathcal D}}
\newcommand{\Ff}{{\mathcal F}}
\newcommand{\Oo}{{\mathcal O}}

\newcommand{\Mm}{{\mathcal M}}

\newcommand{\overvect}[1]{{\stackrel{\rightharpoonup}{#1}}}

\begin{document}

\author[C. Simpson]{Carlos Simpson}
\address{CNRS, Laboratoire J. A. Dieudonn\'e, UMR 6621
\\ Universit\'e de Nice-Sophia Antipolis\\
06108 Nice, Cedex 2, France}
\email{carlos@math.unice.fr}
\urladdr{http://math.unice.fr/$\sim$carlos/} 

\title[Middle convolution]{Katz's middle convolution algorithm}

\begin{abstract}
This is an expository account of Katz's middle convolution operation on local systems
over ${\bf P}^1-\{ q_1,\ldots , q_n\}$. We describe the Betti and de Rham versions,
and point out that they give isomorphisms between different moduli spaces of local systems,
following V\"olklein, Dettweiler-Reiter, Haraoka-Yokoyama. 
Kostov's program for applying the Katz algorithm is to say that in the range where middle convolution 
no longer reduces the rank, one should give a direct construction of local systems. This has been 
done by Kostov and Crawley-Boevey. We describe here an alternative construction using the notion of 
cyclotomic harmonic bundles: these are like variations of Hodge structure except that the Hodge 
decomposition can go around in a circle. 
\end{abstract}


\maketitle


\section{Introduction} \label{introduction}

There is a growing body of literature about Katz's ``middle convolution'' algorithm on local systems on
$\pp ^1- \{ q_1, \ldots , q_n\}$. The purpose of the present paper is expository: we would like to
describe two versions of Katz's construction in complex geometry, the Betti version involving
complex local systems, and the de Rham version involving vector bundles with logarithmic connection. 
Katz's book \cite{Katz}
was written in the framework of $\ell$-adic sheaves, which at first made it difficult to understand for
complex geometers including myself. Subsequently, V\"{o}lklein and Dettweiler-Reiter recast the construction
in complex geometry and algebra. In \S \S 2.7-2.9 of Katz's book, the convolution was defined in a geometric way which
is applicable in any context where one has a Grothendieck formalism and a category of perverse sheaves. Thus, the translation
into complex geometry may be viewed as coming directly from there \cite[5.9]{Katz}. 
Katz then interpreted the convolution as conjugate to a
tensor product, via Fourier transform, and used that to obtain some of the main properties of his construction. 
The complex analogy for this would {\em a priori} bring into play
the notion of irregular connections on a $2$-dimensional variety, a theory which remains poorly understood
(see however \cite{Boalch2} and the recent preprint \cite{AkerSzabo}).
It is possible to do a full treatment of middle convolution staying within the realm of complex geometry but without using
Fourier transform, as has
been shown and exploited by the works of Strambach, V\"{o}lklein, Dettweiler, Reiter, Kostov, Crawley-Boevey, Haraoka, Yokoyama.

Many applications of Katz's theory concern the case of rigid local systems. 
For example, Gleizer has studied explicit solutions \cite{Gleizer},
and Roberts' preprint \cite{Roberts} includes an extensive discussion of how to apply the algorithm
to determine which rigid local systems exist. 
V\"olklein, Dettweiler and Reiter have done extensive work on using Katz's existence results in the rigid case to
construct motivic local systems with interesting monodromy groups, obtaining results on the inverse
Galois problem. 

The middle convolution 
transformation was first applied in the non-rigid case by Kostov. An important invariant which we denote
by $\delta (\overvect{g})$ is the change in rank induced by Katz's transformation. As long as 
$\delta (\overvect{g})< 0$ we can apply middle convolution to reduce the rank (or otherwise, conclude
that the local system couldn't exist). Kostov made the fundamental
observation that when we get into the range $\delta (\overvect{g}) \geq 0$, we should expect that the local system
always exists and look for a direct construction. Kostov applied this to solve the existence problem
in many cases \cite{KostovCRAS}--\cite{KostovConnectedness}, 
such as when $\overvect{g}$ is {\em simple} i.e. the multiplicities of eigenvalues are not
all divisible by the same integer $d>1$, or for generic eigenvalues even if $\overvect{g}$ is
not simple \cite{KostovControl}. Crawley-Boevey looked at the existence question from a point of
view of root systems in \cite{CrawleyBoeveyIHES} where Katz's algorithm plays a role.
In that language, the transformation on local monodromy data is considered as a reflection in a root system,
and a sequence of reflections is used to move up to the positive chamber. Once we are in the positive chamber,
analogous to the condition $\delta (\overvect{g})\geq 0$, Crawley-Boevey found a direct construction of indecomposable
parabolic bundles, and applied a parabolic variant of Weil's theorem to construct flat connections. 

At the end of the present paper, we propose a technique for constructing local systems in the range 
$\delta (\overvect{g})\geq 0$ by using the correspondence between Higgs bundles and local systems in the parabolic case.
This construction is heavily inspired by Kostov's program, and 
is obviously a variant on Crawley-Boevey's indecomposable parabolic bundles. So, it is not really very new but
might present some advantages such as making
clear the role played by the condition $\delta (\overvect{g})\geq 0$. The objects we introduce,
cyclotomic harmonic bundles, might
be interesting in their own right such as for studying the behavior of everything near infinity
in the moduli spaces.

We will 
look at Katz's operation as giving an isomorphism between moduli spaces for different local
monodromy data. Of course this applies to the rigid case too, but as Katz pointed out long ago, for 
local systems on $\pp ^1$ with specified singularities, local rigidity implies the stronger rigidity statement
that there is at most one irreducible representation with the given local data. Thus, in the rigid case the
moduli spaces are single points so even the cardinality is not an interesting invariant.  
Instead, we are motivated by looking for low-dimensional moduli spaces, for which things like
Hitchin's hyperk\"ahler structure, or the Riemann-Hilbert correspondence, could be viewed explicitly.
The phrase ``toy example'' was coined by T. Hausel in \cite{Hausel} to refer to this kind of low-dimensional
case arising from a punctured projective line. He looked at a space of parabolic Higgs bundles of rank $2$ on
$\pp ^1 - \{ q_1,\ldots , q_4\}$. Boalch looked at an example of the middle convolution relating this space to a
space of rank $3$ representations in \cite{Boalch1} \cite{Boalch2}, and considered the Painlev\'e equations for
these cases. It seems like a good idea to pursue the philosophy of looking at low-dimensional cases,
and to get started we need to have a thorough understanding of how the classification based on Katz's 
algorithm works. That's the motivation for this paper.

Conceptually, the middle convolution operation is pretty easy to understand. Let $Y$ and $Z$ denote two copies of
the projective line $\pp ^1$, with reduced effective
divisors $Q_Y\subset Y$ and $Q_Z\subset Z$ (which we will often denote just by
$Q$), both given by the same finite collection of $n$ points $Q = \{ q_1, \ldots , q_n\}$. Let 
$D\subset Z\times Y$ denote the ``diagonal configuration'' consisting of the diagonal $\Delta$
plus the vertical and horizontal divisors given by preprojections of $Q$. 

A {\em convoluter} is a rank one local system on $Z\times Y$ with singularities along $D$. 
In Katz's original setup this would be a rank one $\ell$-adic sheaf. In the complex geometric
``Betti'' and ``de Rham'' situations we consider here, the convolution object is respectively a
rank one representation of $\pi _1(Z\times Y -D)$, or a logarithmic connection on the trivial bundle given
by a logarithmic one-form. In either case we denote the convoluter by $\beta$. 
Let 
$$
\xi : Z\times Y -D \longrightarrow Z - Q_Z,\;\;
\eta : Z\times Y -D \longrightarrow Y - Q_Y
$$
denote the two projections. 
Given an irreducible rank $r$ local system $L$ on $Y-Q$, we can form the ``raw convolution'' defined as the higher
direct image
$$
RC_{\beta}(L) :=
\rr ^1\xi _{\ast}(\eta ^{\ast}(L)\otimes \beta ),
$$
a local system on $Z-Q_Z$. 
Unfortunately, the raw convolution 
will not in general be an irreducible local system, because there are some contributions whenever
the tensor product $\eta ^{\ast}(L)\otimes \beta $ has trivial eigenvalues along the ``horizontal'' piece
$H:= \eta ^{\ast}(Q_Y)$ of the divisor $D$. This is remedied by defining the ``middle convolution''  
to be the {\em middle direct image}
$$
MC_{\beta}(L) := M\rr ^1\xi _{\ast}(\eta ^{\ast}(L)\otimes \beta ),
$$
heuristically defined as the kernel of the map to the quotient 
systems corresponding to the unwanted local cohomolgy groups. 

When studying local systems on $\pp ^1$ with singularities, we are interested in fixing the local type of the
singularities. For the present paper, we will simplify things considerably by making the 
convention that the local monodromy transformations be semisimple (Convention \ref{semisimple}) or the corresponding
statement for the residues of a logarithmic connection (Convention \ref{semisimpleresidues}). 
This allows us to avoid complicated discussions of Jordan normal forms. The reader who is interested may
refer to the original references for discussions of this aspect.

When we discuss moduli spaces we will consider the moduli spaces of local systems with fixed
conjugacy classes of local monodromy. The notion of rigidity considered by Katz takes into
account the fixing of the local conjugacy classes. So, the first and in some sense main question about 
Katz's construction is
to understand what is its effect on the local monodromy transformations. 

We will try to explain the answer,
and how to see why it works that way. This will occupy most of the paper, 
and is the main subject of our exposition.  Of course it has already been treated by Katz in the $\ell$-adic case,
and by Strambach, V\"{o}lklein, Dettweiler, Reiter, Kostov, Crawley-Boevey, Haraoka, Yokoyama and others 
in the complex case. Thus there is nothing new
in our exposition. We hope it will be useful as an explanation allowing readers more easily
to consult the original references.

At the end of the paper, we consider the question of how to construct local systems in 
the range $\delta (\overvect{g}) \geq 0$. We propose a construction which is based on Donaldson-style
Yang-Mills theory \cite{Donaldson} \cite{Hitchin}, 
in which we construct a polystable Higgs bundle with parabolic structure
corresponding to the local monodromy. 
This is similar to the construction of systems of Hodge bundles 
which was used in \cite{products}. Here we introduce a new notion which
makes the problem much easier: {\em cyclotomic harmonic bundles}. These are harmonic bundles which are
fixed under the action of a finite cyclic subgroup $\mu _m \subset \cc ^{\ast}$, for the usual action of $\cc ^{\ast}$
on the space of Higgs bundles. These are related to the Higgs bundles considered
by Hitchin in \cite{HitchinTeichmuller}. A cyclotomic harmonic bundle is very much like a variation of Hodge structure,
in that the bundle decomposes $E= \bigoplus _{p=0}^{m-1} E^p$. The only difference is that the indexation
is really by $p\in \mu _m ^{\ast} \cong \zz / m\zz$, and the Kodaira-Spencer components of the Higgs field 
go between $E^p$ and $E^{p-1}\otimes \Omega ^1_X$ where $p-1$ is taken modulo $m$. Thus $\theta$ is no longer
necessarily nilpotent. Our construction takes place in the maximal case when $m=r$ is the rank of $E$,
and $\theta$ is not nilpotent. This means that the $E^p$ are line bundles. Thus the description of $(E,\theta )$
is elementary. It turns out that incorporating parabolic structures into the picture in order to obtain a required
local monodromy type, the condition $\delta (\overvect{g})\geq 0$ is exactly what is needed for the degrees
of the line bundles to work out correctly and enable the construction. Unfortunately it doesn't work
when the dimension of the moduli space is $2$.

In the last section we 
discuss some questions and directions for further study. 

{\em Acknowledgements:} I would specially like to thank V. Kostov for many helpful discussions, and for a particularly
illuminating talk many years ago in which he explained his utilisation of Katz's algorithm. Also I would like to thank
O. Gleizer for some interesting discussions a while ago. At Princeton last year,
Deligne and Katz raised the question of how to understand what is going on, which prompted the present write-up. 
Many aspects we consider here, such as moduli spaces with fixed conjugacy classes, and logarithmic connections,
showed up in the course of recent joint works with K. Corlette and J. Iyer.

\section{Connections and local systems}
\label{connectionslocalsystems}

Denote by $Y$ a 
smooth projective curve, with $K\subset Y$ a reduced divisor. Write $K = k_1+\ldots + k_n$,
with points $k_i\in Y$. We consider local systems $L$ over $Y-K$. 
If $x\in Y-K$ is a choice of basepoint then a local system corresponds to the {\em monodromy representation}
$$
\rho  : \pi _1(Y-K, x)\rightarrow GL(L_x).
$$
The {\em local monodromy transformations} are the $\rho (\alpha _i)$ where $\alpha _i$ is a loop in
standard form going from $x$ to near $k_i$, once around clockwise, then back to $x$. In this paper,
we will systematically make the convention :

\begin{convention}
\label{semisimple}
The local monodromy transformations are semisimple, i.e. diagonalizable matrices.
\end{convention}

A {\em logarithmic connection} on $(Y,K)$ is a vector bundle $E$ on $Y$,
with a connection operator 
$$
\nabla : E \rightarrow E\otimes \Omega ^1_Y(\log K).
$$
The {\em monodromy} of $(E,\nabla )$ is a local system on $Y-K$, described for example as the sheaf $L=E^{\nabla}$ 
of analytic holomorphic sections $e$ of $E$ with $\nabla (e)=0$. Over $Y-K$ this is a locally constant sheaf
or a local system, and corresponds to a monodromy representation $\rho _{E,\nabla }$. 

The {\em residue} of $(E,\nabla )$ at a point $k_i\in K$ is the pair $(E_{k_i}, {\rm res}(\nabla , k_i))$
consisting of the fiber of $E$ over $k_i$, and the residue of the connection which is an endomorphism
of the fiber. 

We say that ``the residues of $\nabla $ are semisimple'' if these endomorphisms
are semisimple i.e. diagonalizable. Furthermore, in order to insure that the monodromy transformation
satisfies Convention \ref{semisimple}, it is convenient to ask that the eigenvalues of the residues
never differ by integers. Indeed, if there are pairs of eigenvalues differing by integers even in a
semisimple residue, this can typically lead to Jordan blocks of size $>1$ in the monodromy.
This gives the analogous convention for the de Rham case.

\begin{convention}
\label{semisimpleresidues}
The residues of $\nabla$ are semisimple and their eigenvalues don't differ by nonzero integers.
\end{convention}

These conventions greatly reduce the complexity of the notation and arguments required to understand Katz's constructions.
Of course Katz and subsequent authors all considered the more general
case of arbitrary Jordan normal forms, and we refer the reader to those references for a more in-depth look at this
aspect. 

\subsection{Middle cohomology}
\label{middlecoh}

If the residues don't have integer eigenvalues, then Deligne's theory gives an easy description of the 
cohomology of the monodromy local system. We denote by $DR(Y,E)$ the {\em logarithmic de Rham complex}
$$
DR(Y,E) := [E \stackrel{\nabla}{\rightarrow} E\otimes \Omega ^1_Y(\log K)].
$$
The connection $\nabla$ and the divisor $K$ are missing and implicit in this notation. This abuse allows us to shorten
most displays below. Unless otherwise stated, all de Rham complexes are supposed to be logarithmic with respect
to the relevant divisor.

\begin{proposition}
Suppose that the eigenvalues of the residues of $(E,\nabla )$ are never integers. 
Let $L_{\rho}:= E^{\nabla}$ be the monodromy local system on $Y-K$. Then hypercohomology
of the logarithmic de Rham complex
$DR(Y,E)$ calculates the cohomology of $Y-K$ with coefficients in $L_{\rho}$:
$$
\hh ^{\cdot}DR(Y,E)\sim H^{\cdot}(Y-K, L_{\rho} ).
$$
\end{proposition}

In the case where some eigenvalues are integers, the situation is more complicated. 
The same monodromy representation can come from several different logarithmic connections, whose residual
eigenvalues will differ by integers. The cohomology which is calculated by the de Rham complex will
in principle depend on which lift we have chosen. A canonical choice, somewhat different from the
choices coming from lifts, is given by the notion of ``middle cohomology''
and the ``middle de Rham complex''. The reason for the word ``middle'' is that it corresponds to the middle
perversity in intersection cohomology. In the one-dimensional case, as was well understood by Katz, the
notion of intersection cohomology corresponds to the more classical construction $j_{\ast}$,
as opposed to the derived $\rr j_{\ast}$.  In this context the word ``middle'' is more notation
than notion. 

The Betti version is as follows. 
Let $j: Y-K\hookrightarrow Y$ denote the inclusion. Assume $K$ is nonempty. 
If $\rho $ is a representation of the fundamental group of $Y-K$ corresponding to a 
local system $L_{\rho}$ on $Y-J$, then we define the {\em middle cohomology}
$$
MH^i(Y, L_{\rho}) := H^i(Y, j_{\ast} (L_{\rho})).
$$
The non-derived $j_{\ast}(L_{\rho})$ is the degree zero part of the total derived
$\rr j_{\ast}(L_{\rho})$; the other piece is $R^1j_{\ast}(L_{\rho})[-1]$. 
Which gives an exact triangle in the derived category
$$
j_{\ast}(L_{\rho}) \rightarrow \rr j_{\ast}(L_{\rho}) \rightarrow 
R^1j_{\ast}(L_{\rho})[-1] \rightarrow j_{\ast}(L_{\rho})[1] \ldots .
$$
Look at the long exact sequence of hypercohomology for this triangle. 
The hypercohomology of the total $\rr j_{\ast}$ gives the cohomology of $Y-K$.
Also $R^1j_{\ast}(L_{\rho})$ is concentrated at $K$, so it only contributes for global sections.
In particular we have 
$$
MH^0(Y, L_{\rho}) = H^0(Y-K, L_{\rho}) 
$$
and there is a long exact sequence
\begin{equation}
\label{mcohlongexact}
0\rightarrow MH^1(Y, L_{\rho}) \rightarrow H^1(Y-K,L_{\rho}) \rightarrow 
\oplus _{k\in K}R^1j_{\ast}(L_{\rho})_k \rightarrow 
MH^2(Y, L_{\rho}) \rightarrow 0.
\end{equation}
Since $K$ is nonempty, $Y-K$ is homotopic to a one dimensional complex so its $H^2$ with local coefficients
vanishes.

An observation which is important for defining the middle convolution is the following.

\begin{proposition}
\label{nontrivialpointcase}
Suppose that for at least one point $k_i\in K$, the local monodromy has no fixed vectors.
Then
the middle cohomology in degrees zero and two vanishes.
\end{proposition}
{\em Proof:}
It is clear that $H^0(Y,L_{\rho}) = 0$ because there are no flat sections near the point $k_i$, so there
can be no global flat sections. This proves that $MH^0=0$. There is a Poincar\'e-Verdier duality between 
$MH^i(Y,L)$ and $MH^{2-i}(Y,L^{\ast})$, and $L^{\ast}$ also has no fixed vectors at $k_i$. This gives
$MH^2(Y,L)=0$. One can also prove the vanishing by a direct topological argument. 
\eop

As a corollary we obtain the dimension of the middle cohomology group in this case:

\begin{corollary}
\label{mhdimension}
Suppose $L_{\rho}$ is a local system of rank $r$ on $Y-K$, and suppose that for at least one point
$k_i\in K$ the local monodromy has no fixed vectors. Then we have
$$
{\rm dim} MH ^1(Y,L_{\rho}) = r(n-2) - \sum _{j=1}^n {\rm cofix}(L_{\rho}, k_j).
$$
where $n$ is the number of points in $K$ and ${\rm cofix}(L_{\rho}, k_j)$ is the dimension of the
space of cofixed vectors of the local monodromy at $k_i$. 
\end{corollary}
{\rm Proof:}
Note that
$H^0(Y-K,L_{\rho}) = 0$ as pointed out in the proof of 
\ref{nontrivialpointcase}, and $H^2(Y-K,L_{\rho})=0$ because $Y-K$ is homotopically a $1$-dimensional complex.
Thus by calculating the Euler characteristic we have
$$
{\rm dim} H^1(Y-K,L_{\rho})=  r(n-2).
$$
On the other hand, 
$$
{\rm dim} H^1(B^{\ast}_{k_j},L_{\rho}) = {\rm cofix}(L_{\rho}, k_j).
$$
The exact sequence \ref{mcohlongexact} gives the dimension of $MH^1$.
\eop

\subsection{Middle homology}
\label{middlehom}

The {\em middle homology} is 
obtained by duality with the middle cohomology:
$$
MH_i(Y-K, L):= MH^i(Y-K, L^{\ast})^{\ast}.
$$

This is interesting only if the monodromy of $L $ has some eigenvalues equal to $1$ around a point $k_i$. 
The loop around that point, with the eigenvector as coefficient, gives a cycle in $H_1(Y-K, L)$.

Let $F_i\subset L_z$ denote the subspace of vectors fixed by the monodromy transformation $\rho _L(\alpha _i)$.
Since we are assuming that the local monodromy transformations are unipotent, the dimension of $F_i$ is
equal to the multiplicity of $1$ as eigenvalue of $\rho _L(\alpha _i)$. We get a map
$$
\phi : \bigoplus _{i=1}^k F_i \rightarrow H_1(Y-K, L).
$$

\begin{lemma}
\label{injectivity}
Suppose that the monodromy around at least one of the points $k_i$ has no fixed vector. Then
$\phi$ is injective and the first middle homology is 
the cokernel of the map $\phi$:
$$
MH_1(\Gamma , L):= \frac{H_1(Y-K , L)}{\phi  \bigoplus _{i=1}^k F_i }.
$$ 
\end{lemma}
{\em Proof:}
This is dual to Proposition \ref{nontrivialpointcase}.
\eop

\subsection{Middle de Rham cohomology}
\label{middledrcoh}

For the de Rham version of middle cohomology, if $(E,\nabla )$ is a vector bundle with logarithmic
connection on $(Y,K)$, we define the {\em middle de Rham complex} 
$$
MDR(Y,E,\nabla )= \left[ MDR^0(Y,E,\nabla )\rightarrow MDR^1(Y,E,\nabla ) \right]
$$
with
$$
MDR^0(Y,E,\nabla ):= E,
$$
and 
$$
MDR^1(E,\nabla ):= \ker \left(  E\otimes \Omega ^1_Y(\log K) \rightarrow E_K^0 \right) 
$$
where $E_K^0$ is the quotient of the fiber $E_K$ over $K$, corresponding to the $0$-eigenspaces
of ${\rm res}(\nabla , k_i)$ at the points $k_i\in K$. 
The differential is given by $\nabla$ as for the usual de Rham complex. 

We should stress here that this definition is
the right one only under our convention and assumption that the local monodromy, and the residues of 
$\nabla$, are semisimple. 

Define the {\em middle de Rham cohomology} to be the hypercohomology $\hh ^iMDR(Y,E,\nabla )$.

By definition we have a short exact sequence of complexes of sheaves on $Y$,
$$
0\rightarrow MDR(Y,E,\nabla ) \rightarrow DR(Y,E,\nabla ) \rightarrow E_K^0[-1]\rightarrow 0.
$$
This gives the same kind of long exact sequence as before.

\begin{lemma}
\label{middlebettiderham}
Suppose that the residues of $\nabla$ have no nonzero integer eigenvalues. Then the
above short exact sequence for the middle de Rham cohomology coincides after Riemann-Hilbert correspondence,
with the previous exact triangle for the middle Betti cohomology. In particular, if $L_{\rho}$ is the
monodromy local system $E^{\nabla}$ then we have
a natural isomorphism 
$$
\hh ^iMDR(Y,E,\nabla ) \cong MH^i(Y, L_{\rho}).
$$
\end{lemma}
\eop

If there are nonzero integer eigenvalues, on the other hand, then the corresponding subspaces
are fixed for the monodromy transformation but don't appear in the quotient $E^0_K$. In this case
the middle de Rham cohomology will be different from the middle Betti cohomology.

\begin{exercise}
Describe what cohomology is calculated by $DR(Y,E,\nabla )$ and $MDR(Y,E,\nabla )$ when 
the residues of $\nabla$ may have some nonzero integer eigenvalues. 
\end{exercise}

{\em Hint:} It depends on the sign of the eigenvalues.

Proposition \ref{nontrivialpointcase} above thus has the corresponding corollary in the de Rham case.

\begin{corollary}
\label{nontrivialresiduecase}
Suppose that the residues of $\nabla$ are semisimple and 
have no nonzero integer eigenvalues, and suppose that for at least one
point $k_i\in K$, the residue ${\rm res}(\nabla , k_i)$ has all eigenvalues different from $0$. Then
the middle cohomology in degrees zero and two vanishes:
$$
\hh ^0MDR(Y,E,\nabla ) = 0, \;\;\; \hh ^2MDR(Y,E,\nabla ) =0.
$$
The dimension of the middle cohomology in degree $1$ is given by 
$$
{\rm dim} MH ^1(Y,L_{\rho}) = r(n-2) - \sum _{j=1}^n {\rm rk}(E^0_{k_j}).
$$
\end{corollary}
{\em Proof:} 
This is immediate from \ref{nontrivialpointcase} and \ref{middlebettiderham}, and with Corollary \ref{mhdimension}
or its proof we get the dimension count. 
\eop

\subsection{The Betti moduli spaces}
\label{bettimodsp}

We are interested in the moduli of representations with fixed conjugacy classes at the singularities.
The first version to look at is the ``Betti'' moduli space. See \cite{CorletteSimpson}.

Let $q_1,\ldots , q_n$ be $n$ distinct points 
in $Y:=\pp ^1$ and fix a basepoint $z$ different from these. 
Put $\Gamma := \pi _1(Y-\{ q_1,\ldots , q_n\} , z)$. Let $\gamma _1,\ldots , \gamma _n$ denote standard
loops based at $z$ going around the points $q_1,\ldots ,q_n$ respectively.

Fix closed subsets $C_1,\ldots , C_n\subset GL(r)=GL(r,\cc )$
invariant under the conjugation action. In the present paper in keeping with Convention
\ref{semisimple} these will be semisimple conjugacy classes (see below). However the definition can be made with
more general closed subsets which would then have to contain many different conjugacy classes including semisimple ones.
In this case the structure of the moduli space is more complicated, for example it can be nonempty even when the
moduli space for the semisimple conjugacy classes in the closure might be empty. 

Let 
$$
{\rm Rep}(\Gamma , GL(r); C_1\ldots C_n)\subset {\rm Rep}(\Gamma , GL(r))
$$
be the closed subset of representations $\rho : \Gamma \rightarrow GL(r)$ such that 
$$
\rho (\gamma _i)\in C_i.
$$
Since it is a closed subset of an affine variety, it is also affine.
The group $GL(r)$ acts on the representation variety and it preserves our closed subset because we
have assumed that the $C_i$ are conjugation-invariant. Thus we get an action of $GL(r)$ on the affine variety
${\rm Rep}(\Gamma , GL(r); C_1\ldots C_n)$ so we can take the universal categorical quotient
$$
M_B(C_1,\ldots , C_n):={\rm Rep}(\Gamma , GL(r); C_1\ldots C_n)/GL(r).
$$

This has the following usual description on the level of points. Two points $\rho , \rho '$ of the representation
variety are {\em $S$-equivalent} if the closures of their orbits intersect. In the preimage of any point of $M_B$ there
is a unique closed orbit, which shows that this relation is an equivalence relation and the points of 
$M_B(C_1,\ldots , C_n)$
are the $S$-equivalence classes. We have the same description for the action of $Gl(r)$ on 
${\rm Rep}(\Gamma , GL(r))$, and since ${\rm Rep}(\Gamma , GL(r); C_1\ldots C_n)$ is a closed $GL(r)$-invariant
subvariety, the closure of an orbit of $\rho \in {\rm Rep}(\Gamma , GL(r); C_1\ldots C_n)$ is the same
when taken in the bigger representation variety or the closed subset. Therefore, the relation of $S$-equivalence
when we restrict the conjugacy classes, is the restriction of this relation on the full variety.\footnote{Notice that
for this statement, we have used the condition that the $C_i$ are closed subsets. If we tried to do this with
locally closed subsets, for example corresponding to nonsemisimple conjugacy classes but not their closures,
it wouldn't work the same way.} The relation of $S$-equivalence for the full representation variety
is well-understood, see Lubotsky-Magid \cite{LubotskyMagid} for example. 
In particular, two points $\rho , \rho '$ are $S$-equivalent if and only if their semisimplifications are isomorphic.
The semisimplification is again a representation in ${\rm Rep}(\Gamma , GL(r); C_1\ldots C_n)$, and 
the points of $M_B(C_1,\ldots , C_n)$ represent the isomorphism classes of semisimple representations. 

Now restrict our attention to the case of semisimple conjugacy classes $C_1,\ldots , C_n$, that is to say
the conjugacy classes of diagonalizable matrices. The $C_i\subset GL(r)$ are closed subsets,
so the above discussion applies. 

It is convenient to think of a semisimple class as being determined by a divisor on $\ggg _m$.
Write a divisor as $g= \sum _{a\in \ggg _m}m(a)[a]$ where $[a]$ is the point $a$ considered as a reduced
effective divisor and the sum is finite. If $r= {\rm deg}(g):= \sum _am(a)$,
then the divisor $g$ corresponds to the conjugacy class $C(g)\subset GL(r)$ of diagonalizable matrices
having eigenvalues $a$ with multiplicities $m(a)$. A sequence of semisimple conjugacy classes
is then represented by a {\em local monodromy vector} of $n$
divisors $\overvect{g} = (g_1,\ldots , g_n) \in Div(\ggg _m)^n$. 
We come to our main notation for the Betti moduli spaces: 
$$
M_B(\overvect{g}) := M_B(C(g_1), \ldots , C(g_n)),
$$
where the collection of points $Q= q_1+\ldots + q_n$ is implicit but not mentionned. 

The vector or partition consisting of the $m(a)$
is a partition of $r$. Kostov calls this the {\em multiplicity vector}, and the vector of multiplicity
vectors corresponding to $g_1,\ldots , g_n$ is called by Kostov the {\em polymultiplicity vector} or PMV. 
To obtain a geographic understanding one should look only at the PMV, see Roberts \cite{Roberts}.

\subsection{Nitsure's de Rham moduli space}
\label{nitsuredrmod}

We can define the following $2$-functor $\Mm _{DR}(r,d)$ of $\cc$-schemes of finite type $T$.
Put $\Mm _{DR}(r)[T]$ equal to the groupoid of $(E,\nabla )$ where $E$ is a vector bundle of rank $r$
and degree $d$ on
$P\times T$ and 
$$
\nabla : E\rightarrow E\otimes _{\Oo _{P\times T}}\Omega ^1_{P\times T/T}(\log Q\times T)
$$
is a relative logarithmic connection. Standard moduli theory shows that it is an Artin algebraic stack
locally of finite type.

Say that a logarithmic connection $(E,\nabla )$ is {\em semistable} if for any subbundle $F\subset E$
preserved by $\nabla$, we have
$$
\frac{{\rm deg}(F)}{{\rm rk}(F)} \leq  \frac{{\rm deg}(E)}{{\rm rk}(E)} .
$$
Define stability using a strict inequality for strict nonzero subbundles. 
Semistability and stability are open conditions \cite{Nitsure}, and the open substack of semistable objects
$$
\Mm _{DR}^{\rm se}(r,d)\subset \Mm _{DR}(r,d)
$$
is an Artin stack of finite type (it follows from the boundedness in Nitsure's construction \cite{Nitsure}).

Semistability of a logarithmic connection would be a consequence of semistability of the underlying bundle, but
doesn't imply it in general. 
Esnault with Viehweg \cite{EsnaultViehweg2} and
Hertling \cite{EsnaultHertling}, and also
Bolibruch \cite{Bolibruch} have studied the problem of realization of monodromy representations as logarithmic connections
on semistable bundles in the higher genus case, generalizing Bolibruch's well-known work on $\pp ^1$. Our present notion of
semistability of the pair $(E,\nabla )$ is somehow less subtle. 

Nitsure constructs in \cite{Nitsure} the {\em moduli space} which is a universal categorical quotient
$$
\Mm _{DR}^{\rm se}(r,d)\rightarrow M _{DR}(r,d).
$$ 
The points represent $S$-equivalence classes of semistable logarithmic connections, and there
is a unique polystable object in each $S$-equivalence class.

Suppose ${\bf c}_1,\ldots ,{\bf c}_n \subset {\bf gl}(r)$ are closed subsets invariant under the adjoint action of 
$GL(r)$ on its Lie algebra ${\bf gl}(r)$. Then, as before, we obtain a closed substack
$$
\Mm _{DR}^{\rm se}(r,d; {\bf c}_1,\ldots ,{\bf c}_n)\subset \Mm _{DR}^{\rm se}(r,d)
$$
consisting of logarithmic connections $(E,\nabla )$ such that up to choice of basis of $E_{q_i}$,
the residue ${\rm res}(\nabla , q_i)$ lies in ${\bf c}_i$. 

Again, here we will concentrate on the case where each ${\bf c}_i$ is the conjugacy class of a semisimple matrix,
which is closed and $GL(r)$-invariant.  As above, such a conjugacy class may be parametrized by an
effective divisor 
$g_i\in Div(\aaa ^1)$, with ${\rm deg}(g_i)=r$. 

We denote by $\overvect{g}= (g_1,\ldots ,g_n)$ a vector of divisors parametrizing semisimple conjugacy classes,
either $C_i = C(g_i)\subset GL(r)$ in the ``multiplicative case'' $g_i\in Div (\ggg _m)$
or ${\bf c}_i = {\bf c}(g_i) \subset {\bf gl}(r)$ in the ''additive case'' $g_i\in Div(\aaa ^1)$. 

The rank $r$ is recovered from $\overvect{g}$ as the degree of any one of the divisors $g_i$
(they all have to have the same degree). We can also define the {\em trace} of a divisor $g= \sum m(\alpha )[\alpha ]$
to be the sum 
$$
{\rm tr}(g) := \sum m(\alpha )\alpha \in \cc .
$$
If $A\in {\bf c}(g)$ is a matrix in the corresponding conjugacy class then $Tr (A) = {\rm tr}(g)$.
The residue formula for the logarithmic connection on the determinant line bundle of $E$ provides the formula
$$
d = {\rm deg}(E) = \sum _{i=1}^n {\rm tr}(g_i).
$$
This obviously has to be an integer, otherwise the moduli space will be empty.
The degree $d$ of the bundle $E$ may be recovered from $\overvect{g}$. Thus, we are justified in writing
$$
M _{DR}(\overvect{g}) := M _{DR}(r,d; {\bf c}(g_1),\ldots , {\bf c}(g_n)).
$$

\subsection{Deformations and obstructions}
\label{defobstr}

Suppose $(E,\nabla )$ represents a point in $M_{DR}(\overvect{g})$. We compare the deformation and
obstruction theory of $(E,\nabla )$ as logarithmic connection, with that of $(E,\nabla )$ as 
a point in the moduli stack $\Mm _{DR}(\overvect{g})$. The following result was pointed out by
N. Katz in the late 1980's. 

\begin{theorem}
\label{defobs}
The deformation and obstruction theory for the moduli stack $\Mm _{DR}(\overvect{g})$
at a point $(E,\nabla )$ is governed by the middle cohomology groups of the endomorphism
bundle $\hh ^iMDR(Y,End(E))$ for $i=0$ (automorphisms), $i=1$ (deformations) and $i=2$ (obstructions).
\end{theorem}
{\em Proof:}
We give an heuristic but ultramodern explanation. Definitions and explicitations would need to be filled in,
but this should convince the reader why it is true. 

The deformation and obstruction theory of $(E,\nabla )$ as a logarithmic connection, that is as a point in
$\Mm _{DR}(r,d)$, is given by the
$L_{\infty}$ algebra
$$
{\bf D}_{E,\nabla}(\Mm _{DR}):= \hh ^{\cdot} (Y, End(E)\otimes \Omega ^{\cdot}_Y(\log Q)).
$$
Let $R:= {\bf gl}(r) // GL(r)$ denote the moduli stack of conjugacy classes of matrices.
It is a smooth Artin algebraic stack. At a point corresponding to a matrix $A$, its deformation
theory is controlled by the $L_{\infty}$-algebra concentrated in degrees $0$ and $1$
$$
{\bf D}_A(R):= {\bf gl}(r) \stackrel{[-,A]}{\longrightarrow} {\bf gl}(r).
$$
If $A$ is in the semisimple conjugacy class ${\bf c}_i$, then the deformation theory of 
$\langle {\bf c}_i\rangle :={\bf c}_i//GL(r)$
is controlled by 
$$
{\bf D}_A(\langle {\bf c}_i\rangle )= 
[{\bf gl}(r) \stackrel{[-,A]}{\longrightarrow} ({\rm Im}(u\mapsto [u,A]))]
$$
\begin{equation}
\label{substackmap} 
= 
\ker ({\bf D}_A(R) \rightarrow {\bf gl}(r)^A)
\end{equation}
where ${\bf gl}(r)^A$ is the space of cofixed or vectors of the adjoint action of $A$, which is isomorphic
to the space of fixed vectors since $A$ is semisimple. In practical terms, if $A$ is in diagonal
form then the degree one piece of ${\bf D}_A(\langle {\bf c}_i\rangle )$
is the space of off-block-diagonal matrices (for the blocks determined by the eigenvalues of $A$)
and ${\bf gl}(r)^A$ is the space of block-diagonal matrices.

For any point $q_i\in Q$ the construction $(E,\nabla )\mapsto {\rm res}(\nabla , q_i)$ gives a morphism
of moduli stacks
$\Mm _{DR}(r,d)\rightarrow R$.
Putting these together gives a morphism 
$\Mm _{DR}(r,d)\rightarrow R^n$.
On the other hand, a vector of divisors $\overvect{g}$ represents a collection
of conjugacy classes ${\bf c}(g_i)$ which gives the substack 
$$
R(\overvect{g}):= \prod _{i=1}^n \langle {\bf c} (g_i) \rangle \subset  R^n,
$$
and by definition 
$$
\Mm _{DR}(\overvect{g}) = \Mm _{DR}(r,d)\times _{R^n} R(\overvect{g}).
$$

To get the deformation theory for this fiber product, we should take the homotopy fiber product
of the $L_{\infty}$-algebras. The one for $R(\overvect{g})$ is the kernel of a map (\ref{substackmap})
on the one for $R^n$. 

Hinich explains how to go between a sheaf of $L_{\infty}$-algebras, and a global $L_{\infty}$-algebra
\cite{Hinich}. 
In our case, ${\bf D}_{E,\nabla}(\Mm _{DR})$ is the globalization of the sheaf of $L_{\infty}$-algebras
$DR(Y,End(E))$ (the logarithmic de Rham complex along $Q$). 
Going back and forth a few times we see that the deformation theory for 
$\Mm _{DR}(\overvect{g})$ is controlled by the globalization, or hypercohomology, of the kernel
of the map corresponding to (\ref{substackmap}) on the sheaf of $L_{\infty}$-algebras, this map expressed in
local terms at the singularities as
\begin{equation}
\label{localmap}
DR(Y,End(E)) \rightarrow \bigoplus _{i=1}^n End(E) _{q_i}^0
\end{equation}
where the superscript $0$ means the trivial eigenspace for the action of ${\rm res}(\nabla , q_i)$.
The kernel is exactly the middle de Rham complex for $End(E)$. Thus, the deformation theory
is controlled by an $L_{\infty}$-algebra $\hh ^{\cdot}MDR(Y,End(E))$. 
\eop

\begin{corollary}
\label{smoothness}
Suppose $(E,\nabla )\in\Mm _{DR}(\overvect{g})$ is a point corresponding to an irreducible representation. 
Then it is a smooth point of the moduli stack, and has only scalar automorphisms so it is also a smooth point
of the moduli space where the tangent space is $\hh ^1MDR(Y,End(E))$.
\end{corollary}
{\em Proof:}
Decompose into the trace-free part and the scalars: $End(E) = End'(E)\oplus \Oo$. The trace of the obstruction 
map is zero, and $\hh ^1MDR(Y, \Oo )=0$ since $Y\cong \pp ^1$. Thus, the deformations and obstructions
are given by $\hh ^1MDR(Y,End'(E))$ and $\hh ^2MDR(Y,End'(E))$ respectively. Poincar\'e duality for 
the middle cohomology, plus the fact that $End'(E)$ is self-dual, gives 
$\hh ^2MDR(Y,End'(E))\cong \hh ^0MDR(Y,End'(E))=0$ since $E$ has no trace-free endomorphisms because it is irreducible.
Thus the space of obstructions vanishes, and the tangent space is given by $\hh ^1MDR(Y,End(E))= 
\hh ^1MDR(Y,End'(E))$.
\eop

\subsection{Dimension counting}
\label{dimcount}

From the previous discussion of deformations and obstructions, we find the following boiled-down statement.

\begin{proposition}
\label{dimensioncount}
Suppose $(E,\nabla )$ is a stable point in $M _{DR}^s(\overvect{g})$.
Then the dimension of the moduli space at the given point (or any other stable point) 
is obtained by a naive dimension count:
$$
{\rm dim}(M _{DR}^s(\overvect{g})) = \sum _{i=1}^n {\rm dim}({\bf c}(g_i)) - 2 r^2 + 2.
$$
The factor $2(r^2-1)$ corresponds to the fact that the conjugation action factors through $PGL(r)$ and the
product identity lies in $SL(r)$.

The same dimension count holds for the open subset of irreducible representations 
$M_B^{\rm irr}(\overvect{g})$ if it is nonempty. 
\end{proposition}
{\em Proof:}
Apply Corollary \ref{smoothness}. The tangent space is $\hh ^1MDR(Y,End(E))$ and 
$$
\hh ^0MDR (Y,End(E))= \hh ^2MDR(Y,End(E)) = \cc
$$
since $E$ is irreducible.
Obtain the dimension count by using the fact that $MDR(Y,End(E))$
is the kernel of the map (\ref{localmap}), noting
$$
{\rm dim}({\bf c}(g_i)) = r^2 - {\rm dim}End(E)_{q_i}^0,
$$
and calculating the Euler characteristic. 
The Riemann-Hilbert correspondence gives the corresponding
statement for $M_B^{\rm irr}(\overvect{g})$.
\eop

We introduce the {\em defect}. It may be seen as playing a role in the dimension count, but is also foremost 
related
to Katz's algorithm as we shall explain later. For each $g_i$, let $\nu (g_i)$ be the maximal multiplicity
of an eigenvalue. The centralizer of a matrix $A$ in the conjugacy class ${\bf c}(g_i)$ is the set of 
block-diagonal matrices, and the dimension of the conjugacy class is the number of positions which are not
in the block-diagonal pieces. 

Think visually of shifting all of the diagonal blocks to the left of the matrix.
In other words, transpose each square diagonal block with the rectangle consisting of all places to the left of it
in the same rows. In the resulting picture, the square blocks are now arrayed from top to bottom
flush with the left edge of the matrix. The dimension of the conjugacy class is still the number of positions
which are not in these blocks. This leftover part can 
be divided into two regions: a big rectangle of size $r\times (r-\nu (g_i))$ which is everything to the right of the
biggest block, plus a union of other rectangular regions corresponding in each
row to the positions to the right of the edge of the corresponding block for that row, but to the left of the size of the biggest block.
The second piece might be empty, indeed it is empty exactly in the case when the blocks all have the same size. 
We obtain the crude estimate
$$
{\rm dim}({\bf c}(g_i)) \geq r(r-\nu (g_i)),
$$
leading to the crude estimate for the dimension of the moduli space as
$$
{\rm dim}(M _{DR}^s(\overvect{g}))  \geq n r^2 - r\sum _{i=1}^n\nu (g_i) - 2r^2 + 2
= 2 + r\left( (n-2) r - \sum _{i=1}^n\nu (g_i)\right) .
$$
In view of this formula already, it seems reasonable to consider the quantity 
$$
\delta (\overvect{g}) := (n-2) r - \sum _{i=1}^n\nu (g_i).
$$
We call this the {\em defect} because it enters into Katz's algorithm in a remarkably elegant way: if 
$\beta$ is a convolution object corresponding to a choice of maximal-multiplicity eigenvalue for each
$g_i$, then the new rank
of the Katz-transformed local monodromy data is 
$$
r' = r + \delta (\overvect{g}) .
$$
We will want to run Katz's algorithm when $\delta (\overvect{g})$ is negative. We can do so until we get to 
a vector whose defect is positive. 

In terms of the defect, the crude dimension count says
${\rm dim}(M _{DR}^s(\overvect{g}))  \geq 2 + r\delta (\overvect{g})$.
In order to refine the dimension count, introduce the {\em superdefect} denoted locally by
$$
\sigma (g_i) := {\rm dim}({\bf c}(g_i)) - r(r-\eta (g_i)),
$$
and globally by
$$
\sigma (\overvect{g}) := \sum _{i= 1}^n\sigma (g_i).
$$
These quantities, which are always $\geq 0$, are just
the differences between the crude dimension counts and the actual dimensions. Thus we have,
when the stable open set is nonempty,
$$
{\rm dim}(M _{DR}^s(\overvect{g}))  = 2 + r\delta (\overvect{g}) + \sigma (\overvect{g}).
$$

In view of the possibility of applying Katz's algorithm to decrease the rank whenever $\delta < 0$,
the remaining case to investigate is when $\delta (\overvect{g})\geq 0$. Under this hypothesis,
the dimension of the moduli space is always $\geq 2$. The cases of dimension $0$ were the
subject of Katz's original book: they arise when one ends up with a rank one local system, at which point
it will no longer be possible to find a convoluter satisfying Convention \ref{chinontriv}.
The case of dimension $2$ is particularly interesting, although unfortunately our construction of 
\S \ref{harmonicbundles} (Corollary \ref{endresult}) will not apply.

\begin{lemma}
\label{dim2case}
The superdefect $\sigma (g_i)$ vanishes if and only if all of the eigenvalues of $g_i$ have the same multiplicity.
In the domain $\delta (\overvect{g})\geq 0$, the dimension of the moduli space will be exactly $2$,
if and only if $\delta = 0$ and the superdefects $\sigma (g_i)$ all vanish.
\end{lemma}
{\em Proof:}
Note that $\sigma (g_i)$ is the number of places left over in the complement of the pushed-left
diagonal blocks, after taking out the big rectangle. This vanishes only if all of the blocks have the same size.
The last statement follows from $\sigma (\overvect{g}) \geq 0$. 
\eop

If $\overvect{g}$ is an example of the case $\delta = \sigma = 0$ then any multiple (meaning to multiply all
of the divisors $g_i$ by the same amount) is also an example. Thus the examples of this case come in families
which are indexed by an integer $d\geq 0$. 
Following Kostov, write the type of $g_i$ as a partition of $r$, for example $(d,d)$ indicates a divisor of
the form $d[a] + d[b]$ supported at two eigenvalues $a$ and $b$ both with multiplicity $d$. 
Then $\overvect{g}$ has type given by a ``polymultiplicity vector'' $PMV(\overvect{g})$
which is a vector of partitions. 
Applying the definition of the defect we immediately see the following, due to Kostov in 
the paper \cite{KostovKappa} where he investigates explicitly the resulting list of cases:

\begin{lemma}
\label{dim2list}
In the domain $\delta (\overvect{g})\geq 0$, the only cases where 
${\rm dim}(M _{DR}^s(\overvect{g}))  = 2 $ are the following four:
\newline
---for $n=4$, $d=r/2$ and $PMV(\overvect{g}) = ((d,d), (d,d), (d,d), (d,d))$;
\newline
---for $n=3$, $d=r/3$ and $PMV(\overvect{g}) = ((d,d,d), (d,d,d), (d,d,d))$;
\newline
---for $n=3$,  $d= r/4$ and $PMV(\overvect{g}) = ((2d,2d), (d,d,d,d), (d,d,d,d))$; and
\newline
---for $n=3$, $d=r/6$ and $PMV(\overvect{g}) = ((3d,3d), (2d,2d,2d), (d,d,d,d,d,d))$.
\end{lemma}
{\em Proof:}
Apply Lemma \ref{dim2case}.
Vanishing of the superdefect means that $g_i$ is of type $(d_i,\ldots , d_i)$ for some $d_i=r/b_i$.
Vanishing of the defect says $\sum _i (1/b_i) = n-2$, and the only solutions with $b_i$ integers $\geq 2$ are 
$$
2 = \frac{1}{2} + \frac{1}{2} + \frac{1}{2} + \frac{1}{2}, \;
1 = \frac{1}{3} + \frac{1}{3} + \frac{1}{3}, \;
1 = \frac{1}{2} +  \frac{1}{4} + \frac{1}{4}, \;
1 = \frac{1}{2} +  \frac{1}{3} + \frac{1}{6}.
$$
These give the cases of the lemma.
\eop

\begin{exercise}
Classify the possible polymultiplicity vectors for $\overline{g}$ in the cases 
when the moduli space has dimension $4$ and $6$. 
\end{exercise}

\section{The diagonal configuration and its blowing up}
\label{diagonalconfig}

The convolution operation comes from the diagonal configuration consisting of vertical and horizontal
lines plus the diagonal. This kind of configuration is a
recurring theme in Hirzebruch's work \cite{Hirzebruch}. 

In what follows, put $Y:=\pp ^1$ and $Z:= \pp ^1$, and look at the product $Z\times Y$. We have the projections
$\xi :Z\times Y\rightarrow Z$ and
$\eta :Z\times Y\rightarrow Y$.

Fix a subset of distinct points $Q:= \{q_1,\ldots , q_n\}\subset \pp ^1$ and let $Q_Y$ or $Q_Z$ denote this subset
considered as a divisor in $Y$ or $Z$ respectively. 
Let 
$$
D:= Z\times Q_Y \cup Z\times Q_Z \times Y \cup \Delta \subset Z\times Y.
$$
be the divisor obtained by using the divisors $Q_Z$ and $Q_Y$ in the vertical and horizontal directions, and adding the
diagonal. Denote also by $\eta$ and $\xi$ the projections
$$
\xi :(Z\times Y)-D\longrightarrow Z-Q_Z, \;\;\; \eta : (Z\times Y)-D\longrightarrow Y-Q_Y. 
$$
The pair $(Z\times Y,D)$ is the {\em diagonal configuration}. 

The divisor $D$ does not have normal crossings, indeed the diagonal meets the other components in a series of
triple points. 
In order to obtain a variety with normal crossings compactification compatible with the projection maps,
we have to blow up the diagonal configuration at these triple crossing points .
Let $X$ be the resulting variety,
thus we have a birational map
$$
X\rightarrow Y \times Z
$$
obtained by blowing up the points $(q_1,q_1),\ldots , (q_n,q_n)$. Let $J\subset X$ denote the reduced 
inverse image of the divisor $D$. We have a decomposition
$$
J = T + U_1 + \cdots + U_n + H_1 + \cdots + H_n + V_1 + \cdots + V_n,
$$
where:
\newline
---$T$ is the strict transform of the diagonal $\Delta$;
\newline
---$H_i$ is the horizontal strict transform of
$Z \times \{ q_i \}$; 
\newline
---$V_i$ is the vertical strict transform of
$\{ q_i \} \times Y$; 
and
\newline
---$U_i$ is the exceptional divisor lying over $(q_i, q_i)$. 

These intersect as follows: each $U_i$ meets $T$, $H_i$ and $V_i$ in three distinct points.
Also $H_i$ meets $V_j$ for $i\neq j$. These intersections are transverse, and there are no other intersections.

Let $\xi : X\rightarrow Z$ denote the first projection. It is seen as going in the vertical direction, so 
it contracts $U_i+V_i$ to the point $q_i$ and indeed $\xi ^{-1}(q_i)=U_i + V_i$. 
Let $\eta : X\rightarrow Y$ denote the second projection going in
the horizontal direction, so $\eta ^{-1}(q_i)= U_i + H_i$. 

The other divisor components are mapped isomorphically onto the bases of these projections: 
$$
\xi : T\stackrel{\cong}{\rightarrow} Z, \;\;\; 
\xi : H_i\stackrel{\cong}{\rightarrow} Z
$$
and 
$$
\eta : T\stackrel{\cong}{\rightarrow} Y, \;\;\; 
\eta : V_i\stackrel{\cong}{\rightarrow} Y.
$$
These divisor components intersect transversally all fibers of $\xi$ or $\eta$ respectively. 
 
\subsection{Convoluters---the Betti version}
\label{convbetti}

The basic setup of Katz's convolution operation is to take a local system on $Y$, pull it back
to $Y\times Z$ or the blow-up $X$, tensor with a rank one local system, and push forward to $Z$ using
$\rr ^1\xi _{\ast}$. Obviously, the first step in understanding and calculating this, is to 
understand the rank one objects \cite{Budur} \cite{GoldmanXia} over $(X,J)$. We look at the ``Betti'' case of
local systems or representations of the fundamental group. 

The birational blowing-up morphism is an isomorphism outside $J$ and $D$, that is
$$
X-J \stackrel{\cong}{\rightarrow} (Z\times Y)-D.
$$
In particular, local systems on one or the other 
are the same thing. We denote generically by $\beta$ our convoluter,
in this case a local system. It is given by a representation of the fundamental group into $\cc ^{\ast}$,
a representation which factors through the abelianization as
$$
\beta : H_1(X-J, \zz )\rightarrow \cc ^{\ast}.
$$

The loops $\gamma _{H_i}$, $\gamma _{V_i}$ and $\gamma _T$ 
going around the respectively denoted components of $J$, generate the first homology of $H_1(X-J)$. For
our calculations it is convenient to include $\gamma _{U_i}$ also as generators. These are subject to the
following relations:
$$
\gamma _T + \sum _{i=1}^n \gamma _{H_i} = 0, \;\; 
\gamma _T + \sum _{i=1}^n \gamma _{V_i} = 0, \;\; 
\gamma _{U_i} = \gamma _T + \gamma _{H_i} + \gamma _{V_i} \;\; (i=1,\ldots , n), 
$$
coming respectively from a vertical $\pp ^1$ intersecting $T$ and the $H_i$; 
from a horizontal $\pp ^1$ intersecting $T$ and the $V_i$; and 
for $i=1,\ldots , n$ from a small $\Cc ^{\infty}$ deformation of the exceptional $U_i$, 
which meets $T$, $H_i$ and $V_i$ and intersects the undeformed $U_i$ negatively.
That these generate the module of relations, can be seen from a Leray spectral sequence argument.

Denote by $\beta ^{H_i}$ the monodromy of $\beta$ on the loop $\gamma _{H_i}$ and similarly for the
other generators. Thus, specifying a local system of rank one on $X-J$ comes down to specifying
$$
\beta ^{H_i}, \;\; \beta ^{V_i}, \; \beta ^{U_i}, \; \beta ^T \; \in \cc ^{\ast},
$$
subject to the relations
$$
\beta ^T \cdot \prod _{i=1}^n \beta ^{H_i} = 1, \;\;\;
\beta ^T \cdot \prod _{i=1}^n \beta ^{V_i} = 1, \;\;\;
\beta ^{U_i} = \beta ^T \cdot \beta ^{H_i} \cdot \beta ^{V_i}.
$$
Of course the last relations mean that $\beta ^{U_i}$ are redundant.

We will use our convoluters to define a convolution operation, in which the diagonal $\Delta$ plays
a primordial role. It will be important to have nontrivial monodromy around the diagonal. To simplify
notation set
$$
\chi :=\beta ^T
$$
and make the following convention.

\begin{convention}
\label{chinontriv}
The monodromy around the diagonal of our convoluter is nontrivial, that is $\chi \neq 1$.
\end{convention}

\subsection{Convoluters---the de Rham version}
\label{convdr}

We will find it most convenient to restrict to convolution with rank one logarithmic connections on
the trivial bundle, that is de Rham objects of the form $(\Oo _X, d+ \beta )$ where $\beta$ is a 
one form on $X$ with logarithmic poles along $J$. The more general case can be viewed as being 
subsumed by the theory of
parabolic logarithmic $\lambda$-connections \cite{TMochizuki3}, see \S \ref{furtherquestions}.  

In the present case, then, a convoluter is just a section
$$
\beta \in H^0(X,\Omega ^1_X(\log J)).
$$
Given $\beta$ we denote by $\beta ^{H_i}$, $\beta ^{V_i}$, $\beta ^{U_i}$ and $\beta ^T$ its
residues along $H_i$, $V_i$, $U_i$ and $T$ respectively. Note that $H^1(X,\Oo _X) = 0$,
so by Deligne's mixed Hodge theory 
$$
\int : H^0(X,\Omega ^1_X(\log J)) \stackrel{\cong}{\rightarrow} H^1((X-J)^{\rm top},\cc )
$$
and the integrals over $\gamma _{H_i}$ etc.\ are $2\pi \sqrt{-1}$ times the residues $\beta ^{H_i}$ etc.
In particular, the structure of $H_1(X-J, \zz )$ recalled in the previous subsection implies that
$\beta$ is determined by its residues, and these are subject to the equations
$$
\beta ^T + \sum _{i=1}^n \beta ^{H_i} = 0, \;\;\;
\beta ^T + \sum _{i=1}^n \beta ^{V_i} = 0, \;\;\;
\beta ^{U_i} = \beta ^T + \beta ^{H_i} + \beta ^{V_i}.
$$
The de Rham convolution object $\beta$ gives rise to a local system, or Betti convoluter $\varphi$
with
$$
\varphi ^{H_i} = e^{2\pi \sqrt{-1}\beta ^{H_i}}, \ldots . 
$$

The analogue of Convention \ref{chinontriv} is:

\begin{convention}
\label{diagresnontriv}
The residue $\beta ^T$ of $\beta$ along the diagonal is not an integer.
\end{convention}

\section{Middle convolution---Betti version}
\label{bettimiddleconv}

In this section we will work with the divisor $D\subset (Z\times Y)$, and do our computations in
braid-group style \cite{MoishezonTeicher} \cite{StrambachVolklein} \cite{DettweilerReiterMC}. 
One could alternatively use the blowing up $(Z,J)$ and give a treatment similar to the one we 
will give in the de Rham case later.

The discussion of this section is the complex geometric version of Katz's construction.
Katz gave a geometric definition of middle convolution in \cite[2.7-2.9]{Katz}. His formulae there, stated in the context of
perverse $\ell$-adic sheaves, work universally in any geometric context. The complex geometric version was defined more explicitly,
and first exploited
by Dettweiler and Reiter \cite{DettweilerReiterAlgKatz} \cite{DettweilerReiterMC} \cite{DettweilerReiterPreprint}, 
and V\"olklein, Strambach,  \cite{VolkleinBraid} \cite{StrambachVolklein}.  
They write down explicit matrices but the motivation
comes from braid-style computations.   
Kostov proposes an ingenious version of the
construction which doesn't refer to the geometric picture, but instead is based on the possibility of
multiplying the connection matrix by a scalar to get to the case of integer eigenvalues
\cite{KostovSteklov}. And, Crawley-Boevey views the construction, again in algebraic terms, as something about 
root systems.  
Boalch considered a particular example of middle convolution in a non-rigid case \cite{Boalch2}, and the link with Katz's construction was made in
\cite{DettweilerReiterPainleve}. 
In \cite{DettweilerReiterMC} following \cite[Chap. 5.1]{Katz} it is shown that the explicit matrix definition of $MC$ has a geometric
or cohomological interpretation as
a higher direct image---this is the point of view we adopt here. 
The braid-style calculations of group cohomology necessary to get the local form of monodromy out of this geometric definition
were done in 
\cite{DettweilerReiterMC} but using the Pochammer  basis for the group cohomology classes, rather than a standard basis as we shall use here. 
In spite of the numerous references on this subject, we go through the details, where possible
keeping simplifying assumptions for our expository purpose. 

\subsection{Definitions}
\label{sec-definitions}

Recall that 
$$
\eta : Z\times Y\rightarrow Y ,\;\;\; \xi : Z\times Y\rightarrow Z
$$
are the projections, and use the same name for the projections on the open subset $(Z\times Y)-D$. 
A {\em convoluter} $\beta$ is a rank one local system on $(Z\times Y)-D$. 

Define the {\em raw convolution} $RC_{\beta}$ as follows. 
If $L$ is a local system on $Y-Q_Y$
then put
$$
RC_{\beta}(L):= R^1\xi _{\ast}(\beta \otimes (\eta ^{\ast}L)).
$$
It is a local system on $Z-A_Z$. The {\em middle convolution} $MC_{\beta}(L)$ will be a subsystem of  
$RC_{\beta}(L)$ the kernel of the map to some 
natural essentially local systems coming from triviality of certain pieces of the 
local monodromy transformations. It corresponds to replacing the cohomology of the fibers in the
$R^1\xi_{\ast}$ construction, by the middle cohomology discussed in \S \ref{middlecoh} above. 

In order to make explicit calculations, we will adopt the viewpoint of homology rather than
cohomology. Let $RC^{\ast}_{\beta}(L)$ denote the local system obtained by taking the homology of the
fibers with coefficients in $\beta \otimes \eta ^{\ast}L$. Let $MC^{\ast}_{\beta}(L)$ denote the
quotient corresponding to ``middle homology'' defined by duality with middle cohomology (the precise
definition will appear in our discussion below). If we let $\beta ^{\ast}$ and $L^{\ast}$ denote the
dual local systems, then by the duality between homology and cohomology we have 
$$
RC^{\ast}_{\beta}(L) = \left(  RC_{\beta ^{\ast}}(L^{\ast}) \right) ^{\ast}
$$
and similarly for $MC^{\ast}$. Thus it is equivalent if we look at homology, and it is easier to
visualize geometrically classes in homology with local coefficients.

\subsection{Computations in group homology}
\label{sec-compgrouphom}

We would like to calculate the local monodromy transformations of the raw and then middle convolutions. 
In order to do this, we transform the question into a computation of the action of the fundamental group of the base,
on the group homology of the fiber. See \cite{VolkleinBraid} \cite{DettweilerReiterMC} \cite{Dettweiler}
for example. 

In order to speak of fundamental groups, we need
to choose basepoints. 
Choose a basepoint $b\in Y-Q$. For $z\in Z-Q_Z$ the fiber of $X-J$ over 
$z$ is $Y-Q-\{ z\}$. This has $(z,b)$ as basepoint whenever $z\neq b$. In particular,
in order to get a fibration of based spaces we should additionally take the point $b$
out of the base.  For this reason, put $Q_Z^b:= Q_z\cup \{ b\}$. 

On the other hand, we would like to consider the fundamental group of $Z-Q_Z^b$. 
Choose another basepoint 
$c\in Z-Q_Z^b\subset \pp ^1\cong Y$. 
In the fiber over $c\in Z$ we have the complement of $Q_Y$ and should also take out the diagonal 
point $(c,c)$. Thus, let $Q_Y^c := Q_Y\cup \{ c\}$. 

The fiber $Y_c = \{ c\} \times (Y-Q_Y^c)$ of the projection $\xi :  X-J\rightarrow Z-Q_Z^b$ over $c$, 
is an open Riemann surface pointed by
the basepoint $(c,b)$. 

Let $\Gamma :=  \pi _1(Y_c, (c,b))) \cong \pi _1(Y-Q_Y^c , b)$.
It is a free group. The fundamental group of the base 
$\Upsilon := \pi _1(Z-Q_Z^b, c)$
acts on $\Gamma$. We will describe the action in greater detail below. Denote the action by 
$u\mapsto (\gamma \mapsto \mu (u,\gamma ))$
for $u\in \Upsilon$, $\gamma \in \Gamma$. 

Make the convention for group composition that $ab$ means $b$ followed by $a$. That way, a monodromy
representation indicates transport of sections along the path and satisfies $\rho (ab)=\rho (a)\rho (b)$. 

The local system $RC^{\ast}_{\beta}(L)$ restricted to $Z-Q_Z^b$ can be described as follows.
The local system $L$ corresponds to a representation $\rho$ of $\Gamma$ on the vector space $L_b$,
invariant with respect to the action of $\Upsilon$ in the sense that
$$
\rho (\mu (u,\gamma )) =\rho (\gamma )
$$
for any $\gamma \in \Gamma$ and $u\in \Upsilon$. Similarly, the local system $\beta$ 
corresponds to a pair of characters
$$
\beta ^{c,\cdot }:\Gamma \rightarrow \cc ^{\ast},\;\;\;
\beta ^{\cdot , b} : \Upsilon \rightarrow \cc ^{\ast}
$$
and the first of these is again invariant,
$\beta ^{c,\cdot} (\mu (u,\gamma )) = \beta^{c,\cdot} (\gamma )$.
Tensoring together we obtain a representation denoted 
$$
\rho ^{\beta} : \gamma \mapsto \beta ^{c,\cdot}(\gamma )\rho (\gamma ),
$$
again invariant with respect to the action of $\Upsilon$. The local system $RC_{\beta}(L)$
(resp. $RC^{\ast}_{\beta}(L)$)
corresponds to the vector space $H^1(\Gamma ,\rho ^{\beta})$  (resp. $H_1(\Gamma ,\rho ^{\beta})$).
The action of $\Upsilon$ is obtained by the natural action, tensored with the character 
$\beta ^{\cdot ,b}$. This tensorization is due to the fact that the local system $\beta$ is not
trivial on the basepoint section $z\mapsto (z,b)$ over $Z-Q^b_Z$. 

Denote by $V$ the vector space $\cc ^r$ on which the representation $\rho ^{\beta}$  is defined. 
Denote by $H^1(\Gamma , V)$ the cohomology and $H_1(\Gamma , V)$ the homology.

Fix generators for $\Gamma$ as follows: we have loops $\alpha _1,\ldots , \alpha _n$ going clockwise around the points
$q_1,\ldots , q_n$ in the standard way, 
and $\delta$ going clockwise around the point $(c,c)\in \Delta$ in the fiber $ \{ c\} \times (Y-Q_Y^c)$.
The group $\Gamma$ has generators $\alpha _i, \delta$ subject to the single relation
$$
\delta \alpha _1 \cdots \alpha _n = 1.
$$ 
Using this relation any one of the generators could be ignored but it will be more convenient to keep all of them.

The character $\beta^{c,\cdot} $ acts on these generators as follows:
$\beta^{c,\cdot} (\alpha _i) = \beta ^{V_i},
\;\;\; \beta^{c,\cdot} (\delta ) = \beta ^T$.
Thus we have
$$
\rho ^{\beta}(\delta  ) = \beta ^T \cdot 1^r, \;\;\;
\rho ^{\beta}(\alpha _i  ) = \beta ^{V_i}\rho (\alpha _i ),\;\; i= 1,\ldots , n .
$$

The homology $H_1(\Gamma , V)$ is the homology at degree one of the sequence
$$
C_2(\Gamma , V) \rightarrow C_1(\Gamma , V) \rightarrow C_0(\Gamma , V ) = V.
$$
Furthermore, $C_1(\Gamma , V)$ is the $\cc$-vector space formally generated by the symbols
$G(\gamma , v)$ where $\gamma \in \Gamma$ and $v\in V$, subject only to the relation of $\cc$-linearity in the variable
$v$. And $C_2(\Gamma , V)$ is generated by symbols $Q(\gamma , \xi , v)$ where $\gamma , \xi \in \Gamma $ and
$v\in V$. Geometrically, $G(\gamma , v)$ represents a cycle which starts with value $v$ and continues along the
path $\gamma$. And $Q(\gamma , \xi , v)$ represents a simplex whose sides are $\gamma$, $\xi$ and $\xi \gamma$
with coefficient $v$ at the starting point. 

The boundary operators for the complex are
$$
\partial G(\gamma , v) = \gamma (v) - v,\;\;\; 
\partial Q(\gamma , \xi , v) = G(\gamma ,v) + G(\xi , \gamma (v)) - G(\xi \gamma , v), \;\;\; 
\partial \circ \partial = 0.
$$ 

We work with the vector space $C_1(\Gamma , V) / \partial C_2(\Gamma , V)$ denoted just $C_1/\partial C_2$ for short. 
It is finite dimensional,
and its elements are $\cc$-linear combinations of classes denoted $G[\gamma , v]$ which are now subject to the
relations that this symbol is $\cc$-linear in $v$, and that 
$$
G[\xi \gamma , v] = G[\gamma ,v] + G[\xi , \gamma (v)].
$$
If we fix a basis $\{ v_j\}$ for $V$ then from the set of generators of 
$\Gamma$ we obtain a basis for  $C_1/\partial C_2$ consisting of the
$G[\alpha _i, v_j] \;\;\mbox{and}\;\; G[\delta , v_j]$.
It will be useful in what follows to have a formula for multiple products. For example
$$
G[\eta \xi \gamma , v] = G[\gamma , v] + G[\xi , \gamma (v)] + G[\eta , \xi \gamma (v)]
$$
and more generally
$G[\gamma _1\cdots \gamma _m , v] = \sum _{i=1} ^m G[\gamma _i, \gamma _{i+1}\cdots \gamma _m (v)]$.
Similarly for the inverse, the equation
$$
0=G[\gamma \gamma ^{-1}, v] = G[\gamma ^{-1},v] + G[\gamma , \gamma ^{-1}(v)]
$$
gives
$G[\gamma ^{-1},v] =- G[\gamma , \gamma ^{-1}(v)]$.

Now consider the action of $\Upsilon$. What we call the ``natural action'' is the one coming from
the action on the explicit generators written above. This corresponds to tensoring to trivialize the
restriction of the local system on the basepoint section. In the end, since the local system is not
trivial on the basepoint section, we will have to take the natural action tensored with the
character $\beta ^{\cdot , b}$. 

The generator $u_i$ of $\Upsilon$ corresponds to a path where the point $c$ goes around the point $q_i$.
There is some choice about how to arrange this picture, with respect to the picture of the standard
generators of $\Gamma$. Think of the points $q_1,\ldots , q_n$ as lined up in a row, with the basepoint
$b$ off to one side so that the points are arrayed from left to right when viewed from $b$.
The paths $\alpha _i$ go straight from $b$ to $q_i$, once around clockwise, then back to $b$.
On the other hand, let $c$ be on the other side of the row of points $q_i$. We obtain a number of paths
$\delta = \delta _1,\ldots , \delta _n$ going from $b$ to $c$, around clockwise, and back to $b$. These are
defined by saying that the starting and ending path for
$\delta _i$ goes just to the left of the point $q_i$, for $1<i\leq n$ the path
goes between $q_{i-1}$ and $q_i$.  We have the relation
$$
\delta _{i+1} = \alpha _i \delta _i \alpha _i ^{-1}.
$$
In particular the $\delta _i$ are all conjugate to $\delta = \delta _1$, 
which implies that $\rho ^{\beta}(\delta _i)$ are always
multiplication by $\beta ^T$. 

Now, define $u_k$ as the path which sends $c$ straight to $q_k$, around clockwise, and
back to its starting point. This happens on the other side of our picture from
the paths starting at $b$. 

With this picture, the action of $u_k$ doesn't change the $\alpha _j$ for $j\neq k$. On the other hand,
we have a Dehn twist between $\delta _k$ and $\alpha _k$. Notice that the introduction of the different
conjugates $\delta _k$ allows us to represent these Dehn twists uniformly for each $k$; if we try to
write down the formula with $\delta = \delta _1$ it becomes more complicated. 

A geometric look at the picture of $c$ going clockwise around $q_k$ yields:

\begin{proposition}
\label{actioncomp}
The action of $u_k$ on $\Gamma$ is given by
$$
\mu (u_k, \alpha _k ) = \delta _k^{-1}  \alpha _k \delta _k ,
$$
and 
$$
\mu (u_k, \delta _k) =  \delta _k^{-1}\alpha _k^{-1}\delta _k\alpha _k \delta _k . 
$$
\end{proposition}
\eop

\begin{exercise}
Define conjugates $\zeta _i$ where the point $c$ goes in between $q_{i-1}$ and $q_i$, around $b$, and back. 
Describe the action of $\zeta _i$. 
We have the relations $u_1\cdots u_{i-1}\zeta _i u_i \cdots  u_n = 1$ in $\Upsilon$. 
Check that the action defined by the above formulae
for the $u_i$ plus the formulae for $\zeta _i$, satisfies these relations. 
\end{exercise}

\subsection{The local monodromy transformations}
\label{sec-localmonodromy}

Now we would like to compute the eigenvalues of the monodromy transformations. 
This computation is local around one of the points $q_k$. For our present purposes we only ask for the
Jordan normal form of the monodromy transformation. More precise information, in fact the explicit
monodromy matrices with respect to the Pochammer basis, are obtained in \cite{DettweilerReiterMC},
\cite[Lemma 3.3.5, Proposition 3.3.6]{Dettweiler}. Our computation is along these lines but 
we don't need to consider the Pochammer elements. 

To reduce notation put
$$
\chi := \beta ^T, \;\;\; \beta _ i := \beta ^{H_i}.
$$
The values of $\beta ^{V_i}$ don't matter, as we are initially calculating the natural action trivialized over
the basepoint section. Our previous formulae become
$$
\rho ^{\beta}(\delta_k  ) = \chi ,
\;\;\;
\rho ^{\beta}(\alpha _i  ) = \beta _i\rho (\alpha _i ).
$$
For the action of $u_k$, we have
$u_k \cdot G[\alpha _i , v_j] = G[\alpha _i, v_j],\;\; i\neq 1$.
On the other hand,
$$
u_k \cdot G[\alpha _k , v_j] = G[\delta_k ^{-1}\alpha _k\delta_k , v_j] = 
G[\delta_k ^{-1}, \alpha _k \delta_k (v_j)] + G[\alpha _k , \delta_k (v_j)] + G[\delta_k ,v_j]
$$
$$
= G[\alpha _k , \delta_k (v_j)] + G[\delta_k ,v_j]
- G[\delta_k , \delta_k ^{-1}\alpha _k \delta_k (v_j)]
$$
and using the formula for the action of $\delta_k$ which is by multiplication by $\chi$ (in particular it
commutes with the $\alpha _k$),
$$
u_k \cdot G[\alpha _k , v_j] = \chi  G[\alpha _k, v_j] + G[\delta_k , (1 - \alpha _k)(v_j)].
$$
Finally, 
$$
u_k\cdot G[\delta_k , v_j] = G[\delta_k ^{-1}\alpha _k^{-1}\delta_k \alpha _k \delta_k , v_j] 
$$
$$
= -G[\delta_k , \delta_k ^{-1}\alpha _k^{-1}\delta_k \alpha _k \delta _k(v_j)] 
- G[\alpha _k ,  \alpha _k^{-1}\delta_k \alpha _k \delta (v_j)]
+ G[\delta_k , \alpha _k \delta_k (v_j)]
+ G[\alpha _k, \delta _k(v_j)] + G[\delta _k, v_j]
$$
$$
= -\chi G[\delta_k , v_j] 
-\chi  ^2G[\alpha _k ,  v_j]
+ \chi G[\delta _k, \alpha _k (v_j)]
+ \chi G[\alpha _k, v_j] + G[\delta _k, v_j]
$$
$$
= (\chi - \chi ^2)G[\alpha _k, v_j ] + G[\delta_k , (\chi (\alpha _k -1)+ 1)v_j].
$$

\begin{lemma}
\label{matrix}
Suppose $v_j$ is an eigenvector of $\rho (\alpha _k)$ with eigenvalue $r_{k,j}$. Denote
$\chi := \beta ^T$. 
Then the two-dimensional subspace of
$C_1/\partial C_2$ generated by $G[\alpha _k , v_j]$ and $G[\delta_k ,v_j]$  is invariant under the transformation
$u_k$, and on this subspace (with the two generators taken as basis vectors) the 
transformation $u_k$ has matrix
$$
u_k |_{\langle G[\alpha _k , v_j], G[\delta _k,v_j] \rangle } = 
\left(
\begin{array}{cc}
\chi         &  (\chi - \chi ^2)\\
(1-\beta _k r_{k,j})  &  1+\chi (\beta _k r_{k,j}-1)
\end{array}
\right) .
$$
\end{lemma}
{\em Proof:}
In the computations above, the action of $\alpha _k$ is by the representation
$\rho ^{\beta}$, and $v_j$ is an eigenvector of $\rho ^{\beta}(\alpha _k)$ but this time with eigenvalue
$\beta _kr_{k,j}$. Thus we should take the previous formulae and replace 
$\rho ^{\beta}(\alpha _k) v_j$ by $\beta _k r_{k,j}v_j$, which gives the stated matrix.
\eop

\begin{corollary}
\label{eigenvalues}
In the situation of the previous lemma, the eigenvalues of $u_k$ acting on the two dimensional subspace
$\langle G[\alpha _k , v_j], G[\delta_k ,v_j] \rangle$
are $1$ and $\chi \beta _k r_{k,j}$. 
\end{corollary}
{\em Proof:}
The determinant of the matrix in the lemma is 
$$
\chi + \chi ^2 (\beta _k r_{k,j}-1) - (\chi - \chi ^2)(1-\beta _k r_{k,j}) = 
\chi + \chi ^2 \beta _k r_{k,j} - \chi ^2 -\chi + \chi \beta _k r_{k,j} + \chi ^2 -\chi ^2\beta _k r_{k,j}
$$
$$
= \chi \beta _k r_{k,j}.
$$
The trace is $1+\chi \beta _k r_{k,j}$. 
The eigenvalues satisfy two equations which clearly hold for $1$ and $\chi \beta _k r_{k,j}$. 
\eop

Suppose $\chi \beta _k r_{k,j}=1$, then the matrix in the above lemma is
$$
\left(
\begin{array}{cc}
\chi         &  (\chi - \chi ^2)\\
(1-\chi ^{-1})  &  2-\chi 
\end{array}
\right)  =  1 + (\chi - 1)
\left(
\begin{array}{cc}
1      &   - \chi \\
\chi ^{-1}  &  -1 
\end{array}
\right)  ,
$$
that is $1$ plus a rank one matrix whose square is zero. In this case the
$2\times 2$ matrix of Lemma \ref{matrix} is not semisimple. Therefore, keep the following 
restriction on our eigenvalues.

\begin{convention}
\label{chirhobeta}
We assume that $\beta$ has the property that
$\chi \beta _k r_{k,j}\neq 1$ for all eigenvalues $r_{k,j}$ of $\rho (\alpha _k)$. In other words,
the matrix $\chi \rho ^{\beta}(\alpha _k)$ has  only nontrivial eigenvalues.
\end{convention}

\begin{lemma}
\label{cycles}
Suppose Convention \ref{semisimple} so $\rho (\alpha _k)$ is a semisimple matrix, and
Convention \ref{chirhobeta} so that $\chi \beta _k r_{k,j}\neq 1$. 
Then the monodromy transformation of $C_1/\partial C_2$ around $u_k$ looks up to conjugacy 
like 
$$
\chi \rho ^{\beta }(\alpha _k) \oplus 1^m
$$
where $m$ is given by a dimension count.
\end{lemma}
{\em Proof:}
The map
$$
V\oplus V \rightarrow C_1/\partial C_2, \;\; (u,v)\mapsto G[\alpha _k, u] + G[\delta _k, v]
$$
is injective. The action of $u_k\in 
\Upsilon$ preserves the image and there it acts as $\chi \rho ^{\beta }(\alpha _k)\oplus 1^r$.
Thus the multiplicity of an eigenvalue in the action of $u_k$ is at least as
big as its multiplicity in $\chi \rho ^{\beta}(\alpha _k)$.

On the other hand, the images of the maps 
$V\rightarrow H_1(\Gamma , V)$ given by $v\mapsto G[\alpha _i,v]$ for $i\neq k$,
span a subspace on which $u_k \in \Upsilon$ acts trivially, and with the subspace of the previous paragraph,
these two subspaces generate $C_1/\partial C_2$. We obtain a surjective $u_k$-equivariant map from a
representation of the form $\chi \rho ^{\beta }(\alpha _k) \oplus 1^{m'}$ to $C_1/\partial C_2$.
Thus, the multiplicity of a nontrivial eigenvalue in $u_k$ is at most its multiplicity in 
$\chi \rho ^{\beta}(\alpha _k)$. This surjection also shows that the action of $u_k$ on
$C_1/\partial C_2$ is semisimple. 

The condition that the eigenvalues of $\chi \rho ^{\beta}(\alpha _k)$ be all nontrivial means that the multiplicities
are the same as their multiplicities in $u_k$. This gives the direct sum decomposition of the lemma.
\eop

\begin{exercise}
Calculate the action of $\zeta _k$. After going to the action on the 
homology which is the kernel of the boundary map
$$
H_1 = Z_1/\partial C_2 = \ker (C_1/\partial C_2 \rightarrow C_0\cong V),
$$
the operator $\zeta _k$ should act by multiplication by a scalar. After
tensoring with $\beta ^{\cdot , b}$ it should give the identity since the local system $RC^{\ast}_{\beta}(L)$
doesn't depend on the choice of basepoint and hence
extends across $\{ b\}$. 
\end{exercise}

We now consider the action of $\Upsilon$ on $H_1:= \ker (C_1/\partial C_2 \rightarrow C_0)$.

\begin{proposition}
\label{rawconvolution}
Suppose that $\rho$ is an irreducible representation of rank $r>2$. Suppose 
that the monodromy transformations
$\rho (\alpha _i)$ are semisimple (Convention \ref{semisimple}), suppose 
that $\chi \neq 1$ (Convention \ref{chinontriv}), and suppose  that the
eigenvalues $r_{k,j}$ of $\rho (\alpha _k)$ are different from $\chi ^{-1}\beta _k ^{-1}$ (Convention \ref{chirhobeta}).
The dimension of $H_1(\Gamma , \rho ^{\beta})$ is equal to $(n-1)r$. 
The group $\Upsilon$ acting on $H_1(\Gamma , \rho ^{\beta})$ by the raw convolution
representation has the following effect on the generators:
$$
u_k\mapsto \beta ^{V_k}\otimes ( \chi \cdot \rho ^{\beta }(\alpha _k) \oplus 1^{(n-2)r}).
$$
\end{proposition}
{\em Proof:}
The conditions of Proposition \ref{nontrivialpointcase}
apply even to the dual local system, therefore $H_0=0$. In particular, the boundary map is surjective
onto $C_0$. 
Recall that $H_2=0$ because we
are looking at an open curve. 
The Euler characteristic of the complement of $n+1$ points is $(n-1)$
which gives ${\rm dim}H_1 = (n-1)r$.

On the other hand, the natural action of $\Upsilon$ on $C_0 \cong V$ is trivial. Thus, for
the monodromy transformations of the $u_k$, the kernel $H_1$ of the boundary map
contains all of the nontrivial part. For this action, using the dimension count 
and Lemma \ref{cycles},
the matrix of the action of $u_k$ is $\chi \cdot \rho ^{\beta }(\alpha _k) \oplus 1^{(n-2)r}$.
As pointed out at the start of the computation, we then have to tensor with the character $\beta ^{\cdot , b}$ to
get the representation corresponding to the raw convolution.
\eop

\subsection{Middling}
\label{middling}

Suppose that $\rho ^{\beta}$ has some eigenvalues equal to $1$ around a point $q_i$. 
The loop around that point, with the eigenvector as coefficient, gives a cycle in $H_1(\Gamma , \rho ^{\beta})$
which will be  covariant under $\Upsilon$. Going to the middle convolution, replacing $H_1$ by the middle
version, gets rid of these invariant cycles. 

We are assuming that $\chi \neq 1$ so this behavior doesn't occur at the point $c$, and in particular the point 
$c$ serves as a point where there are no fixed vector so we can apply Lemma \ref{injectivity} above. 

Let $F_i\subset L_b$ denote the subspace of vectors fixed by the monodromy transformation $\rho ^{\beta}(\alpha _i)$.
Since we are assuming that the local monodromy transformations are unipotent, the dimension of $F_i$ is
equal to the multiplicity of $1$ as eigenvalue of $\rho ^{\beta}$. We get a map
$$
\phi _i : F_i \rightarrow H_1(\Gamma , \rho ^{\beta}), \;\; \phi _i(f):= G[\alpha _i,f].
$$ 
These put together to give
$\phi : \bigoplus _{i=1}^n F_i \rightarrow H_1(\Gamma , \rho ^{\beta})$.
Recall that
$$
MH_1(\Gamma , \rho ^{\beta}):= \frac{H_1(\Gamma , \rho ^{\beta})}{\phi  \bigoplus _{i=1}^k F_i }.
$$ 
The group $\Upsilon$ acts on $MH_1(\Gamma , \rho ^{\beta})$. As before,
we can calculate with the natural action trivialized on the basepoint section,
which should then be tensored with the character 
$\beta ^{\cdot , b}$ to obtain the middle coconvolution $MC^{\ast}_{\beta}(L)$.

\begin{lemma}
\label{phiequivariant}
The map $\phi$ is equivariant for the natural action on the target,
and with $u_k$ acting trivially on $F_i$ for $i\neq k$, and by multiplication by $\chi$
on $F_k$. Also, $\phi$ is injective. 
\end{lemma}
{\em Proof:}
From the previous subsection, the action of $u_k$ preserves $G[\alpha _i,f]$ for $i\neq k$.
If $f$ is a fixed vector for $\rho ^{\beta}(\alpha _k)$ then it is an eigenvector with 
eigenvalue $\beta _k r_{k,j}= 1$. In the matrix of Lemma \ref{matrix}, we get that
the image of $G[\alpha _k,f]$ is $\chi G[\alpha _k,f]$. Injectivity of $\phi$ follows from Lemma
\ref{injectivity}. 
\eop

This lemma leads to the computation of the monodromy action of $u_k$ on the middle homology.
Write $V = V' \oplus V''$ where $\rho (\alpha _k)$ acts with eigenvalue $\beta _k^{-1}$ on
$V'$, and with eigenvalues distinct from $\beta _k^{-1}$ on $V''$. 
Thus, for $\rho ^{\beta}(\alpha _k)$ the fixed subspace is $V'$ with its complement $V''$.
Recall \ref{rawconvolution} that before tensoring with $\beta ^{\cdot , b}$,
the natural action
of $u_k$ on $H_1(\Gamma , \rho ^{\beta})$ decomposes as $\chi \rho ^{\beta}(\alpha _k) \oplus 1^{(n-2)r}$.
The underlying vector space decomposes as $V' \oplus V'' \oplus \cc ^{(n-2)r}$, and 
$u_k$ acts by $\chi$ on $V'$, by eigenvalues different from $\chi$ on $V''$, and trivially on 
$\cc ^{(n-2)r}$. 

Convention \ref{chirhobeta} says that if $r_{k,j}$ is an eigenvalue of $\rho (\alpha _k)$
then $\chi \beta _k r_{k,j}\neq 1$. The $\beta _kr_{k,j}$ are the eigenvalues of $\rho ^{\beta}(\alpha _k)$.
This condition therefore says that the eigenvalues of $u_k$ on $V''$ are different from $1$. 
Convention \ref{chinontriv} says that $\chi \neq 1$. Therefore the three subspaces in the
above decomposition of $H_1(\Gamma , \rho ^{\beta})$ are distinguished by the eigenvalues of
$u_k$.

Lemma \ref{phiequivariant} now implies that $\phi$ sends $F_i$ into the part $\cc ^{(n-2)r}$,
and sends $F_k$ into the part $V'$. On the other hand, $V'$ is the space of fixed vectors of $u_k$,
isomorphic (and indeed, equal) to $F_k$. Therefore, in the middle homology there is no remaining
eigenspace for $\chi$, the term $V''$ remains intact, and the trivial eigenspace is reduced by
an appropriate amount, to a size given by the dimension count. We can state this as follows.

\begin{proposition}
\label{mcdecomposition}
Assume Conventions \ref{semisimple}, \ref{chinontriv}, and \ref{chirhobeta}. The action 
of $u_k \in \Upsilon$ on the middle homology is semisimple, and we have a decomposition
$$
MH_1(\Gamma , \rho ^{\beta}) = V'' \oplus \cc ^m,
$$
where $V''$ is the direct sum of all eigenspaces of $\rho (\alpha _k)$ for eigenvalues different
from $\beta _k^{-1}$, and $m$ is given by a dimension count. The natural action of $u_k$ on $V''$ is
by $\chi \beta _k \rho (\alpha _k)$, and the natural action on $\cc ^m$ is trivial. The middle convolution
action is obtained by multiplying everything by $\beta ^{V_k}$. 
\end{proposition}
\eop

To put this another way, suppose $a$ is an eigenvalue of $\rho (\alpha _k)$ of multiplicity $m_k(a)$.
Then the corresponding eigenvalue of the action of $u_k$ on the middle convolution is:
\newline
---$\beta ^{V_k}\chi \beta ^{H_k} a = \beta ^{U_k}a$ 
with the same multiplicity $m_k(a)$ when $\beta ^{H_k}a\neq 1$; or
\newline
---$\beta ^{V_k}$ with multiplicity $m_k(a)+{\rm dim}MH_1 - r$ when $a = (\beta ^{H_k}) ^{-1}$.
\newline
This is seen by recalling that $\beta _k:= \beta ^{H_k}$, $\chi := \beta ^T$ 
and $\beta ^{U_k}= \beta ^{V_k} \beta ^T \beta ^{H_k}$. 

\begin{exercise}
The product of all the eigenvalues for all singular points must be $1$. As a reality-check,
see that this remains true for
the middle convolution with the above formulae.
\end{exercise}

\subsection{The Katz transformation on the level of local monodromy}
\label{kclocmon}

We create some notation for describing the effect of the middle convolution operation on local monodromy. 
Let $\LL$ denote an abelian group with the group law written multiplicatively. 
Define $Div(\LL)$ to be the free abelian group generated by points of $\LL$.
An element of $Div(\LL)$ is thus a finite linear combination $g=\sum _{\alpha \in \LL} m(\alpha )\cdot [\alpha ]$ 
with $m(\alpha )\in \zz$
and $m(\alpha ) = 0$ for almost all $\alpha \in \LL$. The divisor is {\em effective} if all the coefficients are
positive $m(\alpha )\geq 0$. 

The elements of $\LL$ are thought of as representing possible eigenvalues, and elements of $Div(\LL)$ represent
conjugacy classes of semisimple matrices with these eigenvalues. The cases of interest are $\LL =\ggg_m$,
which applies to the Betti case of the present chapter, and $\LL = \aaa ^1$ which will apply for the
de Rham case in the next chapter.

Since we are restricting in this paper to the
case of semisimple local monodromy, we can use the simpler $Div(\LL )$ rather than the set of all Jordan normal 
forms \cite{Katz} \cite{KostovSteklov} \cite{CrawleyBoeveyIHES} \cite{CrawleyBoeveyShaw} \cite{Roberts} etc.
The rank of the matrix is the degree of the divisor, that is the sum of the coefficients $m(\alpha )$. 
Denote this by $|g|$.
Define the {\em determinant} to be $det(g):= \prod _{\alpha \in \LL}
\alpha ^{m(\alpha )} \in \LL$, 
well defined since almost all factors are the identity element $1_{\LL}$. For obvious reasons
when the operation of $\LL$ is conventionally denoted additively we write $Tr(g)$ rather than $det(g)$.  

Fix $n$. A {\em local monodromy vector} is an $n$-tuple of elements of $Div(\LL)$, denoted 
$$
\overvect{g}= (g_1,\ldots , g_n) \in Div(\LL )^n,
$$
such that the degrees are the same, $|g_1| = \ldots = |g_n|$. Denote this common degree by $r(\overvect{g})$ 
and call it the {\em rank} of $\overvect{g}$ because it will correspond to the rank of the local system. 
Define the {\em total determinant} to be the product 
$$
Det(\overvect{g}) := det(g_1)\cdots det(g_n).
$$
In order to be a candidate for the local monodromy vector of a local system, we must have $Det(\overvect{g})=1$.

A {\em convoluter} is a 
function 
$$
\beta : H_1(X-J,\zz )\rightarrow \LL
$$
which, in view of the generators and relations for $H_1(X-J,\zz )$, can be thought of as a vector 
$$
\beta  = (\beta ^{H_1}, \ldots , \beta ^{V_1}, \ldots ,\beta ^{U_1},\ldots , \beta ^T)\in \LL ^{3n+1}
$$
subject to the relations 
$$
\beta ^{H_1}\cdots \beta ^{H_n} \cdot \beta ^T = 1, \;\;\; 
\beta ^{V_1}\cdots \beta ^{V_n} \cdot \beta ^T = 1,\;\;\; 
\beta ^{U_i} = \beta ^{H_i}\beta ^{V_i}\beta ^T \;\; (1\leq i \leq n).
$$

As pointed out above in the Betti (\S \ref{convbetti}) and de Rham (\S \ref{convdr}) cases, 
a convoluter contains the data necessary for defining a rank one object on 
$X$, which will also be denoted by $\beta$. 
The coefficients correspond to the local monodromy around the divisors $V_i$, $H_i$, $U_i$ and 
the diagonal $T$ respectively. In this picture the group $\LL$ is the group of possible local  monodromy for
rank one objects, which depends on what kind of object we are considering. 

The {\em Katz operation on semisimple local monodromy} assigns to a local monodromy vector $\overvect{g}$
and a convoluter $\beta$ for the same number of points $n$, a new local monodromy vector 
$\kappa (\beta ,\overvect{g})$. This is defined concretely as follows.

Define the {\em defect} $\delta (\beta , \overvect{g})$, which is going to be the difference between
the rank of the original local system, and the rank of 
the new local system obtained by middle convolution. Write out the coefficients
$$
g_i = \sum _{\alpha} m_i(\alpha )\cdot [\alpha ],
$$
where for clarity we denote by $[a]$ the point $a\in \LL$ considered as a divisor.
The defect is defined as
$$
\delta (\beta , \overvect{g}):= (n-2)r - \sum _{i=1}^n m (\beta ^{H_i,-1}).
$$
Corollary \ref{mhdimension}, applied to the divisor $K:= Q \cup \{ c\}$ with $n+1$ points, says that
\begin{equation}
\label{defectrank}
{\rm dim} MH^1(\Gamma , \rho ^{\beta}) = r + \delta (\beta , \overvect{g}).
\end{equation}
If no term $\beta$ is specified, it means to choose any $\beta$ such that 
$\beta ^{H_i, -1}:= (\beta ^{H_i})^{-1}\in \LL$ is a point
of maximal multiplicity for $g_i$, the resulting $\delta (\overvect{g})$ obviously 
doesn't depend on which choice is made. This is the same formula as considered in \S \ref{dimcount}.

Define the {\em local Katz transformation} at the point $q_i$ by
$$
\kappa _i(\beta ,\overvect{g}) := (m_i(\beta ^{H_i,-1})  + \delta (\beta , \overvect{g}) ) \cdot [\beta ^{V_i}] + 
\sum _{\alpha  \beta ^{H_i}\neq 1}
m_i(\alpha )\cdot [\alpha  \beta ^{U_i}].
$$
The {\em global Katz transformation} is defined by
$$
\overvect{\kappa}(\beta ,\overvect{g}):=
\left( \kappa _1(\beta ,\overvect{g}),\ldots , \kappa _n(\beta ,\overvect{g}) \right) .
$$

\begin{scholium}
\label{bettischolium}
Suppose $\rho$ is a representation of rank $r$ on $Y-Q_Y$ satisfying Convention \ref{semisimple}
that the local monodromy transformations are semisimple. Suppose $\beta$ is a convoluter, a rank one local system
on $X-J$. Assume that $\beta $ satisfies Convention \ref{chinontriv} that $\chi \neq 1$, and that
Convention \ref{chirhobeta} holds: $\chi \rho ^{\beta}(\alpha _k)$ have no trivial eigenvectors.

Let $\overvect{g}\in Div(\ggg _m)^n$ denote the vector of local monodromy data for $\rho$, and define
the defect $\delta (\beta , \overvect{g})$ as above.

Under these conditions, the middle coconvolution $MC^{\ast}_{\beta}(\rho )$ and the middle convolution
$MC_{\beta}(\rho )$ are local systems on $Z-Q_Z\cong Y-Q_Y$ of rank 
$$
r' = r + \delta (\beta , \overvect{g}),
$$
whose local monodromy transformations are semisimple and have local monodromy types 
given by the Katz transformation
$$
\overvect{\kappa} (\beta , \overvect{g}).
$$
\end{scholium}
{\em Proof:}
We have done the computations for the middle coconvolution in the previous subsection. 
The same is true for the middle convolution
by Poincar\'e-Verdier duality. The change in ranks is formula (\ref{defectrank}) above, which 
makes the defect appear in the multiplicity of the new eigenvalue as described after Proposition
\ref{mcdecomposition}.
\eop

\subsection{The Katz morphism on Betti moduli spaces}
\label{katzmorbetti}

This construction extends to giving a morphism on the level of moduli spaces: 

\begin{theorem}
\label{bettimorphism}
Let $M_{B}(\pp ^1,Q; \overvect{g})$ denote the Betti moduli space of local systems
on $\pp ^1-Q$ having semisimple local monodromy transformations corresponding to $\overvect{g}$. 
Suppose $\beta$ is a rank one local system on $(Z\times Y)-D$. 
Suppose that $(\beta ,\overvect{g})$ satisfy Conventions \ref{chinontriv} and \ref{chirhobeta}. Then the
middle convolution construction $L\mapsto MC_{\beta }(L)$ gives a morphism of moduli spaces
$$
MC_{\beta }: M_{B}(\pp ^1,Q; \overvect{g}) \rightarrow M_{B}(\pp ^1,Q; \overvect{\kappa} (\beta , \overvect{g})).
$$
\end{theorem}

This is sort of obvious, although technically speaking it requires some work: we should carry out
the middle convolution construction in the context of local systems of modules over a ring.
The fact that the $H^0$ and $H^2$ terms vanish, so the dimension of $H^1$ never jumps, is the
basic thing which makes it work. Notice that our Conventions \ref{chinontriv}
and \ref{chirhobeta} are only conditions on $\beta , \overvect{g}$, in particular they don't require
defining open subsets of the moduli spaces. 

This type of morphism between moduli spaces 
was considered in \cite{CrawleyBoeveyShaw} and other places. It is clearly related to the theory of
representations of the braid group such as the Burau representation, see \cite{Lawrence} \cite{Marin}.

\subsection{Involutivity}
\label{involut}

One of the main properties of Katz's construction is its involutivity. This implies that it gives an 
isomorphism of moduli spaces. The involutivity is basic to the constitution of an algorithm: one can go
forward to see if a local system with transformed local monodromy data
should exist, and if one is found then one can go backward to give back
a local system with the original local monodromy data.

Katz shows associativity of the convolution operator which allows him to deduce involutivity 
\cite[2.9.7]{Katz}. Later proofs were also given in the algebraic setting by V\"{o}lklein, Dettweiler-Reiter,
and Crawley-Boevey and Shaw.

Katz's proof didn't rely on the Fourier transform interpretation, which nevertheless furnishes a 
conceptual reason for involutivity: convolution can be interpreted as 
a composition of two Fourier transform operators
using also tensor products with rank one systems.
The Fourier transform is involutive by analogy with classical real analysis, so its composition two times 
and also with tensoring by an invertible rank one system, is involutive with an appropriate change of
convoluter as described below.  

It would be interesting to use connections with irregular singularities, and ``wild'' harmonic theory, to 
make this argument precise in the complex geometric setting. This would involve Bloch-Esnault \cite{BlochEsnault},
Sabbah \cite{SabbahIrregular} and Szabo \cite{Szabo}. See also \cite{Boalch2} and the recent preprint \cite{Hien}. 
Very recently Aker and Szabo have contructed an involutive Nahm transform for parabolic Higgs bundles \cite{AkerSzabo} 
which should lead to a complex analytic version of the Fourier transform construction.

For the middle convolution operation, involutivity can already be seen on the
level of local monodromy data. 

\begin{proposition}
\label{involutivedata}
Let $c:X\rightarrow X$ be the automorphism which flips the factors and let $\beta ^{\ast}$ be the
dual local system whose monodromy transformations are the inverses. Then 
$$
\overvect{\kappa} (c^{\ast}\beta ^{\ast}, \overvect{\kappa} (\beta , \overvect{g})) = \overvect{g}.
$$
\end{proposition}
{\em Proof:}
We will be making changes of variables in the sums, so it is convenient to have the following formula
for the Katz transformation in terms of $g_i = \sum _{\alpha}m_i(\alpha )[\alpha ]$ and the defect
$d:= \delta (\beta , \overvect{g})$:
$$
\kappa _i (\beta , \overvect{g}) =  
(m_i(\beta ^{H_i, -1}) + d)[\beta ^{V_i}] - 
m_i(\beta ^{H_i, -1})[\beta ^{V_i}\beta ^T] 
+ \sum _{\alpha}m_i(\alpha )[\alpha \beta ^{U_i}].
$$

Put $\gamma := c^{\ast}\beta^{\ast}$. In particular we have
$\gamma ^{H_i} = \beta ^{V_i,-1}$,  
$\gamma ^{V_i} = \beta ^{H_i,-1}$,  
$\gamma ^{U_i} = \beta ^{U_i,-1}$, and 
$\gamma ^{T} =  \beta ^{T, -1}$.
Write $\overvect{g'} = (g'_1,\ldots , g'_n) := \overvect{\kappa} (\beta , \overvect{g})$ and let
$r'$ be the rank, $m'(\alpha )$ the multiplicities and $d'$ the defect with respect to $\gamma$. 
One calculates that $d' = -d$
so the defects cancel out and at least on the level of ranks we have
${\rm rk}(\overvect{\kappa} (\gamma , \overvect{g}')) = r$.
We can write
$$
\kappa _i(\gamma , \overvect{g}') = 
(m'_i(\gamma ^{H_i, -1}) + d')[\gamma ^{V_i}] 
- (m'_i(\gamma ^{H_i, -1}))[\gamma ^{V_i}\gamma ^T] 
+ \sum _{\alpha '} m'_i(\alpha ') 
[\alpha '\gamma ^{U_i}]
$$
$$
= m_i(\beta ^{H_i, -1}) [\beta ^{H_i, -1}] 
- (m_i(\beta ^{H_i, -1}) + d)[\beta ^{V_i}\beta ^{U_i,-1}] 
+ \sum _{\alpha}m'_i(\alpha ')[\alpha \beta ^{U_i, -1}].
$$
The sum in the last term amounts to looking at $g'_i$ but translated by $\beta ^{U_i, -1}$, in other words
$$
\sum _{\alpha}m'_i(\alpha ')[\alpha \beta ^{U_i, -1}] =
$$
$$
(m_i(\beta ^{H_i, -1}) + d)[\beta ^{V_i}\beta ^{U_i, -1}] 
- m_i(\beta ^{H_i, -1})[\beta ^{V_i}\beta ^T\beta ^{U_i, -1}] 
+ \sum _{\alpha}m_i(\alpha )[\alpha \beta ^{U_i}\beta ^{U_i, -1}].
$$
After some textual cancellation, our full expression becomes 
$$
\kappa _i(\gamma , \overvect{g}') = 
m_i(\beta ^{H_i, -1}) [\beta ^{H_i, -1}] 
- m_i(\beta ^{H_i, -1})[\beta ^{V_i}\beta ^T\beta ^{U_i, -1}] 
+ \sum _{\alpha}m_i(\alpha )[\alpha ] = g_i.
$$
This completes the proof. 
\eop

Katz has also shown by direct calculation that the virtual dimensions of the moduli spaces 
for $\overvect{g}$ and $\overvect{\kappa}(\beta ,\overvect{g})$ are the same. 

Finally, we state the involutivity of the middle convolution morphism itself.
We have seen the involutivity on the level 
of local monodromy data, so it makes sense to look at the composition of the middle convolution morphisms.

\begin{theorem}
\label{involutivity}
The composition 
$$
M_{DR}(P,Q; \overvect{g}) 
\stackrel{MC(\beta )}{\longrightarrow} 
M_{DR}(P,Q; \kappa (\beta , \overvect{g}))
\stackrel{MC(c^{\ast}\beta ^{\ast})}{\longrightarrow} 
M_{DR}(P,Q; \overvect{g}) 
$$
is the identity, if we are in the situation of Theorem \ref{bettimorphism} for both of the morphisms. 
\end{theorem}

We don't describe the proof here but refer to Katz \cite{Katz}, V\"olklein \cite{VolkleinBraid},
Dettweiler-Reiter \cite{DettweilerReiterAlgKatz}, and more
recently Crawley-Boevey and Shaw \cite{CrawleyBoeveyShaw}.

\subsection{Detecting emptiness of the moduli space}
\label{detectingemptiness}

One of the main features of Katz's construction is that it permits us to detect whether a given moduli
space is empty or not in terms of the next moduli space in the algorithm. In other words,
$$
M_B(\overvect{g}) = \emptyset \Leftrightarrow M_B(\kappa (\beta , \overvect{g})) = \emptyset 
$$
assuming Conventions \ref{chinontriv} and \ref{chirhobeta}.

This is specially the case when $M_B(\kappa (\beta , \overvect{g}))$ is not even defined because
one of the divisors in the vector $\kappa (\beta , \overvect{g})$ is no longer effective. It is comforting
to work this case out explicitly. Let $d=\delta (\beta ,\overvect{g})$ denote the defect. The multiplicities
in the local divisors $\kappa _i(\beta , d,g_i)$ are either the same as in $g_i$, or else
they are changed by adding $d$. In particular, if $d\geq 0$ then we will never get to a noneffective divisor.
Thus we may assume that $d<0$. Suppose that $\kappa _i(\beta , d,g_i)$ becomes noneffective.
The only multiplicity which changes is $m_i (\beta ^{H_i,-1})$, which becomes
$$
m_i (\beta ^{H_i, -1}) + d.
$$
In particular, we are in the current situation, only if 
$$
m_i (\beta ^{H_i, -1}) + d < 0.
$$
Plugging in the formula for the defect, we have
$$
m_i (\beta _{H_i}^{-1}) + r(n-2) - \sum _{j=1}^n m_j(\beta ^{H_i, -1})  < 0,
$$
and adding $r$ to both sides and simplifying we get
$$
\sum _{j\neq i}(r - m_j(\beta ^{H_i, -1}))  < r.
$$
This says that 
the sum for $j\neq i$ of the ranks of the matrices $\rho ^{\beta }(\alpha _j) - 1$ is $<r$. 
Since these matrices generate the action of the group algebra on the vector space $V$, under this condition
the action cannot be irreducible. So, $\rho ^{\beta}$ and hence $\rho$ is not irreducible.
Thus, we have the following lemma.

\begin{lemma}
\label{describeempty}
Suppose that $\overvect{g}$ consists of effective divisors, and at least one of the divisors in
$\kappa (\beta , \overvect{g})$ is not effective. In this case, the representation
$\rho$ cannot be irreducible. In the case where $\overvect{g}$ is automatically irreducible, this means that
the moduli space $M_B(\overvect{g})$ is empty. 
\end{lemma}
\eop

\subsection{Running Katz's algorithm (Kostov's program)}
\label{kostovprog}

Kostov invented the protocol of applying Katz's algorithm to the nonrigid case.
Suppose $\overvect{g}$ is a local monodromy vector. Choose a convoluter $\beta$ so that
$(\beta ^{H_i}) ^{-1}$ is an eigenvalue of maximal multiplicity for $g_i$. Thus
$\delta (\beta , \overvect{g})=\delta (\overvect{g})$. If $\delta (\overvect{g}) < 0$ and
if the pair $(\beta , \overvect{g})$ satisfies Conventions \ref{chinontriv} and \ref{chirhobeta},
then we obtain an isomorphism of moduli spaces for $\overvect{g}$ and the Katz-transformed
vector $\overvect{\kappa}(\beta , \overvect{g})$. The rank strictly decreases,
so we can keep going on in the same way, until we get to $r=1$ or more generally to
a case where all of the local monodromy matrices are diagonal; to an impossibility result; 
to the problem discussed in the subsequent paragraph below;
or until we get into the range $\delta \geq 0$. If we hit an impossibility result anywhere
along the way, then
the original moduli space was empty. If we hit $r=1$ then the original moduli space was a point. 
If we get into the range $\delta \geq 0$ then according to Kostov 
we expect that the moduli space should be
nonempty, with a direct construction of some points \cite{KostovSteklov} \cite{KostovControl}. 
Crawley-Boevey and Shaw \cite{CrawleyBoeveyShaw} gave a different construction covering cases not treated
in \cite{KostovSteklov} \cite{KostovControl}, and 
prove
in some cases that the moduli space is a complete intersection. We will discuss a Higgs-bundle version
of the direct construction in \S \ref{explicitcon} below.

The problem with the previous paragraph is that somewhere along the way, we might hit
a vector $\overvect{g}$ for which every choice of $\beta$ corresponding to maximal multiplicities,
dissatisfies either Convention \ref{chinontriv} or Convention \ref{chirhobeta}.
In this case the algorithm no longer makes sense as we have described it. 
Apparently it can be made to work anyway, 
but this goes beyond the scope of the present discussion
and we refer to the papers of Kostov and Crawley-Boevey. 
Instead, we will just point out that
it doesn't happen if the original eigenvalues are sufficiently general. 

In Kostov's notation, a ``nongenericity relation'' is a subset of the eigenvalues counted with multiplicities,
of the same rank $r'\in 1,\ldots , r-1$ at each point $q_i$, such that the product of them all is $1$. 
Any nontrivial sub-local system has monodromy sub-data which give a nongenericity relation. 

Kostov says that a monodromy data vector $\overvect{g}$ is {\em 1-generic} if there is no nongenericity
relation of rank $r'=1$. This is the same as saying that there is no equation 
$a_1 \cdots a_n = 1$ such that $a_i$ is an eigenvalue of $g_i$. 

\begin{lemma}
\label{1gencons}
Suppose $\overvect{g}$ is $1$-generic, and suppose $\beta$ is a convoluter such that
each $(\beta ^{H_i})^{-1}$ is an eigenvalue of $g_i$. Then the pair $(\beta , \overvect{g})$
satisfies Conventions \ref{chinontriv} and \ref{chirhobeta} and we get a middle convolution isomorphism
between moduli spaces.
\end{lemma}
{\em Proof:}
It is trivial that the pair satisfies the conditions. In order to get an isomorphism we also need
to have the same conditions for the inverse 
pair $(c^{\ast}\beta ^{\ast}, \overvect{\kappa}(\beta , \overvect{g}))$. Convention \ref{chinontriv}
for $c^{\ast}\beta ^{\ast}$ is equivalent to Convention \ref{chinontriv} for $\beta$. For Convention 
\ref{chirhobeta} note that the eigenvalues of $\kappa _i(\beta , \overvect{g})$
are either $\varphi = \beta ^{U_i}\alpha$ for eigenvalues $\alpha$ of $g_i$ with $\beta ^{H_i}\alpha \neq 1$,
or else $\varphi = \beta ^{V_i}$. Convention \ref{chirhobeta} for the inverse pair thus requires 
for these $\varphi$
$$
(c^{\ast}\beta ^{\ast})^T (c^{\ast}\beta ^{\ast})^{H_i}\varphi \neq 1.
$$
Recalling that 
$(c^{\ast}\beta ^{\ast})^T (c^{\ast}\beta ^{\ast})^{H_i} = \beta ^{T,-1}\beta ^{V_i,-1}$, the condition becomes
$$
\beta ^{T,-1}\beta ^{V_i,-1}\cdot \beta ^{U_i}\alpha \neq 1, \;\;\; \mbox{for} \;\beta ^{H_i}\alpha \neq 1,
$$
$$
\beta ^{T,-1}\beta ^{V_i,-1}\cdot \beta ^{V_i} \neq 1.
$$
The first is verified by tautology and the second is Convention \ref{chinontriv}. 
\eop

In Katz's original rigid case, a nongenericity relation among eigenvalues of highest multiplicity
automatically causes the local system to become reducible, and meeting such a nongenericity relation anywhere
along the way rules out existence of any irreducible rigid local system. I would like to thank the referee for pointing
out the following very interesting 
example, which shows that there can be a nongenericity relation among other eigenvalues, even for an irreducible
rigid local system. The example consists of a local system of rank $3$ with $3$ singular points having local monodromy
eigenvalues $(a,b,c), \, (u,v,w), \, (g,h,h)$. It is rigid, and exists even with a nongenericity relation of the
form $aug=1$. If there is no other nongenericity relation then the local system cannot be reducible (by looking at the
block of size $2$). One can construct this system by convolution of a hypergeometric system $(a',b'),(u',v'), (g',h')$
with a convoluter having $\beta ^{H_i} = \beta ^{V_i} = x,y,\, \mbox{or} \, z = (h')^{-1}\, (i=1,2,\, \mbox{or}\, 3)$.
As an exercise in applying the Katz transformation, the convoluted system is 
$$
\left( \begin{array}{c} a'xy^{-1}z^{-1} \\ b'xy^{-1}z^{-1} \\ x \end{array} \right) , \; 
\left( \begin{array}{c} u'yx^{-1}z^{-1} \\ v'yx^{-1}z^{-1} \\ y \end{array} \right) , \; 
\left( \begin{array}{c} g'zx^{-1}y^{-1} \\ z \\ z \end{array} \right) .
$$
Thus $aug = a'xy^{-1}z^{-1} u'yx^{-1}z^{-1}g'zx^{-1}y^{-1} = a'u'g'x^{-1}y^{-1}z^{-1}$ can be equal to $1$ by an appropriate choice of 
$x,y,z$.

For nonrigid
local systems the situation is even less clear and we will be happy with the following result.

\begin{proposition}
\label{runalgorithm}
Fix Kostov's polymultiplicity   vector           (PMV)   \cite{KostovCRAS} etc.\
containing  the multiplicities   of eigenvalues  in   the divisors $g_i$.
The variety of all possible  $\overvect{g}$  with   this polymultiplicity vector, is a disjoint union of
$d$ connected open subsets  of tori, where   $d$ is  the pgcd of all the multiplicities in 
$\overvect{g}$. If $\overvect{g}$ is a  sufficiently   general point in any of these connected components, then
we can run Katz's algorithm   until we  hit  either an   empty moduli space for the reason discussed in 
\S \ref{detectingemptiness},   or
the case of all diagonal local    monodromy  (i.e. rank   one system tensored with $\cc ^r$), 
or the case $\delta \geq 0$ which   will be  discussed  in \S \ref{explicitcon} below. 
The monodromy vectors encountered along the  way are   always themselves general points, 
in particular they are $1$-generic.
\end{proposition}
{\em Proof:}
Invertibility of the transformation on local monodromy data (Lemma \ref{involutivedata}) plus its continuity
with respect to the eigenvalue parameters if the PMV is fixed,
imply that for $\overvect{g}$ general in its
connected component, and $\beta$ general in the variety of possible choices given that the $\beta ^{H_i}$ come
from $\overvect{g}$ (that is, general among the possible choices of $\beta ^{V_i}$), the resulting
$\overvect{\kappa }(\beta , \overvect{g})$ is again general in its connected component. Thus,
formally applying a sequence of Katz transformations as for the algorithm, we encounter only general 
local monodromy vectors. 

If the PMV is not simple, that is if the pgcd of all the multiplicities is $d\geq 2$,
then there can exist a nongenericity relation even for general $\overvect{g}$. 
However, the nongenericity relation is always of rank at least $r/d$, and the case $d=r$ is the 
degenerate one with only diagonal matrices. Thus, a
general $\overvect{g}$ in any connected component is
always $1$-generic, except in the degenerate diagonal case.
\eop

In case of a non-simple PMV, the variety
in the previous proposition has some components where there is a nongenericity relation. 
If the moduli space has dimension $2$, when we get to $\delta = 0$ and $\sigma = 0$ Kostov
shows in \cite{KostovKappa} that all local systems are reducible for the nongeneric components. 
The case of dimension $2$ is somewhat special and is not covered by our construction in \S \ref{explicitcon}. 

Roberts studies the geographical implications of Katz's algorithm in the rigid case \cite{Roberts}, and
it would be good to extend his results to the nonrigid case.

\section{Middle convolution---the de Rham version}
\label{drmiddleconv}

The de Rham version involves replacing local systems by
logarithmic connections \cite{Nitsure} \cite{BlochEsnault} \cite{Inaba}.
Middle convolution in the ``Fuchsian'' case of connections on the trivial bundle has been
extensively considered \cite{KostovCRAS} etc., \cite{HaraokaYokoyama} \cite{HaraokaYokoyama2} 
\cite{DettweilerReiterMC} \cite{DettweilerReiterPreprint}
\cite{Gleizer} \cite{CrawleyBoeveyAdditive}.
In our treatment we don't distinguish between trivial and nontrivial underlying bundles, so in a certain
sense we consider less information than these references, on the other hand our approach places things
in an abstract setting. 

In order to use the logarithmic de Rham complex, it is essential to have a morphism between smooth projective
varieties with normal crossings divisors. Thus we use the blowing-up $X$ with its divisor $J\subset X$ described 
in \S \ref{diagonalconfig}.  
The second projection gives a map $\xi : (X,J)\rightarrow (Z,Q)$ in good position, meaning that
the inverse image of $Q$ is the divisor $U + V \subset J$ which has normal crossings. 

For  a vector bundle with logarithmic connection $(E,\nabla )$ on $(Y,Q_Y)$ and a de Rham convoluter
$\beta \in H^0(X,\Omega ^1_X(\log J))$, define 
a vector bundle with logarithmic connection on $X$:
$$
(F,\nabla _F):= \eta ^{\ast}(E,\nabla )\otimes (\Oo _X,d+\beta ).
$$
The divisor $HT:=  H_1+\ldots + H_n+T \subset J$ is transverse to the fibers of $\xi$. In a relative version
of the discussion of \S , \ref{middledrcoh} we can define the middle relative
de Rham complex with respect to $HT$, by the exact sequence
\begin{equation}
\label{mdrdef}
0\rightarrow MDR(X/Z,F;HT)\rightarrow DR(X/Y,F)\rightarrow \Ff ^0_{HT/Z}[-1]\rightarrow 0.
\end{equation}
For $z\in Z$, denote by $X_z$ the fiber of $\xi$ over $z$, and let 
$$
MDR(X_z,F; HT_z):= MDR(X/Z,F;HT)|_{X_z}
$$
with similar notation for the full de Rham complex. 
Over points $z\in Z-Q$ this is the same thing as the middle de Rham complex for $(X_z,HT_z)\cong (Y, Q+\{ z\} )$ 
considered in \S \ref{middledrcoh}.  
In order to have a good base-change theory, we impose the following.

\begin{convention}
\label{goodbasechange}
For every $z\in Z$, the degree $0$ and $2$ hypercohomology groups of the restriction $MDR(X_z,F; HT_z)$ vanish.
\end{convention}

This condition implies that $\rr ^1\xi _{\ast}MDR(X/Z, F; HT)$ is locally free over $Z$ with fiber over
a point $z$ equal to $\hh ^1MDR(X_z,F; HT_z)$. It has a logarithmic Gauss-Manin connection denoted by
$\nabla _{GM,{\rm mid}}$, and we define the {\em de Rham middle convolution} as
$$
MC_{\beta}(E,\nabla ):=(\rr ^1\xi _{\ast}MDR(X/Z, F; HT),  \nabla _{GM,{\rm mid}}),
$$
a vector bundle with logarithmic connection on $(Z,Q_Z)$. 

The restriction of the quotient term in (\ref{mdrdef}) to a point $z\in Z$ is just a skyscraper sheaf
placed in cohomological degree $1$, so it has no $\hh ^0$ or $\hh ^2$. 
The long exact sequence for the higher derived direct image of the exact sequence (\ref{mdrdef}) 
therefore gives the following. 

\begin{lemma}
\label{mdrshortexact}
Suppose that Convention \ref{goodbasechange} holds. Then the same vanishing holds for the
full de Rham complex, the $\rr ^1\xi _{\ast}DR(X/Z,F)$ is again a vector bundle compatible with base change,
and we have a short exact sequence
\begin{equation}
0\rightarrow 
\rr ^1f_{\ast}MDR(X/Z,F;HT) \rightarrow \rr ^1f_{\ast} DR(X/Z,F)\rightarrow 
\rr ^0f_{\ast}(\Ff ^0_{HT/Z})\rightarrow 0.
\end{equation}
This short exact sequence is compatible with the Gauss-Manin connections $\nabla_{GM,{\rm mid}}$ on the
left and $\nabla _{GM}$ in the middle. 
\end{lemma}
\eop

The classical definition of the Gauss-Manin connection is as the connecting map for the short exact sequence
of complexes
\begin{equation}
\label{gmshortexact}
0\rightarrow DR(X/Z,F)\otimes \xi ^{\ast}\Omega ^1_Z(\log Q)[-1]
\rightarrow DR(X,F)\rightarrow DR(X/Z,F)\rightarrow 0.
\end{equation}
When $q\in Q$ is a singular point, the de Rham complex $DR(X,F)$, which by convention means the
logarithmic de Rham complex with respect to $J$, can be restricted to a complex 
$DR(X,F)|_{X_q}$ on the fiber $X_q\subset J$. We obtain a restriction of (\ref{gmshortexact})
to $X_q$.
Note that $\Omega ^1_Z(\log Q)_q\cong \cc$ 
and the residue of $\nabla _{GM}$ at $q$ is the endomorphism 
$$
\hh ^1DR(X_q,F)\rightarrow
\hh ^2(DR(X_q,F)\otimes _{\cc}\Omega ^1_Z(\log Q)_q[-1]) = \hh ^1DR(X_q,F)
$$
induced by the connecting map for the restriction of (\ref{gmshortexact}).

The expression as a connecting map is not very convenient for calculating the eigenvalues. 
The calculation was done by Katz in \cite{KatzICM} (thanks to H. Esnault for pointing out this reference). 
Without going through all of the details, here
is the conclusion. In our case, $q= q_i$ for some $i=1,\ldots ,n$, and
the singular fiber $X_q$ consists of two components 
$X_q = U_i\cup V_i$ meeting in a point $w_i := U_i\cap V_i$.
We have a short exact sequence
\begin{equation}
\label{components}
0 \rightarrow DR(U_i, F_{U_i}(-w_i)) \rightarrow DR(X_q, F_{X_q}) \rightarrow DR(V_i, F_{V_i}) \rightarrow 0.
\end{equation}
Note that $HT$ meets $X_q$ in a collection of smooth points distinct from the crossing point $w_i$. 
Thus the exact sequence defining the middle de Rham complex is compatible with (\ref{components}), and
we have the same short exact sequence for middle de Rham complexes 
\begin{equation}
\label{mdrcomponents}
0 \rightarrow MDR(U_i, F_{U_i}(-w_i), HT_{U_i}) 
\rightarrow MDR(X_q, F_{X_q}, HT_{X_q}) \rightarrow MDR(V_i, F_{V_i}, HT_{V_i}) \rightarrow 0.
\end{equation}

We refine Convention \ref{goodbasechange} to apply to each of the components:

\begin{convention}
\label{goodbasechange2}
For $j=0,2$ we require that 
$$
\hh ^j MDR(U_i, F_{U_i}(-w_i), HT_{U_i})=0, \;\;\;
\hh ^j MDR(V_i, F_{V_i}, HT_{V_i}) =0. 
$$
\end{convention}

Assuming Convention \ref{diagresnontriv}, this condition for all the $q_i$
implies Convention \ref{goodbasechange}. For points $z\in Z-Q$, Lemma \ref{nontrivialresiduecase} 
provides the required vanishing. 

\begin{proposition}
Assuming Convention \ref{goodbasechange2}, we get a short exact sequence from (\ref{mdrcomponents}) on the 
level of $\hh ^1$. 
The residues of $\nabla _F$ along $U_i$ and $V_i$ give endomorphisms of $F_{U_i}$ and $F_{V_i}$. 
These fit into a diagram
$$
\begin{array}{ccccccccc}
0 & \rightarrow & \hh ^1 MDR(U_i, F_{U_i}(-w_i)) & \rightarrow & \hh ^1 MDR(X_q, F_{X_q}) &
\rightarrow & \hh ^1 MDR(V_i, F_{V_i}) & \rightarrow & 0 
\\
 & & \downarrow & & \downarrow & & \downarrow & & 
\\
0 & \rightarrow & \hh ^1 MDR(U_i, F_{U_i}(-w_i)) & \rightarrow & \hh ^1 MDR(X_q, F_{X_q}) &
\rightarrow & \hh ^1 MDR(V_i, F_{V_i}) & \rightarrow & 0 
\end{array}
$$
where the endomorphism of $\hh ^1DR(X_q, F_{X_q})$ is the residue of the middle Gauss-Manin
connection $\nabla _{GM, {\rm mid}}$ at $q=q_i$. For brevity the notations $HT_{U_i}$ etc. have
been omitted.
\end{proposition}

In our case,
the endomorphisms of the left and right terms will be semisimple. This will imply that the residue 
of $\nabla _{GM, {\rm mid}}$ is
semisimple, as long as we know that the endomorphisms on the left and right don't have any common eigenvalues.
We can state this as the following lemma. 

\begin{lemma}
\label{residuenablagm}
Let $\psi _{U_i}$ and $\psi _{V_i}$ be the endomorphisms of 
$\hh ^1 MDR(U_i, F_{U_i}(-w_i),HT_{U_i})$ and $\hh ^1 MDR(V_i, F_{V_i},HT_{V_i})$ respectively, 
determined by the
endomorphisms of $F_{U_i}$ and $F_{V_i}$ given by the residues of $\nabla$ along $U_i$ and $V_i$. 
Suppose that these endomorphisms are semisimple, and don't have any common eigenvalues.
Then the residue of $\nabla _{GM,{\rm mid}}$ at $q_i$ is semisimple 
and isomorphic to $\psi _{U_i} \oplus \psi _{V_i}$.
\end{lemma}
\eop

On $U_i$ and $V_i$ we have a logarithmic structure also at the point $w_i$.
However, this point is not included in the ``middle'' part which is just $HT_{U_i}$ or $HT_{V_i}$.

We now turn to the fact that $F$ is the pullback of $(E,\nabla )$, tensored with $(\Oo _X,d+\beta )$.
From the above discussion, the main problem is to calculate the restrictions
$$
\eta ^{\ast}(E, \nabla )|_{U_i},\;\;\; 
\eta ^{\ast}(E, \nabla )|_{V_i}, \;\;\; \beta |_{U_i}, \;\;\; \beta |_{V_i}.
$$
We can define the restriction of a logarithmic one-form $\beta$ to $V_i$ as follows.
It depends on the pullback of the coordinate function $t$ from $Z$ (where $t(q_i)=0$). 
Set $b:= {\rm res}(\beta , V_i)$,
then $\beta - b \frac{dt}{t}$ is a logarithmic form having zero residue along $V_i$, thus it is
in the kernel of the residue map on logarithmic forms which maps by restriction to $\Omega ^1_{V_i}(\log )$.
Define
$$
\beta |^t_{V_i} := (\beta - b \frac{dt}{t}) |_{V_i}.
$$
It is a logarithmic form on $V_i$ whose residues along $H_j\cap V_i$ 
are just $\beta ^{H_j}$, for $j\neq i$. This determines the restriction, and it has residue
at the intersection point 
$$
{\rm res}(\beta |^t_{V_i}, w_i) = -\sum _{j\neq i}\beta ^{H_j}.
$$
The similarly-defined restriction $\beta |^t_{U_i}$ is a logarithmic form on $U_i$ whose residues
along $H_i\cap U_i$ and $T\cap U_i$ are respectively $\beta ^{H_i}$ and $\beta ^T$, so
$$
{\rm res}(\beta |^t_{U_i}, w_i) = -\beta ^{H_i} - \beta ^T.
$$
The relation $\beta ^T + \sum _i \beta ^{H_i} = 0$ gives
${\rm res}(\beta |^t_{V_i}, U_i\cap V_i) +{\rm res}(\beta |^t_{U_i}, U_i\cap V_i)  = 0$, 
characteristic of the fact that these restrictions correspond to a single logarithmic form
on $X_{q_i}= U_i\cup V_i$.

Now restrict $(F,\nabla _F)$ to $V_i$. Since $\eta |_{V_i}$ is the identity, 
$\eta ^{\ast}(E,\nabla )|_{V_i}\cong (E,\nabla )$.
The restriction of $F$ is therefore
$$
(F,\nabla _F)|_{V_i}=
(\eta ^{\ast}E,\eta ^{\ast}\nabla + \beta )|_{V_i} = 
(E,\nabla + \beta |^t_{V_i}).
$$
The residual endomorphism induced by $\nabla _F$ is just scalar multiplication by $\beta ^{V_i}$.

Next look at the restriction of $(F,\nabla _F)$ to $U_i$. 
It clearly depends only on the local form of $(E,\nabla )$ near the point $q_i$. We may even localize in an
analytic neighborhood, and so assume that $E$ has the form of a trivial bundle $\Oo ^r$ and
the connection is given by $\nabla = d + A\frac{dy}{y}$. We use the notation $y$ for our 
coordinate on $Y$ at the point $q_i$ (which should be the same as $t$ under $Y\cong Z$), 
also considered as a function on $Z\times Y$ or $X$. 

Now $t$ and $y$ give coordinates on $Z\times Y$. The ratio $u=y/t$ is a coordinate on $X$, 
in a neighborhood of the point $H_i\cap U_i$ along $U_i$. On $U_i$ it corresponds to the linear 
coordinate which takes the values
$u(U_i\cap H_i) = 0$,
$u(U_i\cap T) = 1$,
$u(U_i\cap V_i) = \infty $.  

The relation
$\frac{dy}{y} =\frac{du}{u} + \frac{dt}{t}$
allows us to calculate the restriction
$$
(\frac{dy}{y} )|^t_{U_i} = \frac{du}{u}.
$$
The residue of the pullback of $\frac{dy}{y}$ along $U_i$ is equal to $1$. 
Using $(E,\nabla ) \cong (\Oo ^r, d+A\frac{dy}{y})$ we get that the restriction of the pullback to
$U_i$ is
$$
(\eta ^{\ast}E,\eta ^{\ast}\nabla + \beta )|_{U_i} = 
(\Oo _{U_i}^r, d + A \frac{du}{u} + \beta |^t_{U_i}),
$$
and the residue of $\eta ^{\ast}\nabla + \beta $ along $U_i$ is $A+ \beta ^{U_i}$ which is an 
endomorphism of this bundle preserving the logarithmic connection. Here, in canonical terms
$\Oo ^r$ corresponds to the fiber $E_{q_i}$ and $A$ corresponds to the residue of $\nabla $ at $q_i$.

Apply Lemma \ref{residuenablagm} to calculate the residue of $\nabla _{GM, {\rm mid}}$. 
By Convention \ref{semisimpleresidues}, $A$ is semisimple with eigenvalues
never differing by a nonzero integer. 
Invoking either Corollary \ref{nontrivialresiduecase} using the fact that the residue of $\beta |^t_{U_i}$ at $T_{U_i}$
is a nonzero scalar (Convention \ref{diagresnontriv}), or just by direct computation, we have
$$
\hh ^iDR(\Oo _{U_i}(-U_i\cap V_i)^r, d + A \frac{du}{u} + \beta |^t_{U_i}) = 0,\;\;\; i=0,2.
$$
A direct computation using the fact that $\Omega ^1_{U_i}(\log )(-w_i)\cong \Oo _{U_i}$ gives
$$
\hh ^1DR(\Oo _{U_i}(-U_i\cap V_i)^r, d + A \frac{du}{u} + \beta |^t_{U_i}) 
\cong H^0 (\Omega ^1_{U_i}(\log )(-w_i)\cong \Oo _{U_i}) ^r
= \cc ^r,
$$
and the action of the residue of $\eta ^{\ast}\nabla + \beta $ is given by the matrix 
$A+ \beta ^{U_i}$.

The residue is
nontrivial on the diagonal (Convention \ref{diagresnontriv}) so the middle condition at $U_i\cap T$ has no effect, 
and the middle condition at the point $U_i\cap H_i$ removes the zero eigenspace of the residue there, 
that is to say the zero eigenspace of the matrix
$A+\beta ^{H_i}$. Introduce the following notation: if $M$ is a semisimple matrix then $M^{\dag}$ is
the same endomorphism but only of the sum of eigenspaces different from zero. Thus the contribution from 
$U$ to the residue of $\nabla _{GM,{\rm mid}}$ on the middle direct image is
$$
{\rm res}(\nabla _{GM, {\rm mid}})_U = (A+ \beta ^{H_i})^{\dag} + (\beta ^{U_i}- \beta ^{H_i}).
$$
Recall that $\beta ^{U_i} = \beta ^{H_i} + \beta ^T + \beta ^{V_i}$, giving
$$
{\rm res}(\nabla _{GM, {\rm mid}})_U = (A+ \beta ^{H_i})^{\dag} + \beta ^{V_i}+ \beta ^{T}.
$$

The contribution from $V$ is the cohomology of $\nabla + \beta |^t_{V_i}$, with middle condition at the
points $H_j\cap V_i$ for $j\neq i$ and no middle condition at $w_i = U_i\cap V_i$. This contribution occurs
with a single eigenvalue which is the residue of $\beta$, in our notation $\beta ^{V_i}$.
Let $d_i$ denote the dimension of this cohomology group and ${\bf 1}^{d_i}$ is the identity matrix of rank $d_i$. 
If Lemma \ref{residuenablagm}  can be applied then we conclude that
the full residue of the Gauss-Manin connection on the middle convolution is given by 
\begin{equation}
\label{gmresult}
{\rm res}(\nabla _{GM, {\rm mid}}) = \left[ (A+ \beta ^{H_i})^{\dag} + \beta ^{V_i}+ \beta ^{T}\right]
\oplus \left[ \beta ^{V_i} {\bf 1}^{d_i}\right] .
\end{equation}

In order to be able to apply Lemma \ref{residuenablagm} we need to know that the eigenvalues of the two pieces
don't coincide. 
We also need something for the first sentence in Corollary \ref{nontrivialresiduecase}. 
The following condition is analogous to Convention \ref{chirhobeta} from the Betti case.

\begin{convention}
\label{alphabetabeta}
For any eigenvalue $\alpha$ of the residue $A={\rm res}(\nabla , q_i)$, we have
$$
\alpha + \beta ^{H_i} + \beta ^T \not \in \zz , \;\;\; \alpha + \beta ^{H_i}\not \in \zz - \{ 0\} .
$$
\end{convention}

The first condition is equivalent to saying $\alpha - \sum _{j\neq i} \beta ^{H_j} \not \in \zz$, and 
if each $-\beta ^{H_j}$ is an eigenvalue of the residue at $q_j$ then this condition would be a  
consequence of 1-genericity as in \ref{1gencons}. The second condition will hold whenever we need to
choose $-\beta ^{H_i}$ from among the eigenvalues of ${\rm res}(\nabla , q_i)$,
by Convention \ref{semisimpleresidues}.  

\begin{lemma}
Suppose Conventions \ref{semisimpleresidues}, \ref{diagresnontriv} and \ref{alphabetabeta}
hold. Then the eigenvalues of $(A+ \beta ^{H_i})^{\dag}+ \beta ^{V_i}+ \beta ^{T}$ are distinct from $\beta ^{V_i}$, so
Lemma \ref{residuenablagm} can be applied as above (\ref{gmresult}) with
$$
d_i = (n-2)r - \sum _{j\neq i}m_j(-\beta ^{H_j}).
$$
\end{lemma}
{\em Proof:}
The eigenvalues of $(A+ \beta ^{H_i})^{\dag}+ \beta ^{V_i}+ \beta ^{T}$ are of the form 
$\alpha + \beta ^{H_i}+ \beta ^{T} +\beta ^{V_i}$
for $\alpha$ eigenvalues of $A$. Under Convention \ref{alphabetabeta} these are different from $\beta ^{V_i}$.
For the dimension of the piece coming from $V_i$, note that  
$$
{\rm res}(\nabla + \beta |^t_{V_i}, w_i) = 
A + {\rm res}(\beta |^t_{V_i}, w_i) = A -\sum _{j\neq i}\beta ^{H_j} = A + \beta ^{H_i} + \beta ^T.
$$ 
Convention \ref{alphabetabeta} says that the eigenvalues here are never integers, 
also the residues at points of $H_{V_i}$ are never nonzero integers. Thus
Corollary \ref{nontrivialresiduecase} applies and we can calculate the dimension $d_i$
by using the Euler characteristic which gives the formula as stated. The terms in the sum over $i\neq j$ 
come from the middle conditions at the points $H_j\cap V_i$; there is no middle condition at the
remaining point $w_i$.  
\eop

As in \S \ref{nitsuredrmod} and \S \ref{kclocmon}
above, let $\overvect{g} = (g_1,\ldots , g_n)$ denote the residual data for $(E,\nabla )$
with $g_i\in Div (\aaa ^1)$ effective divisors representing the multiplicity vectors of
the eigenvalues. Asking that the residues lie in conjugacy classes ${\bf c}(g_i)$ 
insures Convention \ref{semisimpleresidues} automatically, and Conventions 
\ref{diagresnontriv} and \ref{alphabetabeta} are conditions only on the pair $(\beta , \overvect{g})$. 
The result of Lemma \ref{residuenablagm} applied as in (\ref{gmresult}) says exactly that the vector of 
residual data for $\nabla _{GM,{\rm mid}}$ is given by $\overvect{\kappa}(\beta , 
\overvect{g})$. We can sum up as follows. 

\begin{scholium}
\label{drscholium}
Suppose $(E,\nabla )$ is a logarithmic connection on $(Y,Q_Y)$ with semisimple residues
(Convention \ref{semisimpleresidues}) corresponding to a vector
$\overvect{g}\in Div(\aaa ^1)^n$, and suppose $\beta$ is a de Rham convoluter (\S \ref{convdr}). 
Suppose that Convention \ref{diagresnontriv}) that $\beta ^T\not \in \zz$, and
Convention \ref{alphabetabeta} that $\alpha + \beta ^{H_i} + \beta ^T \not \in \zz$ and 
$\alpha + \beta ^{H_i}\not \in \zz - \{ 0\}$ for any eigenvalue
$\alpha$ of ${\rm res}(\nabla , q_i)$. These are conditions on $(\beta , \overvect{g})$ only.
Then the de Rham middle convolution $MC_{\beta}(E,\nabla )$ is a logarithmic connection on 
$(Z,Q_Z)$ with semisimple residues whose vector of residual data is given by the Katz transformation
$\overvect{\kappa}(\beta , \overvect{g})$. 
\end{scholium}

\begin{theorem}
\label{drkatziso}
Suppose $Q\subset \pp ^1$ is a set of $n$ points, $\overvect{g}\in Div(\aaa ^1)^n$
is a vector of semisimple residual data, and $\beta$ is a de Rham convoluter. Suppose
Conventions \ref{diagresnontriv} and \ref{alphabetabeta} hold. Then middle convolution induces
an isomorphism 
$$
MC_{\beta} : M_{DR}(\pp ^1, Q; \overvect{g}) \stackrel{\cong}{\rightarrow}
M_{DR}(\pp ^1, Q; \overvect{\kappa}(\beta , \overvect{g}).
$$
This isomorphism is involutive like in \S \ref{involut} (but with $-c^{\ast}\beta$ instead of
$c^{\ast}\beta ^{\ast}$) and is compatible with the isomorphism of Theorem \ref{bettimorphism} 
via the Riemann-Hilbert 
correspondence \cite[Theorem 1.2]{DettweilerReiterMC}. 
\end{theorem}

\section{Harmonic bundles and parabolic structures}
\label{harmonicbundles}

There is a notion of {\em parabolic bundle on $Z$ with parabolic structures at the $q_i$}. We don't repeat the
definition here. These will be called ``parabolic bundles'' for short.
If $E$ is a parabolic bundle then for each $q_i$ and each $\alpha \in \rr$ we have an associated graded
vector space $Gr _{\alpha , q_i}(E)$ which is finite-dimensional. Multiplication by a local coordinate at $q_i$ gives
an isomorphism 
$$
Gr _{\alpha , q_i}(E) \cong Gr _{\alpha +1, q_i}(E).
$$
Define the {\em residue} of $E$ at $p_i$ to be the associated-graded direct sum
$$
{\rm res}(E; q_i):= \bigoplus _{0\leq \alpha < 1} Gr _{\alpha , q_i}(E).
$$
For a fixed $\lambda \in \cc$ there is a notion of {\em logarithmic $\lambda$-connection} $\nabla$ on a parabolic
bundle $E$. 
The logarithmic structure is
with respect to the divisor $Q = q_1 + \ldots + q_n$. 
For any $\alpha _1,\ldots , \alpha _n$ 
it induces a $\lambda$-connection
$$
\nabla : E_{\alpha _1,\ldots , \alpha _n}\rightarrow E_{\alpha _1,\ldots , \alpha _n}\otimes \Omega ^1_Z(\log Q).
$$
With the canonical isomorphism 
$\Omega ^1_Z(\log Q)_{q_i}\cong \cc $, a $\lambda$-connection on the parabolic bundle $E$ induces an 
endomorphism called the {\em residue}
$$
{\rm res}(\nabla ; q_i): {\rm res}(E;q_i) \rightarrow {\rm res}(E;q_i).
$$
The residue of the pair $(E,\nabla )$ at a point $q_i$ is defined as the pair of an $S^1$-graded vector space with 
endomorphism
$$
{\rm res}(E, \nabla ; q_i):= \left(  {\rm res}(E ; q_i), {\rm res}(\nabla ; q_i) \right) .
$$
When necessary, we introduce an index to denote the piece ${\rm res}(\nabla ; q_i)_{\alpha}$
acting on $Gr _{\alpha , q_i}(E)\subset {\rm res}(E ; q_i)$. 

If $F\subset E$ is a parabolic subbundle compatible with $\nabla$ then we can consider its {\em parabolic degree},
and the {\em parabolic slope} is the parabolic degree divided by the rank.
We say that $(E, \nabla )$ is {\em stable} (resp. {\em semistable}) if for any strict parabolic subbundle 
compatible with $\nabla$, the parabolic slope of $F$ is strictly less than (resp. less than or equal to) that of $E$.

Given a parabolic bundle $E$, define its {\em parabolic type at $q_i$} to be the divisor on $S^1$
$$
{\rm type}(E,q_i) := \sum _{0\leq \alpha < 1}(\dim Gr_{\alpha ,  q_i}(E))\cdot [\alpha ] .
$$
This follows the discussion in \S \ref{kclocmon} with $\LL = S^1$. 

Up to isomorphism, the residue at $q_i$ of a parabolic logarithmic $\lambda$-connection
${\rm res}(E, \nabla ; q_i)$ is classified by specifying its type $\sum g_{\alpha}[\alpha ]$ and
for each $\alpha$, specifying the Jordan normal form of an endomorphism of a vector space of dimension $g_{\alpha}$.

We say that the residues of $(E,\nabla )$ are {\em semisimple} if the Jordan normal forms are diagonalizable.
This means that on each $Gr_{\alpha , q_i}$  the action of ${\rm res}(\nabla ; q_i)_{\alpha}$
is semisimple or diagonalizable, so it corresponds to a divisor on $\aaa ^1$.
Altogether, when the residues are semisimple, the isomorphism class of the residue of $(E,\nabla )$ 
at $q_i$ is determined by a divisor of total degree $r$ on $S^1\times \aaa ^1$.

Recall that we have a notion of {\em tame harmonic bundle} on $Z-Q$. 
A harmonic bundle consists of a flat connection, and an equivariant harmonic map. The tameness condition means
that locally near the singularities, the the flat sections of the connection on sectors have polynomial 
growth with respect to the harmonic metric. Measuring the growth rate of sections leads to parabolic structures.
The flat connection decomposes as 
$$
D = D' + D'' = (\partial + \overline{\theta}) + (\delbar + \theta )
$$
where $\partial + \delbar$ is a unitary connection, $\theta + \overline{\theta}$ antipreserves the metric,
and $\theta$ is holomorphic. 
Fix $\lambda \in \cc$ which allows us to define a differential $\lambda$-connection
$$
D_{\lambda} := \lambda D' + D''.
$$
The $(0,1)$ piece, which has contributions from both terms $\lambda D'$ and $D''$, 
is a usual holomorphic structure giving rise to a 
holomorphic vector bundle $E$. The $(1,0)$ piece is a holomorphic $\lambda$-connection on $E$.

Measuring the growth rate of sections in a holomorphic frame, leads to an extension of $E$ as a parabolic 
bundle over $Z$ with parabolic structure along $Q$, again denoted $E$. The connection $\nabla = (D_{\lambda})^{1,0}$
is a logarithmic $\lambda$-connection on the parabolic bundle $E$. The parabolic logarithmic $\lambda$-connection
$(E,\nabla )$ is {\em polystable}, in other words
a direct sum of stable objects of the same slope. 

Conversely, given a parabolic logarithmic $\lambda$-connection $(E,\nabla )$, if it is polystable
then there exists an essentially unique structure
of harmonic bundle given by a harmonic metric on $E$ over $X-Q$ with the appropriate growth rates determined by
the parabolic structure. The metric connection is unique and the metric is unique up to a positive real scalar on
each stable piece. The Higgs case is $\lambda = 0$. 

In keeping with Convention \ref{semisimple}, we would like to insure that the monodromy transformations are semisimple.
We furthermore assume that the residues of $\theta$ vanish. This amounts to restricting to representations
where the local monodromy eigenvalues are in $S^1\subset \ggg _m$ together with trivial filtered local system
structures at the singularities. 

\begin{convention}
\label{parrestriv}
The residue of the Higgs field $\theta$ on the associated-graded of the parabolic bundle 
$Gr_{\alpha , q_i}(E)$ at any point 
$q_i\in Q$ is equal to zero.
\end{convention}

\begin{lemma}
If $(E,\theta )$ is a polystable
parabolic Higgs bundle of degree zero
satisfying Convention \ref{parrestriv} then the monodromy transformations of the corresponding local system are
semisimple, with eigenvalues in $S^1$ corresponding to $e^{2\pi i\alpha }$ for $\alpha$ the parabolic weights.
\end{lemma}
{\em Proof:}
This follows from the local considerations shown in \cite{hbnc}.
\eop

The moduli space of parabolic logarithmic $\lambda$-connections can be constructed, see
many references on parabolic bundles included in the bibliography below. 
This moduli space becomes isomorphic (as a real analytic space possibly
with singularities) to a space of harmonic bundles which can be constructed as in Hitchin's original 
case \cite{Hitchin}, see also 
Fujiki \cite{Fujiki}. As $\lambda$ varies we get a family over the affine line, which is the
{\em nonabelian Hodge filtration} on the moduli space. Glueing with the complex conjugate we get the
twistor space for Hitchin's hyperk\"ahler structure \cite{Hitchin} \cite{Fujiki}. In our case of an open
curve, some further work is needed, see \cite{Nakajima} for example.
We should deal with the transformation of residual types which occurs when we change
$\lambda$ \cite{hbnc}, a situation which appears to reflect some kind of weight-two phenomenon corresponding
to the punctures. 
In the case of quasi-unipotent monodromy we should be able to deal with the problem
by looking at local systems with unipotent monodromy on a DM-curve. 
This general moduli problem will not be considered any more
here, but constitutes a good source of further questions: first and foremost we would like to have
Katz isomorphisms between these moduli spaces coming from a parabolic middle convolution. 

\subsection{Cyclotomic harmonic bundles}
\label{cycloharbun}

Our construction of Higgs bundles will be based on a trick to insure stability. Recall that $\cc ^{\ast}$
acts on the space of parabolic Higgs bundles by $t: (E,\theta )\mapsto (E,t\theta )$. If we assume that the
residue of $\theta$ is unipotent (or even equal to zero if we want to keep with Convention \ref{semisimple})
then this action preserves the residue of the parabolic structure and $\theta$, so by \cite{hbnc} it
preserves the local type of monodromy transformations of the corresponding local system. 

Recall that a {\em complex variation of Hodge structure} is a harmonic bundle which is a fixed point for the
full action of $\cc ^{\ast}$ or equivalently for the action of $S^1$. A variant is to look at the action of a
finite cyclic subgroup of roots of unity $\mu _m \subset \cc ^{\ast}$. Recall \cite{Hitchin} \cite{canonical}
that the action of $S^1$ preserves the harmonic metric structure.  A harmonic bundle which is a fixed point
of the action of $\mu _m$ will be called a {\em cyclotomic harmonic bundle}, where $m\geq 2$ is considered as
fixed for now---later we can say ``$m$-cyclotomic'' if we need to
specify $m$. 

When we say that $(E,\theta )$ is a fixed point this means that it is provided with an additional structure
of an action of $\mu _m$ on $E$ such that for any $u\in \mu _m$ and $e\in E$ we have
$(u\theta )(ue) = u(\theta e)$.

The structure of a cyclotomic harmonic bundle is very similar to the structure of a complex variation of 
Hodge structure. The group of characters of $\mu _m$ is canonically isomorphic to $\zz / m\zz$ because we have
defined $\mu _m$ as coming from a privileged embedding $\mu _m\subset \cc ^{\ast}$. The action of $\mu _m$ on $E$
provides a decomposition according to characters
$$
E = \bigoplus _{p\in \zz /m\zz } E^p,
$$
and the formula $(u\theta )(ue) = u(\theta e)$ then says
$$
\theta : E^p \rightarrow E^{p+1}\otimes \Omega ^1_Z(\log Q).
$$
The only difference with the case of variations of Hodge structure is that $p+1$ is taken in
the quotient group $\zz /m\zz$ so $\theta$ includes a piece of the form
$$
\theta : E^m \rightarrow E^{1}\otimes \Omega ^1_Z(\log Q).
$$
Thus, $\theta$ is no longer required to be a nilpotent transformation. 
In the parabolic case, the decomposition is compatible with the parabolic structure. 

The action of $\mu _m\subset S^1$ preserves the differential operators of the harmonic bundle,
so by averaging we can always choose a $\mu _m$-invariant harmonic metric. The decomposition of
$E$ is then orthogonal and the pieces $E^p$ are preserved by the metric connection $\partial + \delbar$.
The complex conjugate $\overline{\theta}$ goes from $E^p$ to $A^{0,1}(E^{p-1})$. Thus, in all respects this
looks like a complex variation of Hodge structure except that the Hodge decomposition is viewed circularly
and the Kodaira-Spencer components can go all the way around the circle.

Hitchin in \cite{HitchinTeichmuller} gave a construction of a subspace of Higgs bundles which 
corresponded to an analogue of Teichm\"uller space. In Hitchin's construction,
a basic variation of Hodge structure is modified by adding
a new term in the Higgs field. In this sense, the notion of cyclotomic Higgs field is
a variant on \cite{HitchinTeichmuller}. Hitchin's Teichm\"uller Higgs bundles can sometimes be
cyclotomic: in the notation of \cite{HitchinTeichmuller} if $\alpha _m\neq 0$ but 
$\alpha _{m-1} = \ldots = \alpha _2 = 0$ then the Higgs field written down there is cyclotomic.

The correspondence between Higgs bundles and local systems is compatible with the action of $\mu _m$,
and this helps with the stability condition. Suppose $(E,\theta )$ is an $m$-cyclotomic Higgs bundle, that
is a bundle with action of $\mu _m$ (or equivalently a decomposition as above) compatible with the action on $\theta$.
We say that it is {\em cyclotomically semistable} (resp. {\em cyclotomically stable}) if for any $\mu _m$-invariant
sub-Higgs bundle, the slope is smaller (resp. strictly smaller) than the slope of $E$. 

\begin{proposition}
\label{polystable}
Suppose $(E,\theta )$ is a cyclotomically stable parabolic cyclotomic Higgs bundle of parabolic degree zero.
Then $(E,\theta )$ is polystable as a regular parabolic Higgs bundle, 
and it has a $\mu _m$-equivariant harmonic metric with growth rates corresponding to the parabolic structure.
\end{proposition}
{\em Proof:}
The $\beta$-subbundle is $\mu_m$-invariant, so 
cyclotomic semistability implies semistability. The socle is $\mu_m$-invariant, so cyclotomic
polystability implies polystability. Then \cite{cvhs}, Theorem 1 which took into account the possibility
of having the action of a group such as $\mu _m$, provides an invariant harmonic metric. 
The growth rates are governed by the choice of initial metric, as discussed in \cite{cvhs} \S 10 and \cite{hbnc}.
\eop
 
It is interesting to note that the cyclotomic Higgs bundles play a special role in the compactification of
the de Rham moduli space. Recall that the compactification puts at infinity a divisor obtained by
dividing $M_{Higgs} - \{ \theta \,\,  \mbox{nilpotent} \}$ by the action of $\cc ^{\ast}$. If we take the quotient in the
sense of stacks, then the compactification becomes a DM stack and the stacky points with automorphism 
group $\mu _m$ are exactly the cyclotomic Higgs bundles.

\subsection{The maximal case}
\label{maxcase}

Traditionally one of the easiest cases is when the $E^p$ are line bundles. For this, take $m=r$ equal to 
the rank. 

\begin{lemma}
\label{maximal}
Suppose $E=\bigoplus E^p$ is an $r$-cyclotomic parabolic Higgs bundle of rank $r$ with $\theta$ not nilpotent. This means that
the $E^p$ are parabolic line bundles and
every component $\theta ^p : E^p \rightarrow E^{p+1}\otimes \Omega ^1_Z$ is nonzero. Then 
$(E,\theta )$ is cyclotomically stable, hence polystable as an ordinary parabolic Higgs bundle. 
\end{lemma}
{\em Proof:}
Non-nilpotence of $\theta$ requires that all the components be nonzero, in particular all of the bundles $E^p$ are
nonzero. Since their number is equal to the rank, they must be line bundles. 
Suppose $F\subset E$ is a $\mu _r$-invariant saturated subbundle. It decomposes as a direct sum of $F^p\subset E^p$.
If any one of the $F^p$ is nonzero then all of them are nonzero because every component  $\theta^p$ is
nonzero and they go around in a circle. If $F$ is saturated we get $F^p=E^p$. Thus, there are no
$\mu _r$-invariant subbundles of rank strictly between $0$ and $r$, so $E$ is vacuously cyclotomically stable.
By Proposition \ref{polystable}, $(E,\theta )$ is polystable in the ordinary sense.
\eop

The structure of a non-nilpotent $r$-cyclotomic parabolic Higgs bundle is particularly easy to understand.
It consists of a collection of parabolic line bundles $E^1,\ldots , E^r$, together with morphisms of parabolic bundles
$$
\theta ^p : E^p \rightarrow E^{p+1}\otimes \Omega ^1_Z(\log Q).
$$
Convention \ref{parrestriv} says that we want ${\rm res}(\theta ^p,q_i)$ to induce the zero map on the 
associated-graded $Gr_{\alpha , q_i}$ for any $\alpha$ and any singular point $q_i\in Q$. In view of the
fact that we are dealing with parabolic line bundles, there is only one weight on each side, and
the residue map is automatically zero unless the two weights are the same. This will appear in our
criterion below. 

\subsection{Explicit construction}
\label{explicitcon}

A parabolic line bundle on $(Y,Q)$ is always of the form 
$$
E = \Oo _Y(k)(a_1q_1 + \ldots + a_nq_n)
$$
with $k\in \zz$ and $a_i\in [0,1)$. Written in this way, the vector $[k; a_1,\ldots , a_n]$ is
uniquely determined and uniquely determines $E$, and we may use it as notation. 
Suppose we are given two parabolic line bundles $E= [k; a_1,\ldots , a_n]$ and 
$E' = [k'; a'_1,\ldots , a'_n]$. A map from $E$ to $E'$ will consist of a holomorphic map 
$$
f : \Oo  _Y(k) \rightarrow \Oo _Y(k'),
$$
such that if $a_i > a'_i$ then $f(q_i)=0$. In this description we use the fact that $a_i,a'_i\in [0,1)$,
in particular $|a_i-a'_i|<1$ always. 

In order to have a map which furthermore induces the zero map on associated graded spaces at each $q_i$, 
$$
Gr_{\alpha , q_i}(f) = 0: Gr_{\alpha , q_i}(E)\rightarrow Gr_{\alpha , q_i}(E'),
$$
we should require that $f(q_i)=0$ also when $a_i = a'_i$. Thus, the description
of these maps (which we call zero-residue maps) is that $f(q_i)=0$ when $a_i \geq a'_i$.

Finally, we get to a description of a map
$$
f: E \rightarrow E'\otimes \Omega ^1_Y(\log Q)
$$
inducing zero on the residues as required by Convention \ref{parrestriv}. 
Recall that $\Omega ^1_Y(\log Q) = \Oo _Y(n-2)$. Thus, such a map $f$
is the same thing as a holomorphic section 
$$
f\in \Gamma (Y, \Oo _Y(k'-k + n-2))
$$
such that $f(q_i)=0$ whenever $a_i \geq a'_i$, or we can also say
$$
f\in \Gamma (Y,\Oo _Y(k'-k + n-2 - \# \{ i,\; a_i\geq a'_i\} )).
$$
We get the following lemma. 

\begin{lemma}
\label{parmapcriterion}
Suppose $E= [k; a_1,\ldots , a_n]$ and 
$E' = [k'; a'_1,\ldots , a'_n]$ are parabolic line bundles on $(Y,Q)$. 
Then there exists a nontrivial zero-residue map
$$
f: E \rightarrow E'\otimes \Omega ^1_Y(\log Q)
$$
if and only if 
$$
\# \{ i,\; a_i\geq a'_i\} \leq k'-k + n-2.
$$
If equality holds then the map $f$ has no zeroes other than as required for the residues,
and is unique up to a scalar constant. If the inequality is strict then $f$ has other zeroes.
\end{lemma}
\eop

Suppose we are given a vector of divisors $\overvect{g} = (g_1,\ldots ,g_n)\in Div(S^1)n$,
with $deg(g_i)=r$. 
Choose the standard section of $\rr \rightarrow \rr /\zz $ to identify $S^1\cong [0,1)$. 
Write 
$$
g_i = \sum _{\alpha\in [0,1)} m_i(\alpha )[\alpha ].
$$

A sequence $a_{i,1},a_{i,2},\ldots ,a_{i,r} \in [0,1)$ is called an {\em arrangement of $g_i$}
if each $\alpha $ occurs in the sequence with multiplicity $m_i(\alpha )$. 
This is to say that the sequence of exponentials of the $a_{i,j}$
is a possible sequence of eigenvalues along the diagonal, for a matrix 
in the conjugacy class $C(g_i)$.

An arrangement is called {\em good} if the number of indices $t$ with $a_{i,t}\geq a_{i,t+1}$ is
minimal.  Let $T(g_i)$ be the minimal number of such indices, thus the arrangement is good if
$$
\# \{ t, \; a_{i,t}\geq a_{i,t+1} \} = T(g_i).
$$
Here, and always below, the indices are taken modulo $r$, for example if $t=r$ then $a_{i,t+1}= a_{i,1}$.
Enumerate the indices $t$ as above, in increasing order $t_{1} < \ldots  < t_{p}$. We can thus write our
arrangement as a ``sawtooth'':
$$
\begin{array}{c}
a_{i,1} < a_{i,2} < \ldots < a_{i,t_1} \\
a_{i,t_1} \geq a_{i,t_1+1} \\
a_{i,t_1+1} < \ldots < a_{i,t_2} \\
a_{i,t_2} \geq a_{i,t_2+1} \\
\vdots \\
a_{i,t_p} \geq a_{i,t_p+1} \\
a_{i,t_p+1} < \ldots < a_{i,r} 
\end{array}
$$
with $a_{i,r} < a_{i,1}$ unless $t_p$ happens to be $r$. 
Now let 
$$
g_{i,j}:= [a_{i,t_j+1}] + [a_{i,t_j+2}] + \ldots + [a_{i,t_{j+1}-1}] + [a_{i,t_{j+1}}],
$$
with $t_{j+1} := t_1$ when $j=p$ and the terms in $g_{i,p}$ adapted appropriately. 
These are reduced effective divisors, that is each eigenvalue
occurs with multiplicity at most $1$, because the sequences are strictly increasing in between the $t_j$.
And we have a decomposition 
$$
g_i = g_{i,1} + \ldots + g_{i,p}.
$$
Notice that 
$$
t_{j+1} = {\rm deg}(g_{i,j}) + t_{j} \;\; (\mbox{modulo} \; r).
$$
Conversely, given a pair of $p$-uples written $(t_1,\ldots , t_p; g_{i,1},\ldots , g_{i,p})$ 
with $t_j$ an increasing sequence in $1,\ldots ,r$ and the $g_{i,j}$ giving a decomposition of
$g_i$ into reduced effective divisors, we get an arrangement. The arrangement is good if $p$ is minimal and
equal to the maximal multiplicity in $g_i$. There is a one-to-one correspondence between such notations and
arrangements for $g_i$.

\begin{lemma}
\label{minimalmax}
The minimal number of $t$'s is equal to the maximum multiplicity in the divisor $g_i$,
$$
T(g_i)= \max _{\alpha}m_i(\alpha ).
$$
\end{lemma}
{\em Proof:}
It is easy to see that for any $\alpha$ we have $T(g_i)\geq m_i(\alpha )$. On the other hand, we can
clearly choose a decomposition into reduced effective 
divisors $g_i = g_{i,1}+ \ldots  + g_{i,p}$ with $p = \max_{\alpha}m_i(\alpha )$. Thus the minimal $p$ is
equal to the maximum of the $m_i(\alpha )$.
\eop

Now suppose that for each $i=1,\ldots , g$ 
we have chosen a good arrangement $a_{i,j}$ for $g_i$. 
Let 
$$
(t_{i,1},\ldots , t_{i,p_i}; g_{i,1},\ldots , g_{i,p_i})
$$
be the notation established above with $p_i=T(g_i)$.
For any sequence $k_1,\ldots , k_r$ define parabolic line bundles
$$
E^j := [k_j; a_{1,j},\ldots ,a_{n,j}].
$$
In order to construct a cyclotomic Higgs bundle $E=\bigoplus _{j=1}^r E^j$,
we investigate the possible choice of $k_1,\ldots , k_r$ such that there are nontrivial zero-residue maps 
$$
\theta ^j: E^j \rightarrow E^{j+1}\otimes \Omega ^1_Y(\log Q)
$$
including the case $j=r, j+1=1$. 
Let $\tau _j$ denote the cardinality
$$
\tau _j := \# \{ i , \;\; a_{i,j}\geq a_{i,j+1}\} .
$$
For any sequence of $k_j$ put
$$
z_j := k_{j+1} - (\tau _j + k_j + 2-n).
$$

\begin{lemma}
\label{zmeaning}
With the above notations, there exist zero-residue maps $\theta ^j$ if and only if
$z_j \geq 0$ for $j=1,\ldots , r$. In this case, $z_j$ is the number of extra zeros of $\theta ^j$
beyond what are required by the zero-residue condition. The $z_j$ are subject to the relation
\begin{equation}
\label{zrelation}
z_1 + \ldots + z_r = \delta (\overvect{g}),
\end{equation}
so there exists a possible choice of $z_j$ or equivalently of $k_j$ if and only if $\delta (\overvect{g})\geq 0$.
\end{lemma}
{\em Proof:} The first statements come from Lemma \ref{parmapcriterion}.
From the definition of $T(g_i)$ and $\tau _j$ we get
$$
\sum _{j=1}^r \tau _j = \sum _{i=1}^n T(g_i).
$$
Thus 
$$
z_1 + \ldots + z_r = r (n-2) - \sum _{i=1}^n T(g_i) = \delta (\overvect{g}).
$$
\eop

Given $z_j\geq 0$ subject to the relation $z_1+\ldots + z_r = \delta (\overvect{g})$, and given $k_1$,
we obtain the remaining $k_2,\ldots k_r$ from the formula for $z_j$. Construct the parabolic bundles $E^j$
and nontrivial zero-residue maps $\theta ^j$. This yields an $r$-cyclotomic Higgs bundle $(E,\theta )$.

\subsection{The degree}
\label{degree}

As pointed out in Lemma \ref{maximal}, the $(E,\theta )$ constructed this way is cyclotomically stable,
hence polystable as a Higgs bundle. To finish the construction of a local system we need to insure
that its parabolic degree vanishes. The parabolic degree of $E^j = [k_j; a_{1,j},\ldots ,a_{n,j}]$
is 
$$
{\rm deg}^{\rm par}(E^j) = k + a_{1,j} + \ldots + a_{n,j}.
$$
Adding up gives the parabolic degree of $E$: 
$$
{\rm deg}^{\rm par}(E)=\sum _{j=1}^r k_j + \sum _{i,j}a_{i,j} = \sum _{j=1}^r k_j + \sum _{i, \alpha }m_i(\alpha )\alpha .
$$
By induction,
$$
k_{j+1} = k_1 + z_1 + \ldots + z_j + \tau _1 + \ldots  + \tau _j + j(2-n).
$$
We have
$$
{\rm deg}^{\rm par}(E)= P+ k_1r + \sum _{j=1}^r (r-j)(z_j + \tau _j)
$$
where
$$
P:= \sum _{i, \alpha }m_i(\alpha )\alpha  + \sum _{j=1} ^r j(r-j)(2-n) 
$$
represents the piece which doesn't depend on the choice of
arrangements or of $k_j$. The condition $Det(\overvect{g})= 1$ says that
$P\in \zz$. 

Recall that $\tau_j$ is the number of $i$ such that $a_{i,j}\geq a_{i,j+1}$. This is the same as the number
of $i$ such that $j\in \{ t_{i,1},\ldots , t_{i,p_i}\}$. The terms involving $\tau _j$ can
be recast as a sum over
the elements $t_{i,j}$. We conclude the following formula for the parabolic degree of the Higgs bundle
we have constructed:
\begin{equation}
\label{concludeformula}
{\rm deg}^{\rm par}(E)= 
P + k_1r + \sum _{i=1}^n \sum _{j= 1}^{p_i} (r-t_{i,j}) +  \sum _{j=1}^r (r-j)z_j.
\end{equation}

\begin{theorem}
Suppose given a vector of divisors $\overvect{g}\in Div(S^1)^n$ such that $Det(\overvect{g})= 1$. 
Suppose that the defect is strictly positive,
$$
r(n-2)-\sum T(g_i) = \delta (\overvect{g}) > 0.
$$
Fix any collection of good arrangements $a_{i,j}$ for $\overvect{g}$. 
Then it is possible to choose the $k_j$ subject to the constraint
$$
k_{j+1} \geq \tau _j + k_j + 2-n,
$$
such that ${\rm deg}^{\rm par}(E) = 0$.
\end{theorem}
{\em Proof:}
Fixing the collection of good arrangements, the terms in (\ref{concludeformula}) involving $t_{i,j}$ are fixed. 
Because of the strictly positive defect,
there is a nontrivial choice of $z_1,\ldots , z_j$. Geometrically this means that we have a choice as to 
how many zeros $\theta _j$ can have.

Put $z_r := \delta (\overvect{g})-1$
and $z_j=0$ for all but one value of $j = j'$ in which case $z_{j'}= 1$. Make the convention here that
if $j'=r$ then  $z_r := \delta (\overvect{g})$ instead. With this choice we get
$$
\sum _{j=1}^r (r-j)z_j = r-j',
$$
and by choosing $j'$ appropriately this can take on any value between $0$ and $r-1$. In particular, modulo $r$
it can take on all values. By adjusting $k_1$ appropriately, ${\rm deg}^{\rm par}(E)$ can take on any integer
value.
\eop

\begin{theorem}
Suppose given a vector of divisors $\overvect{g}\in Div(S^1)^n$ such that $Det(\overvect{g})= 1$. 
Suppose that the defect is zero,
$$
r(n-2)-\sum T(g_i) = \delta (\overvect{g})= 0,
$$
but the superdefect is strictly positive $\sigma (\overvect{g}) > 0$. 
Then it is possible to choose a good arrangement $\{ a_{i,j}\}$ for $\overvect{g}$
and $k_1$, which determines the remaining $k_j$ in the zero-defect case by the constraint
$$
k_{j+1} = \tau _j + k_j + 2-n,
$$
such that ${\rm deg}^{\rm par}(E) = 0$.
\end{theorem}
{\em Proof:}
In the case where the defect is zero, we are constrained to have $z_j=0$. In particular,
once we fix $k_1$ then the others are determined. 
We have the simplified formula
$$
{\rm deg}^{\rm par}(E)= 
P + k_1r + \sum _{i=1}^n \sum _{j= 1}^{p_i} (r-t_{i,j}). 
$$
On the other hand, the fact that the superdefect is nonzero means that for some $i$ there is
at least one eigenvalue $\alpha '$ which appears with multiplicity $m_i(\alpha ') < p_i= \max _{\alpha }m_i(\alpha )$.
In particular, for any arrangement which we denote now generically by $A$, we have at least one interval 
not containing the eigenvalue $\alpha '$. Define the following operation on arrangements: find
an interval $t_{i,j}+1,\ldots , t_{i,j+1}$ containing $\alpha '$ but such that the preceding interval
$t_{i,j-1}+1,\ldots , t_{i,j}$ doesn't contain $\alpha '$. Move $\alpha '$ from the one to the other.
We get a new arrangement $\partial A$ with the property that all $t_{i',j'}(\partial A)$ are the same as for $A$,
except 
$$
t_{i,j}(\partial A) = t_{i,j}(A) + 1.
$$
Note that $\partial A$ will always be a good arrangement whenever $A$ is good. 
From this and the above formula we find
$$
{\rm deg}^{\rm par}(E (\partial A)) = {\rm deg}^{\rm par}(E (A))  -1.
$$
In particular, iterating the operation $A\mapsto \partial A$ and modifying $k_1$ 
we find that ${\rm deg}^{\rm par}(E (A))$
can take on all integer values as $A$ runs through all the good arrangements.
\eop

Putting together these two theorems we get:

\begin{corollary}
\label{endresult}
Suppose $\overvect{g}\in Div(S^1)^n$  is a vector of local monodromy data with eigenvalues in $S^1\subset \ggg _m$.
Suppose the defect is positive $\delta (\overvect{g}) \geq 0$. If $\delta = 0$ then suppose that
the superdefect is strictly positive; this is equivalent to supposing 
that the virtual dimension of the moduli space is
at least $4$. Then there exists a parabolic $r$-cyclotomic Higgs bundle $(E,\theta )$ of parabolic degree $0$,
cyclotomically stable and polystable in the usual sense,
corresponding to a local system with local monodromy data $\overvect{g}$. 
\end{corollary}

The only cases left to be treated are when the moduli space has dimension $2$. There are four families as
listed in Lemma \ref{dim2list}. These cases are considered by Kostov in \cite{KostovKappa}. As he notes there,
the determinant of the vector $\frac{\overvect{g}}{d}$ is a $d$-th root of unity. For $d>1$ and primitive
root of unity, it looks like there cannot be an $r$-cyclotomic Higgs bundle; however a solution exists 
\cite{KostovControl} and one might hope to construct
an $m$-cyclotomic Higgs bundle for smaller $m$ and with some component bundles $E^j$ of rank $2$. When the root
of unity is not primitive, Kostov shows that there are no irreducible solutions.

\section{Further questions}
\label{furtherquestions}

It would be good to have the full middle-convolution 
theory for the general setup of parabolic logarithmic $\lambda$-connections
\cite{TMochizuki} \cite{TMochizuki2} \cite{TMochizuki3} \cite{SabbahTwistor}. This raises some nontrivial 
questions such as defining the middle higher direct image in the parabolic setting, obtaining a base-change
result analogous to Convention \ref{goodbasechange2}, and showing polystability of the middle convolution. 
It was my original goal to treat these questions here but that turned out to be very difficult.

Aker and Szabo have communicated to me their recent preprint \cite{AkerSzabo} in which they do the Nahm transform
(essentially the same as Fourier transform) for parabolic Higgs bundles with irregular singularities
having poles of order $\leq 2$ at infinity. This should allow one to obtain
the middle convolution for parabolic Higgs bundles by following Katz's original method.

If the weights of a parabolic structure are rational, i.e. for every point $q_i$ the parabolic type is a divisor
concentrated over roots of unity in $S^1$, then
as discussed in \cite{Borne} \cite{Cadman} \cite{MatsukiOlsson} \cite{IyerSimpson}, 
the parabolic bundle may also be viewed as a bundle on a
Deligne-Mumford stack $Z[\frac{Q}{m}]$
obtained by assigning an integer $m$ to the points $q_i$. Here $m$ should be chosen 
to be divisible by all the denominators of the rational weights which occur. In this case, a logarithmic
$\lambda$-connection on the parabolic bundle may also be viewed as a logarithmic connection on the corresponding
DM-bundle.  An intermediate case between the non-parabolic case we have discussed in \S \ref{drmiddleconv}
and the general case of parabolic logarithmic $\lambda$-connections, would be the case of
parabolic logarithmic $\lambda$-connections with rational weights. Also assuming that the residue of the
connection on $Gr_{\alpha , q_i}(E)$ is the scalar $\alpha$, these objects would be equivalent to
local systems on the DM-stacks $Z[\frac{Q}{m}]$. It should be possible to have a theory of Katz's middle
convolution for these objects. The blown-up surface $(X,J)$ would be provided with a stack structure and
the singular fibers would be twisted curves \cite{AbramovichVistoli} \cite{AbramovichGraberVistoli}.

The moduli spaces have numerous additional structures. 

\begin{conjecture}
The isomorphisms between different moduli spaces given by the middle convolution map, preserve
the Hodge filtration, the $\cc ^{\ast}$ action, and 
the Hitchin hyperk\"ahler structure when this is defined, that is when
the eigenvalues of the local monodromy transformations are
in $S^1$. 
\end{conjecture}

This conjecture can probably be proven by Aker and Szabo with their method \cite{AkerSzabo}, indeed they show that
the Nahm transform preserves the hyperk\"ahler structure of the moduli spaces. 

V\"olklein points out in \cite{VolkleinBraid} 
that the Katz isomorphisms between various Betti moduli spaces are compatible with the
action of the braid group of the points $q_1,\ldots , q_n \in \pp ^1$. Similarly, the cohomological formulation
immediately implies that for $n\geq 4$ the Katz isomorphisms between different de Rham moduli spaces are compatible with
the nonabelian Gauss-Manin connection, i.e. the isomonodromic deformation equations. 
This was used by Boalch to get information about finite Painlev\'e orbits in \cite{Boalch1} \cite{Boalch2}.
It would be interesting to look further at the dynamics of the braid action 
and the isomonodromy equations.

\begin{question}
Which Hodge types can occur at variations of Hodge structure in the moduli spaces? 
How does the Hodge type change under middle convolution?
\end{question}

To what extent do we get unexpected or exceptional automorphisms of
moduli spaces, due to the possibility of running Katz's algorithm in several different ways?
In particular, one could start in the range $\delta \geq 0$, do a series of middle convolutions which go
out of this range, then another series to go back. In some cases this should change the local monodromy vector,
so we should obtain isomorphisms $M_B(\overvect{g}_1)\cong M_B(\overvect{g}_2)$ for $g_1\neq g_2$ 
in the range $\delta \geq 0$. Say that these two local monodromy vectors are {\em middle-convolution equivalent}
in this case. 

\begin{question}
What is the quotient of the set of local monodromy vectors with $\delta \geq 0$, by the relation
of middle-convolution equivalence? In each middle convolution equivalence class, does the operation
of going out and back again provide any nontrivial automorphisms of $M_B$? 
\end{question}

\begin{question}
Is there a Torelli theorem 
saying that the isomorphism class of $M_B(\overvect{g})$ and/or $M_{DR}(\overvect{g})$
possibly with additional structures such as the Hodge filtration, the hyperk\"ahler metric, etc., 
determines the middle convolution equivalence class of $\overvect{g}$ (and maybe the collection
of points $Q$ depending on how much structure we are considering)?
\end{question}

\begin{problem}
Generalize Roberts' observations on the geography of the Katz algorithm \cite{Roberts} to the nonrigid case.
\end{problem}

It would be good to
compare explicitly what is happening in our presentation, which basically follows Kostov's notation and setup, with
the notation and setup used by Crawley-Boevey. Note that in Crawley-Boevey's point of view, the Katz operations are root
reflections, and he uses several reflections in a row to get into a positive Weyl chamber before 
giving an explicit construction. This is obviously basically the same procedure as what we are doing here.
It would be good to compare the numbers, and
also to recover Roberts' results and observations \cite{Roberts} in the Crawley-Boevey formulation. 

What is the exact relationship between our de Rham version of the middle convolution, and the algebraic
operations on Fuchsian systems considered by Kostov, Haraoka-Yokoyama, Crawley-Boevey?

Theorem \ref{defobs} says that the middle cohomology
of $End(E)$ gives the deformation and obstruction theory for the moduli space of
representations with fixed conjugacy classes on a curve.
Remembering that the middle cohomology is really intersection cohomology, this suggests that we should ask
for the geometric interpretation of the intersection cohomology of $End(E)$ in the higher dimensional case.
More precisely, is there a natural derived moduli stack of local systems generalizing Kapranov \cite{Kapranov}
based on the intersection cohomology?  And, what kind of geometric objects does this derived moduli stack 
parametrize? 

\subsection{Low-dimensional cases}

One of the main reasons for looking at moduli spaces of representations on the punctured Riemann sphere
is that these give many more examples with small dimension, than are obtained from Hitchin's original 
case of compact Riemann surfaces. This was first pointed out by Hausel \cite{Hausel} with his ``toy example''.
We have constructed local systems whenever the virtual dimension 
is $\geq 4$. In some sense the first case to look at is dimension $2$, which has to be one of the cases listed in 
Lemma \ref{dim2list}. Unfortunately, our technique of construction broke down in this case, but we can hope
to have a variant. 

The explicit techniques applied by Gleizer in the rigid case \cite{Gleizer} should be applicable to 
low-dimensional cases. 

It would be interesting to compute as explicitly as possible all of the various structures and properties
for some concrete low-dimensional cases. For example, what does the $\cc ^{\ast}$ action on the moduli space
look like?
Some things to study in low-dimensional cases would be: 
compactifications \cite{Hausel} and their 
the dynamics
\cite{DaskalopolousDostoglouWentworth}
\cite{DaskalopolousDostoglouWentworth2}
\cite{DaskalopolousWentworth} \cite{BiswasMitraNag},
the Hitchin system \cite{Hitchin} \cite{HauselThaddeus},
the relationship with Painlev\'e equations \cite{HitchinPoincelet} \cite{Boalch1} \cite{Boalch2} \cite{DettweilerReiterPainleve},
jumps and wall-crossing phenomena such as in \cite{Thaddeus}
\cite{Nakajima}, 
real structures and Toledo invariants \cite{BradlowGarciaPradaGothen} \cite{BurgerIozziWienhard}
\cite{MarkmanXia} \cite{Xia}.


\begin{thebibliography}{MM}


\bibitem{AbramovichVistoli}
D. Abramovich, A. Vistoli.
Compactifying the space of stable maps.
{\em J. Amer. Math. Soc.} {\bf 15} (2002), 27-75.



\bibitem{AbramovichGraberVistoli} 
D. Abramovich, T. Graber, A. Vistoli.  
{\em Algebraic orbifold quantum products}, 
Orbifolds in mathematics and physics (Madison, WI, 2001), 1--24, Contemp. Math., 
\textbf{310}, Amer. Math. Soc., Providence, RI, 2002.

\bibitem{AkerSzabo}
K. Aker, S. Szabo. Algebraic Nahm transform for parabolic Higgs bundles on $\pp ^1$. Preprint {\tt math.AG/0610301}.

\bibitem{BalajiBiswasNagaraj}
V. Balaji, I. Biswas, D.  Nagaraj.
Principal bundles over projective manifolds with parabolic structure over a divisor. 
{\em Tohoku Math. J.} {\bf 53} (2001), 337-367.

	
\bibitem{Belkale}
P. Belkale.
Local systems on $\pp ^1-S$ for $S$ a finite set. 
{\em Compositio Math.} {\bf 129} (2001), 67-86. 


\bibitem{BenHamedGavrilov}
B. Ben Hamed, L. Gavrilov.
Families of Painleve VI equations having a common solution.
Preprint {\tt math.CA/0507002}.


\bibitem{Biquard}
O. Biquard. Sur les fibr\'es paraboliques sur une surface complexe.
{\em J. London Math. Soc.} {\bf  53} (1996), 302-316.


\bibitem{Biswas} 
I. Biswas. 
{\em Parabolic bundles as orbifold bundles}, Duke Math. J., {\bf 88} (1997), 305-325.


\bibitem{Biswas2}
I. Biswas. 
Flat connections on a punctured sphere and geodesic polygons in a Lie group. 
{\em J. Geom. Phys.} {\bf 39} (2001), 129-134.

\bibitem{Biswas3} 
I. Biswas. {\em Chern classes for parabolic bundles}, 
J. Math. Kyoto Univ.  37  (1997),  no. \textbf{4}, 597--613.

\bibitem{Biswas4}
I. Biswas. 
Stable principal bundles and reduction of structure group.
Preprint {\tt math.AG/0608569}. 

\bibitem{BiswasMitraNag}
I. Biswas, M. Mitra, S. Nag.
Thurston boundary of Teichmüller spaces and the commensurability modular group. 
{\em Conform. Geom. Dyn.} {\bf 3} (1999), 50-66.



\bibitem{BlochEsnault}
S. Bloch, H. Esnault.
Local Fourier transforms and rigidity for $\Dd$-modules. 
{\em Asian J. Math.} {\bf 8} (2004),  587-605. 

\bibitem{Boalch1}
P. Boalch. Painlev\'e equations and complex reflections. 
{\em Proceedings of the Conference in Honor of F. Pham (Nice, 2002)}. 
{\sc Ann. Inst. Fourier} 
{\bf 53} (2003), 1009-1022. 

\bibitem{Boalch2}
P. Boalch.
From Klein to Painlev\'e via Fourier, Laplace and Jimbo. 
{\em Proc. London Math. Soc.} {\bf 90} (2005), 167-208. 

\bibitem{BodenYokogawa}
H. Boden, K. Yokogawa. Moduli spaces of parabolic Higgs bundles and parabolic $K(D)$ pairs over smooth curves, I.
{\em Internat. J. Math.} {\bf 7} (1996), 573-598.

\bibitem{Bolibruch}
A. Bolibruch. The Riemann-Hilbert problem on a compact Riemannian surface.
{\em Proc. Steklov Inst. Math.} {\bf 238} (2002), 47-60.

\bibitem{BolibruchMalkeMitschi}
A. Bolibruch, S. Malek, C. Mitschi. 
On the generalized Riemann-Hilbert problem with irregular singularities.
{\em Expo. Math.} {\bf 24} (2006), 235-272. 

\bibitem{Borne} 
N. Borne. Fibr\'es paraboliques et champ des racines. Preprint {\tt math.AG/0604458}.

\bibitem{Bouw}
I. Bouw.
Reduction of the Hurwitz space of metacyclic covers. 
{\em Duke Math. J.} {\bf 121} (2004), 75-111. 

 	

\bibitem{Bradlow}
S. Bradlow.  Special metrics and stability for holomorphic bundles with global sections.
{\em J. Differential Geom.} {\bf 33} (1991), 169-214.

\bibitem{BradlowDaskalopolous}
S. Bradlow, G. Daskalopolous. Moduli of stable pairs for holomorphic bundles over Riemann surfaces,
II. {\em Int. J. Math} {\bf 4} (1993), 903-925. 

\bibitem{BradlowDaskalopolousWentworth}
S. Bradlow, G. Daskalopolous, R. Wentworth. Birational equivalence of vortex moduli. 
{\em Topology} {\bf 35} (1996), 731-748.

\bibitem{BradlowGarciaPradaGothen}
S. Bradlow, O. Garc\'ia-Prada, P. Gothen.
Surface group representations and $U(p,q)$-Higgs bundles. {\em J. Diff. Geom.} {\bf 64} (2003), 111-170.

\bibitem{Budur}
N. Budur. 
Unitary local systems, multiplier ideals, and polynomial periodicity of Hodge numbers. 
Preprint {\tt math.AG/0610382}. 

\bibitem{BurgerIozziWienhard}
M. Burger, A. Iozzi, A. Wienhard. 
Surface group representations with maximal Toledo invariant.
Preprint {\tt math.DG/0605656}. 

\bibitem{Cadman} 
C. Cadman. {\em Using stacks to impose tangency conditions on curves}, Preprint {\tt math.AG/0312349}.

\bibitem{Cailotto}
M. Cailotto.
{\em Algebraic connections on logarithmic schemes.}
C. R. Acad. Sci. Paris S\'er. I Math. \textbf{333} (2001),  1089--1094.


\bibitem{CohenOrlik}
D. Cohen, P. Orlik. Arrangements and local systems. {\em Math. Res. Lett.} {\bf 7} (2000), 299-316.

\bibitem{canonical}
K. Corlette.  Flat $G$-bundles with canonical metrics. {\em J. Diff. Geom.} {\bf 28} (1988), 361-382.

\bibitem{rigid}
K. Corlette. Rigid representations of K\"ahlerian fundamental groups. {\em J. Diff. Geom.} {\bf 33}
(1991), 239-252.

\bibitem{CorletteSimpson}
K. Corlette, C. Simpson. 
Preprint, in preparation.


\bibitem{CrawleyBoeveyAdditive}
W. Crawley-Boevey. On matrices in prescribed conjugacy classes with no common invariant subspace and sum zero. 
{\em Duke Math. J.} {\bf 118} (2003), 339-352.


\bibitem{CrawleyBoeveyIHES}
W. Crawley-Boevey.
Indecomposable parabolic bundles and the existence of matrices in prescribed conjugacy class 
closures with product equal to the identity. 
{\em Publ. Math. I.H.E.S.} {\bf 100} (2004), 171-207. 
 	
\bibitem{CrawleyBoeveyShaw}
W.  Crawley-Boevey, P. Shaw.
Multiplicative preprojective algebras, middle convolution and the Deligne-Simpson problem. 
{\em Adv. Math.} {\bf 201} (2006), 180-208.  	

\bibitem{CullerShalen}
M. Culler, P. Shalen.
Varieties of group representations and splittings of $3$-manifolds.
{\em Ann. of Math.} {\bf 117} (1983), 109-146.


 
\bibitem{DaskalopolousDostoglouWentworth}
G. Daskalopolous, S. Dostoglou, R. Wentworth.  On the Morgan-Shalen compactification of the $SL(2,\cc )$ character
varieties of surface groups. {\em Duke Math. J.} {\bf 101} (2000), 189-207. 



\bibitem{DaskalopolousDostoglouWentworth2}
G. Daskalopolous, S. Dostoglou, R. Wentworth. 
Character varieties and harmonic maps to $ R$-trees. {\em Math. Res. Lett.} {\bf 5} (1998), 523-533.

\bibitem{DaskalopolousWentworth}
G. Daskalopoulos, R. Wentworth. 
(i)\, The Yang-Mills flow near the boundary of Teichm\"uller space. {\em Math. Ann.} {\bf 318} (2000), 1-42.


\bibitem{DeligneRS} 
P. Deligne. {\em Equations diff\'erentielles a points singuliers
reguliers}. {\sc Lect. Notes in Math.} $\bf{163}$, 1970.


\bibitem{DeligneMostow}
P. Deligne, G. Mostow. 
{\em Commensurabilities among lattices in $PU(1,n)$}. {\sc Ann. of Math. Studies}
{\bf 132},  Princeton University Press (1993).

\bibitem{Dettweiler}
M. Dettweiler. 
Galois realizations of classical groups and the middle convolution.
Habilitationsschrift, Universit\"at Heidelberg (2005). 
Preprint {\tt math.NT/0605381}. 

\bibitem{DettweilerReiterOdd}
M. Dettweiler, S. Reiter. 
On rigid tuples in linear groups of odd dimension. 
{\em J. Algebra} {\bf 222} (1999), 550-560. 

 	
\bibitem{DettweilerReiterAlgKatz}
M. Dettweiler, S. Reiter.  
An algorithm of Katz and its application to the inverse Galois problem. 
{\em Algorithmic methods in Galois theory}. {\sc J. Symbolic Comput.} {\bf 30} (2000), 761-798. 

\bibitem{DettweilerReiterMC}
M. Dettweiler, S. Reiter. 
On the middle convolution. Preprint {\tt math.AG/0305311}.

\bibitem{DettweilerReiterPreprint}
M. Dettweiler, S. Reiter. 
Middle convolution of Fuchsian systems and the construction of rigid differential systems. Preprint (2004).

\bibitem{DettweilerReiterPainleve}
M. Dettweiler, S. Reiter.
Painlev\'e equations and the middle convolution.
Preprint {\tt math.AG/0605384}.

\bibitem{DettweilerReiterG2}
M. Dettweiler, S. Reiter. 
On exceptional rigid local systems.
Preprint {\tt math.AG/0609142}.

\bibitem{Donaldson}
S. Donaldson. A new proof of a theorem of Narasimhan and Seshadri. 
{\em J. Diff. Geom.} {\bf 18}  (1983), 269-277.

 	
\bibitem{Dwork}
B. Dwork.
On systems of ordinary differential equations with transcendental parameters. 
{\em J. Differential Equations} {\bf 156} (1999), 18-25. 

\bibitem{EsnaultHertling}
H. Esnault, C. Hertling. Semistable bundles on curves and reducible representations of the fundamental group.
{\em Int. J. Math.} {\bf 12} (2001), 847-855.

\bibitem{EsnaultViehweg} 
H. Esnault, E. Viehweg.  
Logarithmic De Rham complexes
and vanishing theorems, {\em Inventiones} {\bf 86} (1986), 161-194.

\bibitem{EsnaultViehweg2}
H. Esnault, E. Viehweg. 
Semistable bundles on curves and irreducible representations of the fundamental group. 
{\em Algebraic geometry: Hirzebruch 70 (Warsaw, 1998)} {\sc Contemp. Math.} {\bf 241}, AMS (1999), 129-138.

\bibitem{Fujiki}
A. Fujiki. 
Hyper-K\"ahler structure on the moduli space of flat bundles.
{\em Prospects in complex geometry (Katata and Kyoto, 1989)} 
Springer {\sc Lecture Notes in Math.} {\bf 1468}
(1991), 1--83.

\bibitem{Gleizer}
O. Gleizer.
Some explicit solutions of the additive Deligne-Simpson problem and their applications. 
{\em Adv. Math.} {\bf 178} (2003), 311-374. 

\bibitem{GoldmanXia}
W. Goldman, E. Xia. 
Rank One Higgs bundles and representations of fundamental groups of Riemann surfaces.
Preprint {\tt 
math.DG/0402429}.


\bibitem{Golyshev}
V. Golyshev.
Riemann-Roch variations. 
{\em Izv. Math.} {\bf 65} (2001), 853-881.

\bibitem{GorodentsevKuleshov}
A. Gorodentsev, S. Kuleshov. 
Helix theory. 
{\em Mosc. Math. J.} {\bf 4} (2004), 377-440, 535. 

\bibitem{Gothen}
P. Gothen. Components of spaces of representations and stable triples. 
{\em Topology} {\bf 40} (2001), 823-850.



\bibitem{GromovSchoen}
M. Gromov, R. Schoen. 
Harmonic maps into singular spaces and $p$-adic superrigidity for lattices in groups of rank one.
{\em I.H.E.S. Publ. Math.} {\bf 76} (1992), 165-246.


\bibitem{Haraoka}
Y. Haraoka.
Evaluation of Stokes multipliers for a certain system of differential equations 
corresponding to a rigid local system. 
{\em Funkcial. Ekvac.} {\bf 46} (2003), 187-211. 
 	
 	
\bibitem{HaraokaIntegral}
Y. Haraoka. 
Integral representations of solutions of differential equations free from accessory parameters. 
{\em Adv. Math.} {\bf 169} (2002), 187-240.

\bibitem{HaraokaYokoyama}
Y. Haraoka, T. Yokoyama. 
Construction of rigid local systems and integral representations of their sections. 
{\em Math. Nachr.} {\bf 279} (2006), 255-271. 


\bibitem{HaraokaYokoyama2}
Y. Haraoka, T. Yokoyama.
On rigidity of Pfaffian systems coming from Okubo systems.
{\em Kyushu J. Math.} {\bf 55} (2001), 189-205.


\bibitem{Hausel}
T. Hausel. Compactification of moduli of Higgs bundles. {\em J. Reine Angew. Math.} {\bf 503} (1998), 169-192.

\bibitem{HauselThaddeus}
T. Hausel, M. Thaddeus. 
Mirror symmetry, Langlands duality, and the Hitchin system. 
{\em Inventiones} {\bf 153} (2003), 197-229.

\bibitem{Hien}
M. Hien. Periods for irregular singular connections on surfaces.
Preprint {\tt math.AG/0609439}, to appear {\em Math. Ann.}. 


\bibitem{Hinich}
V. Hinich. 
Descent of Deligne groupoids.
{\em Internat. Math. Res. Notices} {\bf 1997}, 223-239.


\bibitem{Hirzebruch}
F. Hirzebruch. 
Arrangements of lines and algebraic surfaces.
{\em Arithmetic and geometry, Vol. II}, 
{\sc Progr. Math.} {\bf 36}
Birkh\"auser (1983), 113--140.

\bibitem{Hitchin}
N.  Hitchin. 
The self-duality equations on a Riemann surface. {\em Proc. London Math. Soc.} {\bf 55} (1987), 59-126.


\bibitem{HitchinPoincelet}
N.  Hitchin. 
Poncelet polygons and the Painlev\'e equations. 
{\em Geometry and analysis (Bombay, 1992)}, 
{\sc Tata Inst. Fund. Res.}, Bombay (1995), 151-185. 

\bibitem{HitchinTeichmuller}
N.  Hitchin. 
Lie groups and Teichm\"uller space.
{\em Topology} {\bf 31} (1992), 449-473.

\bibitem{Illusie} 
L. Illusie. 
{\em R\'eduction semi-stable et d\'ecomposition de complexes de de Rham \`a coefficients},  Duke Math. J., 
{\bf 60} (1990), 139-185.


\bibitem{IllusieKatoNakayama}
L. Illusie, K. Kato, C. Nakayama.
{\em Quasi-unipotent logarithmic Riemann-Hilbert correspondences}, J. Math. Sci. Univ. Tokyo, {\bf 12} (2005), 1-66.

\bibitem{Inaba} 
M. Inaba.
{\em Moduli of parabolic connections on a curve and Riemann-Hilbert correspondence},
Preprint \verb}math.AG/0602004}.


\bibitem{InabaIwasakiSaito}
M. Inaba, K. Iwasaki, Masa-Hiko Saito.
Moduli of Stable Parabolic Connections, Riemann-Hilbert correspondence and Geometry of Painlev\'e 
equation of type VI, Part I.
Preprint {\tt math.AG/0309342}.

\bibitem{IyerSimpson}
J. Iyer, C. Simpson.  
A relation between the parabolic Chern characters of the de Rham bundles.
Preprint {\tt math.AG/0603677}. 



\bibitem{JostYangZuo}
J.  Jost, Y. Yang, K. Zuo. The cohomology of a variation of polarized 
Hodge structures over a quasi-compact K\"ahler manifold.
Preprint {\tt math.AG/0312145}.

\bibitem{JostYangZuo2}
J.  Jost, Y. Yang, K. Zuo.
Cohomologies of unipotent harmonic bundles over quasi-projective varieties I: The case of noncompact curves.
Preprint {\tt math.AG/0505144}.

\bibitem{Kapranov}
M. Kapranov.
Injective resolutions of $BG$ and derived moduli spaces of local systems. 
{\em J. Pure Appl. Algebra} {\bf 155} (2001), 167-179.

\bibitem{FKato} 
F. Kato.
{\em The relative log Poincar\'e lemma and relative log de Rham theory}, Duke Math. J., {\bf 93} (1998), 179-206.

\bibitem{KatoNakayama}
K. Kato, C. Nakayama.
{\em Log Betti cohomology, log \'etale cohomology, and log de Rham cohomology of log schemes over $\cc$}.
Kodai Math. J. \textbf{22} (1999), 161--186.

\bibitem{KatzOda}
N. Katz, T. Oda.
{\em On the differentiation of de Rham cohomology classes with respect to parameters}.
J. Math. Kyoto Univ. {\bf 8}, 2 (1968), 199-213.

\bibitem{KatzICM} 
N. Katz. {\em The regularity theorem in algebraic geometry}, 
Actes du Congr\'es International des Math\'ematiciens (Nice, 1970), Tome \textbf{1}, pp. 437--443. 
Gauthier-Villars, Paris, 1971. 



\bibitem{KatzESDE}
N. Katz. 
{\em Exponential Sums and Differential Equations.}
{\sc Annals of Mathematics Studies} {\bf 124},
Princeton University Press (1990).

\bibitem{Katz}
N. Katz. 
{\em Rigid local systems.}
{\sc Annals of Mathematics Studies} {\bf 139},
Princeton University Press (1996).




\bibitem{Konno}
H. Konno. Construction of the moduli space of stable parabolic Higgs bundles on a Riemann surface. {\em J. Math.
Soc. Japan} {\bf 45} (1993), 253-276.



\bibitem{KorevaarSchoen}
N. Korevaar, R. Schoen. Sobolev spaces and harmonic maps for metric space targets. 
{\em Comm. Anal. Geom.} {\bf 1} (1993), 561-659.

\bibitem{KostovCRAS}
V. Kostov. On the Deligne-Simpson problem. {\em C.R.A.S.} {\bf 329} (1999), 657-662.

\bibitem{KostovKappa}
V. Kostov. 
The Deligne-Simpson problem for zero index of rigidity. 
{\em Perspectives of complex analysis, differential geometry and mathematical physics (St. Konstantin, 2000)},
World Sci. Publ., (2001), 1-35. 	

\bibitem{KostovSteklov}
V. Kostov. On the Deligne-Simpson problem. 
{\em Proc. Steklov Inst. Math.} {\bf 238} (2002), 148-185. 
 	
\bibitem{KostovControl}
V. Kostov.
On some aspects of the Deligne-Simpson problem. 
{\em J. Dynam. Control Systems} {\bf 9} (2003), 393-436. 

\bibitem{KostovAndWeak}
V. Kostov.
On the Deligne-Simpson problem and its weak version. 
{\em Bull. Sci. Math.} {\bf 128} (2004), 105-125. 
 	

\bibitem{KostovSurvey}
V. Kostov.
The Deligne-Simpson problem---a survey. 
{\em J. Algebra} {\bf 281} (2004), 83-108. 

\bibitem{KostovConnectedness}
V. Kostov.
The connectedness of some varieties and the Deligne-Simpson problem. 
{\em J. Dyn. Control Syst.} {\bf 11} (2005),  125-155. 


\bibitem{LaumonSurvey}
G. Laumon.
Exponential sums and $l$-adic cohomology: a survey. 
{\em Israel J. Math.} {\bf 120} (2000), 225-257. 


\bibitem{Laumon}
G. Laumon. 
La transformation de Fourier g\'eom\'etrique et ses applications. 
{\em Proceedings of ICM, Kyoto 1990}, Math. Soc. Japan (1991), 
437-445. 


\bibitem{Lawrence}
R. Lawrence. 
Homological representations of the Hecke algebra.
{\em Comm. Math. Phys.} {\bf 135} (1990), 141-191.

\bibitem{Li}
Jiayu Li. Hermitian-Einstein metrics and Chern number inequalities on parabolic stable bundles over K\"ahler 
manifolds.  {\em Comm. Anal. Geom.} {\bf 8} (2000), 445-475.

\bibitem{LiNarasimhan}
J. Li, M.S. Narasimhan. Hermitian-Einstein metrics on parabolic stable bundles. 
{\em Acta Math. Sin.} {\bf 15} (1999), 93-114.

\bibitem{LiNarasimhan2}
Jiayu Li, M. S. Narasimhan. A note on Hermitian-Einstein metrics on parabolic stable bundles. 
{\em Acta Math. Sinica} {\bf 17} (2001), 77-80.

\bibitem{Litcanu} 
R. Litcanu. 
Lam\'e operators with finite monodromy---a combinatorial approach. 
{\em J. Differential Equations} {\bf 207} (2004), 93-116. 
 	

\bibitem{LubotskyMagid}
A. Lubotsky, A. Magid. 
{\em Varieties of representations of finitely generated groups.}
{\sc Mem. Amer. Math. Soc.} {\bf 58} (1985).

	

\bibitem{Manin} 
Y. Manin. {\em  Moduli Fuchsiani}, Annali Scuola Normale Sup. di Pisa Ser. III
{\bf 19} (1965), 113-126.

\bibitem{Marin}
I. Marin.
Caract\'erisations de la repr\'esentation de Burau. 
{\em Expo. Math.} {\bf 21} (2003), 263-278. 

\bibitem{MarkmanXia}
E. Markman, E. Xia. The moduli of flat ${\rm PU}(p,p)$-structures with large Toledo invariants.  
{\em Math. Z.} {\bf 240} (2002), 95-109.

\bibitem{MaruyamaYokogawa} 
M. Maruyama, K. Yokogawa. 
{\em Moduli of parabolic stable sheaves},
Math. Ann. 293 (1992), no. \textbf{1}, 77--99.



\bibitem{MatsukiOlsson} 
K. Matsuki, M. Olsson. 
{\em  Kawamata-Viehweg 
vanishing and Kodaira vanishing for stacks}, Math. Res. Lett.
{\bf 12} (2005), 207-217.
 
\bibitem{Mehta}
M. Mehta. 
Birational equivalence of Higgs moduli. 
{\em Internat. J. Math.} {\bf 16} (2005), 365-386.


\bibitem{TMochizuki}
T. Mochizuki. Asymptotic behaviour of tame nilpotent harmonic bundles with trivial parabolic structure. 
{\em J. Diff. Geom.} {\bf 62} (2002), 351-559.


\bibitem{TMochizuki2}
T. Mochizuki.
Asymptotic behaviour of tame harmonic bundles and an application to pure twistor $D$-modules.
Preprint {\tt math.DG/0312230}. 


\bibitem{TMochizuki3} Mochizuki, T. 
Kobayashi-Hitchin correspondence for 
tame harmonic bundles and an application. 
Preprint {\tt math.DG/0411300}.

\bibitem{MoishezonTeicher}
B. Moishezon, M. Teicher. Braid group technique in complex geometry. I. Line arrangements in $\cc \pp ^2$.
{\em Braids (Santa Cruz, 1986)}, {\sc Contemp. Math.} {\bf 78}, A.M.S. (1988), 425-555.




\bibitem{Nakajima}
H. Nakajima. Hyper-K\"ahler structures on moduli spaces of parabolic Higgs bundles on Riemann surfaces. 
{\em Moduli of vector bundles (Sanda, 
Kyoto, 1994)} {\sc Lect. Notes Pure Appl. Math.} {\bf 179} (1996), 199-208.



\bibitem{NetoSilva}
O. Neto, P. Silva. 
On regular holonomic systems with solutions ramified along $y\sp k=x\sp n$. 
{\em Pacific J. Math.} {\bf 207} (2002), 463-487. 


\bibitem{Nitsure}
N. Nitsure. 
Moduli space of semistable pairs on a curve. {\em Proc. London Math. Soc.} {\bf 62} (1991), 275-300.


\bibitem{NitsureLog}
N. Nitsure. 
Moduli of semistable logarithmic connections. {\em J. Amer. Math. Soc.} {\bf 6} (1993), 597-609. 

\bibitem{Noohi}
B. Noohi. 
Fundamental groups of algebraic stacks. 
{\em J. Inst. Math. Jussieu} {\bf 3} (2004), 69-103.



\bibitem{Oh} Ohtsuki, M. 
{\em A residue formula for Chern classes associated with logarithmic connections}, 
Tokyo J. Math. 5 (1982), no. \textbf{1}, 13--21.  

\bibitem{Panov} D. Panov. Polyhedral K\"ahler manifolds. {\em Doctoral thesis}, 
Ecole Polytechnique (2005).



\bibitem{ReznikovContinuous}
A. Reznikov. Continuous cohomology of the group of volume-preserving and symplectic diffeomorphisms, measurable transfer 
and higher asymptotic cycles. {\em Selecta Math.} {\bf 5} (1999), 181-198. 

\bibitem{ReznikovH2}
A. Reznikov. The structure of K\"ahler groups. I. Second cohomology. {\em Motives, polylogarithms and 
Hodge theory, Irvine, 1998}, 
{\sc Int. Press Lect. Ser.} {\bf 3} (II) (2002), 717-730.

\bibitem{Roberts}
D. Roberts. Rigid Jordan tuples. Preprint available at \verb}http://cda.morris.umn.edu/~roberts/}.


\bibitem{SabbahIrregular}
C. Sabbah. Harmonic metrics and connections with irregular singularities. 
{\em Ann. Inst. Fourier} {\bf 49} (1999), 1265-1291.


\bibitem{SabbahTwistor}
C. Sabbah.
Polarizable twistor $\Dd$-modules. 
{\em Ast\'erisque} {\bf 300}, (2005).


\bibitem{Schmitt}
A. Schmitt. Projective moduli for Hitchin pairs. 
{\em Internat. J. Math.} {\bf 9} (1998), 107-118. 


\bibitem{Seshadri} 
Seshadri, C. S.
{\em Moduli of vector bundles on curves with parabolic structures.}
Bull. Amer. Math. Soc. \textbf{83} (1977), 124--126.

\bibitem{cvhs}
C. Simpson.
Constructing variations of Hodge structure using Yang-Mills theory and applications to uniformization.
{\em J. Amer. Math. Soc.} {\bf 1} (1988), 867-918.

\bibitem{hbnc}
C. Simpson. Harmonic bundles on noncompact curves. {\em J. Amer. Math. Soc.} {\bf 3}  (1990),  713-770.

\bibitem{products}
C. Simpson. 
Products of matrices. {\em Differential geometry, global analysis, and topology (Halifax, 1990)},
{\sc CMS Conf. Proc.} {\bf 12} (1991), 157-185. 

  
\bibitem{SteerWren}
B. Steer, A. Wren. The Donaldson-Hitchin-Kobayashi correspondence for 
parabolic bundles over orbifold surfaces.
{\em Canad. J. Math.} {\bf 53} (2001), 1309-1339.

\bibitem{StrambachVolklein}
K. Strambach, H. V\"{o}lklein. 
On linearly rigid tuples. 
{\em J. Reine Angew. Math.} {\bf 510} (1999), 57-62. 
 	
\bibitem{Szabo}
S. Szabo. 
Nahm transform for integrable connections on the Riemann sphere.
Preprint {\tt math.DG/0511471}. 


\bibitem{Thaddeus}
M. Thaddeus. Variation of moduli of parabolic Higgs bundles. {\em J. Reine Angew. Math.} {\bf 547} (2002), 1-14. 

\bibitem{TubaWenzel}
I.  Tuba, H. Wenzl. 
Representations of the braid group $B\sb 3$ and of ${\rm SL}(2, Z)$. 
{\em Pacific J. Math.} {\bf 197} (2001), 491-510. 

\bibitem{vanderPutBourbaki}
M.  van der Put. 
Recent work on differential Galois theory. 
{\em S\'eminaire Bourbaki 1997/98}. {\sc Ast\'erisque} {\bf 252} (1998), Exp. No. 849, 341-367. 

	
\bibitem{vanderPut}
M. van der Put.
Grothendieck's conjecture for the Risch equation $y'=ay+b$. 
{\em Indag. Math. (N.S.)} {\bf 12} (2001), 113-124. 

 	
\bibitem{VolkleinTransformation}
H. V\"{o}lklein.
A transformation principle for covers of $\pp ^1$. 
{\em J. Reine Angew. Math.} {\bf 534} (2001), 156-168. 

 	
\bibitem{VolkleinBraid}
H. V\"{o}lklein. 
The braid group and linear rigidity. 
{\em Geom. Dedicata} {\bf 84} (2001), 135-150. 
 	
\bibitem{Xia}
E. Xia. The moduli of flat $U(p,1)$ structures on Riemann surfaces. {\em Geom. Dedicata} {\bf 97} (2003), 33-43. 

\bibitem{Yokogawa}
K. Yokogawa. Compactification of moduli of parabolic sheaves and moduli of parabolic Higgs sheaves.
{\em J. Math. Kyoto Univ.} {\bf 33} (1993), 451-504.

\bibitem{Yokogawa2}
K. Yokogawa.
Infinitesimal deformation of parabolic Higgs sheaves. {\em Internat. J. Math.} {\bf 6} (1995), 125-148. 

\bibitem{Yokoyama}
T. Yokoyama.
Construction of systems of differential equations of Okubo normal form with rigid monodromy. 
{\em Math. Nachr.} {\bf 279} (2006), 327-348. 




\end{thebibliography}
\end{document}